\documentclass[11pt,a4paper]{article}

\usepackage[utf8]{inputenc}
\usepackage[OT2,T1]{fontenc}

\usepackage[a4paper]{geometry}
\usepackage[french,english]{babel}

\usepackage{fancyhdr}
\usepackage[final]{graphicx}
\usepackage{numprint} 

\usepackage[figurewithin=none]{caption}
\usepackage{zref}
\usepackage{dsfont}
\usepackage{enumitem}
\usepackage{ifthen}
\usepackage{setspace}
\usepackage{amsmath,amssymb,makeidx,amsthm}
\usepackage{stmaryrd}
\usepackage{mathrsfs}
\usepackage[all]{xy} 
\usepackage{multirow} 
\usepackage{array}
\usepackage{mathtools}
\usepackage{tikz-cd}
\usepackage{comment}

\usepackage{url}
\usepackage[hypertexnames=false]{hyperref}

\newcommand{\R}{{\mathbb R}}

\newcommand{\Z}{{\mathbb Z}}
\newcommand{\Q}{{\mathbb Q}}
\renewcommand{\C}{{\mathbb C}}

\newcommand{\F}{{\mathbb F}}
\DeclareSymbolFont{cyrletters}{OT2}{wncyr}{m}{n}
\DeclareMathSymbol{\Sha}{\mathalpha}{cyrletters}{"58}
\newcommand{\rar}{\rightarrow}

\newcommand{\Sch}{\ensuremath{\mathbf{Sch}}}

\newcommand{\Ab}{\ensuremath{\mathbf{Ab}}}

\newcommand{\mfk}[1]{\mathfrak{#1}}
\newcommand{\Nekovar}{Nekov\'a\v r~}

\renewcommand{\emptyset}{\varnothing}
\newcommand{\mrm}[1]{\mathrm{#1}}

\newcommand{\OO}{\mathcal{O}}
\newcommand{\Qbar}{\overline{\mathbb{Q}}}

\newcommand{\Fr}{\mrm{Frob}}
\newcommand{\Tate}[2]{\ensuremath{\mathcal{T}_{#1}}#2}

\renewcommand{\tilde}[1]{\widetilde{#1}}
\newcommand{\Aut}[1]{\ensuremath{\operatorname{Aut}(#1)}}

\newcommand{\GL}[1]{\ensuremath{\operatorname{GL}_2(#1)}}
\newcommand{\SL}[1]{\ensuremath{\operatorname{SL}_2(#1)}}
\newcommand{\PGL}[1]{\ensuremath{\operatorname{PGL}_2(#1)}}

\newcommand{\diam}[1]{\ensuremath{\langle #1\rangle}}

\DeclareMathOperator{\Sp}{Spec}

\usepackage{calc}

\makeatletter
\let\oldlabel\label
\renewcommand{\label}[1]{%
  \zref@labelbylist{#1}{special}
  \oldlabel{#1}
}

\zref@newlist{special}
\zref@newprop{section}{\arabic{section}.\arabic{subsection}}
\zref@addprop{special}{section}
\newcounter{propo}
\renewcommand{\thepropo}{\thesubsection.\arabic{propo}}
\makeatletter \@addtoreset{propo}{subsection}\makeatother

\makeatletter \@addtoreset{rema}{subsection}\makeatother

\newcounter{mytheo}
\renewcommand{\themytheo}{\Alph{mytheo}}

\newcounter{contentsin}
\renewcommand{\thecontentsin}{\arabic{contentsin}}

\newcommand{\lem}[2][None]{\refstepcounter{propo}\ifthenelse{\equal{#1}{None}}{}{\label{#1}} \smallskip \noindent \textbf{Lemma \thepropo\;} {\it #2}\smallskip}

\newcommand{\theo}[2][None]{\refstepcounter{propo}\ifthenelse{\equal{#1}{None}}{}{\label{#1}} \smallskip \noindent \textbf{Theorem \thepropo\;} {\it #2} \smallskip}
\newcommand{\theoi}[2][None]{\refstepcounter{theop}\ifthenelse{\equal{#1}{None}}{}{\label{#1}}\smallskip \noindent \textbf{Theorem \thetheop\;} {\it #2} \smallskip}
\newcommand{\cori}[2][None]{\refstepcounter{theop}\ifthenelse{\equal{#1}{None}}{}{\label{#1}}\smallskip \noindent \textbf{Corollary \thetheop\;} {\it #2} \smallskip}

\newcommand{\prop}[2][None]{\refstepcounter{propo}\ifthenelse{\equal{#1}{None}}{}{\label{#1}} \smallskip \noindent \textbf{Proposition \thepropo\;} {\it #2} \smallskip}

\newcommand{\cor}[2][None]{\refstepcounter{propo}\ifthenelse{\equal{#1}{None}}{}{\label{#1}} \smallskip \noindent \textbf{Corollary \thepropo\;} {\it #2} \smallskip}

\newcommand{\defi}[2][None]{\refstepcounter{propo}\ifthenelse{\equal{#1}{None}}{}{\label{#1}} \smallskip \noindent {\bf Definition \thepropo\;} {#2} \smallskip}

\newcommand{\rem}[2][None]{\refstepcounter{propo}\ifthenelse{\equal{#1}{None}}{}{\label{#1}}\smallskip \noindent {\bf Remark \thepropo\;} {#2} \smallskip}
\newcommand{\rems}[1]{\smallskip \noindent {\bf Remarks\;} {#1} \smallskip}

\newcommand{\demo}[1]{\noindent {\it Proof.\;--\;} #1\hfill$\Box$ \smallskip}

\newcommand{\qnn}[2][None]{\refstepcounter{propo}\ifthenelse{\equal{#1}{None}}{}{\label{#1}} \smallskip \noindent \textbf{Question \thepropo\;} {\it #2} \smallskip}

\newcommand{\mytheo}[2][None]{\refstepcounter{mytheo}\ifthenelse{\equal{#1}{None}}{}{\label{#1}} \smallskip \noindent \textbf{Theorem \themytheo\;} {\it #2}\smallskip}
\newcommand{\mycor}[2][None]{\refstepcounter{mycor}\ifthenelse{\equal{#1}{None}}{}{\label{#1}} \smallskip \noindent \textbf{Corollary \themycor\;} {\it #2}\smallskip}
\newcommand{\conjin}[2][None]{\refstepcounter{contentsin}\ifthenelse{\equal{#1}{None}}{}{\label{#1}} \smallskip \noindent \textbf{Conjecture \thecontentsin\;} {\it #2}\smallskip}
\newcommand{\questin}[2][None]{\refstepcounter{contentsin}\ifthenelse{\equal{#1}{None}}{}{\label{#1}} \smallskip \noindent \textbf{Question \thecontentsin\;} {\it #2}\smallskip}
\newcommand{\theoin}[2][None]{\refstepcounter{contentsin}\ifthenelse{\equal{#1}{None}}{}{\label{#1}} \smallskip \noindent \textbf{Theorem \thecontentsin\;} {\it #2}\smallskip}
\newcommand{\defin}[2][None]{\refstepcounter{contentsin}\ifthenelse{\equal{#1}{None}}{}{\label{#1}} \smallskip \noindent \textbf{Definition \thecontentsin\;} {\it #2}\smallskip}

\title{Congruences modulo $23$ to $y^2=x^3-23$ are trivial}

\author{Elie Studnia}

\begin{document}

\maketitle

\section*{Introduction}

The primary goal of this paper is, as its title suggests, to prove the following result: 

\mytheo[mytheo-23]{Let $E/\Q$ be the elliptic curve with Weierstrass equation $y^2=x^3-23$. Let $F/\Q$ be an elliptic curve such that $F[23](\Qbar)$ and $E[23](\Qbar)$ are isomorphic $\mrm{Gal}(\Qbar/\Q)$-modules. Then $E$ and $F$ are isogenous. }

Theorem \ref{mytheo-23} lies at the crossroads of two important open questions regarding the arithmetic of elliptic curves: the \emph{Frey--Mazur conjecture} and \emph{Serre's uniformity question}. 

Let $N \geq 1$ be an integer and $E,F$ be two elliptic curves over a field $K$ of characteristic not dividing $N$. We say that $E$ and $F$ are \emph{congruent modulo $N$} if the $G_K := \mrm{Gal}(\overline{K}/K)$-modules $E[N](\overline{K})$ and $F[N](\overline{K})$ are isomorphic. As an example, $E$ and $F$ are always congruent modulo $N$ if $K$ is separably closed, since both modules are isomorphic to $(\Z/N\Z)^{\oplus 2}$ as abelian groups with a trivial Galois action. Another example (over an arbitrary base field $K$) is when there is an isogeny $\varphi: E \rar F$ with degree prime to $N$: we say in this case that the congruence is \emph{trivial}. Theorem \ref{mytheo-23} can be rephrased as follows: any congruence over $\Q$ to the elliptic curve with Weierstrass equation $y^2=x^3-23$ modulo $23$ is trivial. 

The Frey--Mazur conjecture originates from the following question, asked by Mazur in \cite{FreyMazur}: 

\questin[mazur-question]{Does there exist an integer $N \geq 7$ and a non-trivial congruence of elliptic curves over $\Q$ modulo $N$?}

In \cite{KO}, Kraus and Oesterl\'e give an affirmative answer for $N=7$, and provide concrete examples of such congruences. Since then, more congruences have been found; by work of Kani--Rizzo \cite{KR}, Halberstadt--Kraus \cite[Proposition 6.3]{HalbKraus}, Fisher \cite{Fisher-711, Fisher-13}, we know that there are infinite families of such non-trivial congruences for $N \in \{7,11,13\}$.\footnote{The papers cited add another constraint, which will appear later in this introduction: the isomorphism of Galois modules has a prescribed behavior with respect to the Weil pairings on both sides.} This question has also been studied when $N$ is composite or for $N < 7$, see for instance \cite{Fisher-Expl,Frengley-12,Frengley-Large} and references therein. 

However, it follows from a result of Bakker and Tsimerman \cite{BaTs} that there are no infinite families of non-trivial congruences modulo large enough primes $N$. In fact, congruences modulo larger prime values of $N$ seem exceedingly rare: we know of two pairs of non-trivial congruences modulo $17$ found by Cremona and Fisher (cf e.g. \cite{Fisher-17}). By a result of Cremona and Freitas \cite[Theorem 1.3]{Cremona-Freitas}, there are in fact \emph{no} non-trivial congruences modulo primes $N \geq 19$ between elliptic curves of conductor at most $500000$. 

The Frey--Mazur conjecture makes this observation more precise. It is a variant of conjectures of Frey \cite[Conjecture 5]{Frey-Conj} and seems to have first been written by Darmon \cite[Conjecture 4.3]{Darmon-Serre}. 

\conjin[curvewise-frey-mazur]{(``Curve-wise'' Frey--Mazur conjecture) Let $E/\Q$ be an elliptic curve. There is a constant $C_E > 0$ such that for every prime $p > C_E$, any congruence modulo $p$ to the elliptic curve $E$ over $\Q$ is trivial.}

Given that we do not know of any non-trivial congruence of elliptic curves over $\Q$ modulo any prime $p \geq 19$, the following much stronger version is often found in the literature (cf. e.g. \cite[Conjecture 1.1]{Fisher-17} or the discussion following \cite[Theorem 1.3]{Cremona-Freitas}):

\conjin[optimistic-frey-mazur]{(``Optimistic'' Frey--Mazur conjecture) Conjecture \ref{curvewise-frey-mazur} holds with $C_E=18$. }

The ``curve-wise'' Frey--Mazur conjecture is already a quite delicate problem, and it does not seem known for any elliptic curve $E/\Q$. In fact, we are not aware of any pair $(E/\Q, p)$ such that $p > 13$ is prime and all the elliptic curves over $\Q$ congruent to $E$ modulo $p$ are provably determined, unless the map $\rho_{E,p}: G_{\Q} \rar \Aut{E[p](\Qbar)}$ is not surjective, and the set of $j$-invariants of elliptic curves $F/\Q$ such that the image of $\rho_{F,p}$ is conjugate to that of $\rho_{E,p}$ has been found to be finite and explicitly determined.\footnote{While this is somewhat anachronistic, this restriction covers the results of \cite{HalbKraus-Torsion}, because of subsequent progress on the uniformity question due to Bilu--Parent \cite{BP} and Bilu--Parent--Rebolledo \cite{BPR}, which we will discuss below. One can hope that future progress on the uniformity question will supersede Theorem \ref{mytheo-23} and the other results of this paper.} 

This brings us to the second inspiration for this work, which is Serre's uniformity question. Recall Serre's celebrated open image theorem \cite{Serre-image}:

\theoin[serre-open-image]{Let $E/\Q$ be a non-CM elliptic curve. There exists a constant $B_E > 0$ such that for every prime $p > B_E$, the map $\rho_{E,p}: G_{\Q} \rar \Aut{E[p](\Qbar)}$ is surjective. }

Serre's uniformity question is whether the above constant $B_E$ can be made independent of (``uniform'' with respect to) $E$:

\questin[serre-uniformity]{Is it true that Theorem \ref{serre-open-image} holds for a constant $B_E=B$ that does not depend on $E$?}

Extensive numerical evidence and a considerable amount of work suggest that one should be able to answer in the affirmative to Question \ref{serre-uniformity} with $B=37$, and, in fact that the following more precise Conjecture holds. 

\conjin[optimistic-serre-uniformity]{(``Optimistic'' Serre uniformity) Let $E/\Q$ be an elliptic curve and $p \geq 17$ be a prime such that $\rho_{E,p}: G_{\Q} \rar \Aut{E[p](\Qbar)}$ is not surjective. Then one of the following holds:
\begin{itemize}[noitemsep,label=\tiny$\bullet$]
\item $E$ has CM,
\item $p=17$ and $j(E) \in \{-\frac{17^2\cdot 101^3}{2},-\frac{17\cdot 373^3}{2^{17}}\}$,
\item $p=37$ and $j(E) \in \{-7 \cdot 11^3, -7 \cdot 137^3 \cdot 2083^3\}$.
\end{itemize}}

What we currently know about Serre's uniformity question can be summarized in the following result, due to works of Serre \cite{Serre-image}, Mazur \cite{FreyMazur}, Bilu--Parent \cite{BP}, Bilu--Parent--Rebolledo \cite{BPR}, Balakrishnan--Dogra--M\"uller--Tuitman--Vonk \cite{BDMTV21}, Le Fourn--Lemos \cite{LFL}, Furio--Lombardo \cite{FL}.

\theoin[serre-uniformity-progress]{
Let $E/\Q$ be an elliptic curve and $p \geq 17$ be a prime. If $\rho_{E,p}: G_{\Q} \rar \Aut{E[p](\Qbar)}$ is not surjective, then one of the following holds:
\begin{itemize}[noitemsep,label=\tiny$\bullet$]
\item $E$ has CM,
\item $p=17$ and $j(E) \in \{-\frac{17^2\cdot 101^3}{2},-\frac{17\cdot 373^3}{2^{17}}\}$,
\item $p=37$ and $j(E) \in \{-7 \cdot 11^3, -7 \cdot 137^3 \cdot 2083^3\}$,
\item $p \geq 19$ and the image of $\rho_{E,p}$ is equal to the normalizer of a non-split Cartan subgroup. 
\end{itemize}}

Conjecture \ref{optimistic-serre-uniformity} predicts that the last case is contained in the first one, but how to prove it remains largely an open question. The main recent progress on that topic comes from the \emph{quadratic Chabauty method}, as demonstrated in \cite{BDMTV19,BDMTV21} in the cases where $p \in \{13,17\}$. It was shown in \cite{LFD} that the method should work in principle for general $p$, although the main result of \emph{loc.cit.} is weaker than what would be needed to prove Conjecture \ref{optimistic-serre-uniformity}. 

The link between Theorem \ref{mytheo-23} and Serre's uniformity question is the following: the elliptic curve $E_{23}$ with Weierstrass equation $y^2=x^3-23$ has complex multiplication by imaginary quadratic field $\Q(e^{2i\pi/3})$, which is inert at $23$. Therefore, if $F/\Q$ is any elliptic curve congruent to $E_{23}$, the image of $\rho_{F,23}$ is conjugate to that of $\rho_{E_{23},23}$, which is equal to the normalizer of a non-split Cartan subgroup. Assume Conjecture \ref{optimistic-serre-uniformity} (even only in the case $p=23$). Then $F$ has complex multiplication, hence potentially good reduction everywhere. By elementary local arguments (see Section \ref{subsect:local-standard}), $F$ has the same conductor as $E_{23}$ (which is $108 \cdot 23^2 \leq 500000$). Theorem \ref{mytheo-23} then follows from the result of Cremona--Freitas \cite[Theorem 1.3]{Cremona-Freitas}, or from a direct search of elliptic curves of conductor $108 \cdot 23^2$ in the LMFDB. 

Since Conjecture \ref{optimistic-serre-uniformity} is open even for $p=23$, this is not a valid proof of Theorem \ref{mytheo-23}, and we need to proceed in a different way. Our argument follows the strategy introduced by Mazur in \cite{FreyMazur}, and used to great effect ever since in the resolution of various questions related to the arithmetic of elliptic curves: it is the strategy employed in all the works leading to the proof of Theorem \ref{serre-uniformity-progress} (except for \cite{BDMTV21}), and is also the strategy underlying the known results on torsion of elliptic curves over number fields (such as \cite{Merel,Parent-Torsion}). 

To explain this strategy, we start from a general elliptic curve $E/\Q$ with conductor $N$ and pick a prime $p \geq 7$. As shown in \cite{Studnia-moduli}, there are smooth projective $\Z[1/Np]$-schemes $X_E^{\alpha}(p)$ (indexed by $\alpha \in \F_p^{\times}$ corresponding to the action of the isomorphism of Galois-modules on the Weil pairing -- which will not affect the argument) with geometrically connected fibres; the curves $X_E^{\alpha}(p)$ have a closed subscheme of \emph{cusps} $Z_E^{\alpha}(p)$, which is finite \'etale over $\Z[1/Np]$, such that the rational points of $X_E^{\alpha}(p)_{\Q} \backslash Z_E^{\alpha}(p)_{\Q}$ correspond to couples $(F,\iota)$, where $F/\Q$ is an elliptic curve and $\iota: E[p](\Qbar) \rar F[p](\Qbar)$ is an isomorphism of $G_{\Q}$-modules such that $\langle \iota(P),\,\iota(Q)\rangle_{F[p]} = \langle P,\,Q\rangle_{E[p]}^{\alpha}$ for $P,Q \in E[p](\Qbar)$. 

Mazur's strategy consists in finding a morphism $f: X_E^{\alpha}(p) \rar A$ to a quotient of the Jacobian of $X_E^{\alpha}(p)$ satisfying the following two properties:

\begin{itemize}[noitemsep,label=\tiny$\bullet$]
\item $f$ is a \emph{formal immersion} at the cusps in every prime characteristic $\ell \nmid Np$,
\item $A$ has finitely many rational points. 
\end{itemize} 

Mazur's \emph{formal immersion argument} (originally exposed as \cite[Corollary 4.3]{FreyMazur}) implies that, in these conditions, every point in $(X_E^{\alpha}(p) \backslash Z_E^{\alpha}(p))(\Q)$ lifts to a point of $(X_E^{\alpha}(p) \backslash Z_E^{\alpha}(p))(\Z[1/2Np])$. Translating this back into the original problem, this means that, for a pair $(F,\iota)$ as above, the primes dividing the denominator of $j(F)$ divide $2Np$. By local arguments, this implies that $F$ has good reduction away from $2Np$. 

In the context of Theorem \ref{mytheo-23}, one has $(E,N,p)=(E_{23},108\cdot 23^2,23)$, and we can use the local arguments of Section \ref{subsect:local-standard} to show that $F$ has conductor $N$. As above, \cite[Theorem 1.3]{Cremona-Freitas} implies that $E_{23}$ and $F$ are isogenous. 

While this is much simpler than the arguments of \cite{FreyMazur} --- where the classification of cyclic isogenies of prime degree of rational elliptic curves is reduced to the class number one problem --- or of \cite{BP,BPR,LFL,FL} --- the authors find a contradiction between upper bounds on the $j$-invariants obtained using the Runge method (an example is \cite[Proposition 6.1]{LFL}) with lower bounds found by transcendence methods for non-CM rational points (an example being \cite[Theorem 4.1]{FL}) --- this reliance on what is ultimately an exhaustive search for elliptic curves of conductor $108p^2$ makes this step unsuitable for applications to large primes $p$, which are perhaps more interesting for the study of Conjectures \ref{curvewise-frey-mazur}, \ref{optimistic-frey-mazur}. Still, it is not clear to us what alternate arguments could be used (see the end of this introduction for a brief discussion). 

It remains to construct the morphism $f$ and the abelian variety $A$ satisfying the above properties (if only in the case $(E,p)=(E_{23},23)$), and this is the main focus of this paper. 

Let $E/\Q$ be an elliptic curve and $p \geq 7$ be a prime. In order to determine whether the Jacobian $J_E^{\alpha}(p)$ of $X_E^{\alpha}(p)$ has a quotient with finitely many rational points, one possibility is to decompose its Tate module and find out whether the $q$-divisible Selmer group (attached to some quotient) is finite. Understanding the Tate module of the Jacobian of a general curve is a difficult problem, but the curve $X_E^{\alpha}(p)$ is quite special: while it is not \emph{per se} a modular curve as usually defined, it is isomorphic over $\Qbar$ to the (connected) classical modular curve $X(p)^{\text{conn.}}_{\Qbar}$. Moreover, $X_E^{\alpha}(p)$ is a (compactified) moduli space for certain elliptic curves, and, as shown in \cite[Corollary 1]{Studnia-moduli}, it admits \emph{Hecke correspondences} that define \emph{Hecke operators} $T_{\ell}: J_E^{\alpha}(p) \rar J_E^{\ell\alpha}(p)$ for $\ell \neq p$ prime. Note that $T_{\ell}$ does not define an endomorphism of $J_E^{\alpha}(p)$ but rather ``moves between''  the $J_E^{\alpha}(p)$ --- which causes some technical difficulties to be circumvented.     

This information is (almost) enough to determine the Tate module of $J_E^{\alpha}(p)$ (which is Corollary \ref{decomp-tate-twisted-gq}). Unfortunately, despite the apparent simplicity of the description present in this Corollary, it seems very difficult to extract useful arithmetic information from these representations. Indeed, most of the tools needed to do so require that we know some \emph{automorphic} information on the Galois representation. Since the representations appearing in \emph{loc.cit.} are tensor products of modular forms of weight two with large-dimensional Artin representations that do not factor through a solvable extension of a totally real field (at least in the ``generic'' case where $\rho_{E,p}$ is surjective), understanding their automorphic properties seems currently beyond reach \footnote{With the exception of one rather intriguing mod $p$ automorphy result \cite[Theorem 4.1.6]{Studnia-thesis}, although its significance remains unclear.}.   

Hence, we assume that $\rho_{E,p}$ is not surjective: by Theorem \ref{serre-uniformity-progress}, the most interesting case is the one where the image of $\rho_{E,p}$ is equal to the normalizer $N_C$ of a non-split Cartan subgroup $C$, which we now assume. We also assume that $p \geq 11$: this implies that the number field $K$ with absolute Galois group $\rho_{E,p}^{-1}(C)$ is imaginary quadratic and inert at $p$. The Tate module of $J_E^{\alpha}(p)$ then admits a refined decomposition given in Corollary \ref{tate-module-twisted-cartan-gq}. This situation is far more promising, since these representations are self-dual up to Tate twist, and they are known to possess entire $L$-functions verifying a functional equation, and, by the work of Kato \cite{KatoBSD} and Kings--Loeffler--Zerbes \cite{KLZ15} we can (under certain assumptions) attach to them an Euler system satisfying a reciprocity law. Combining these results with large Galois image results (mainly \cite{bigimage} and its improvements in \cite[Theorem A]{Studnia-Euler}), we deduce (as Proposition \ref{BSD-that-interests-us}) that many quotients of $J_E^{\alpha}(p)$ satisfy the Birch and Swinnerton--Dyer conjecture in analytic rank zero. These quotients are usually not of $\mrm{GL}_2$-type, and this instance of the Birch and Swinnerton--Dyer conjecture (in analytic rank zero) was previously not known.  

The issue is that the sign of the functional equation for most of the Galois representations appearing in the Tate module of $J_E^{\alpha}(p)$ is $-1$, so that the $L$-function has a forced root at the center of the functional equation, and --- assuming the Birch and Swinnerton--Dyer conjecture --- the corresponding quotient of $J_E^{\alpha}(p)$ has infinitely many rational points. \footnote{While it is somewhat disappointing, this result is perhaps not so surprising, given that the ``normalizer of the non-split Cartan'' case of the uniformity question has not been solved so far.} We do not carry out this computation here; we refer the interested reader to \cite[Sections 5.3, 5.4]{Studnia-thesis}. Suffice to say for now that, among others, when $E/\Q$ has good reduction at $p$, for \emph{every} quotient $A$ of $J_E^{\alpha}(p)$, the $L$-function of $A$ vanishes at the center of the functional equation, so we expect $A(\Q)$ to be infinite.  \footnote{While the result of the sign calculation seems more difficult to interpret in the generic case (see \cite[Table 4.1]{Studnia-thesis}), the lack of quotients of rank zero is easier to explain: when $\rho_{E,p}$, it follows from Corollary \ref{decomp-tate-twisted-gq} that the isogeny quotients of $J_E^{\alpha}(p)$ are cut out a finite set $T$ of Hecke operators independent from $E$. Hence, if $X_E^1(p)$ has a quotient of rank zero, the classical connected modular curve $X(p)^{\text{conn.}}_{\Qbar}$ has a $\Qbar$-point $x$ whose image in the Jacobian $J(p)^{\text{conn.}}$ is killed by some $t \in T$. The Mordell--Lang conjecture (cf. e.g. \cite[Th\'eor\`eme 1.1]{Remond}) implies that only finitely many such $\Qbar$-points exist. Thus, for every $p > 13$, the set of $j(E)$ such that $X_E^1(p)$ has a quotient of rank zero is finite.}

The simpler of the two situations in which $J_E^{\alpha}(p)$ has a quotient $A$ in which no factor of $L(A,s)$ has a forced vanishing is if the image of an inertia group at $p$ by $\rho_{E,p}$ is \emph{smaller than expected}, i.e. it is \emph{properly contained} in $C$, which implies among others that $p \equiv 2, 5 \mod{9}$. One can view this as a local variant of the condition studied in \cite{LFL,FL}. We make this assumption from now on for the rest of this introduction; by a result of Kraus \cite{Kraus-thesis}, this condition can be checked using only the minimal discriminant of $E$. In particular, if $p \equiv 5 \mod{9}$ (resp. $p \equiv 2 \mod{9}$), the entire discussion applies to the elliptic curve with Weierstrass equation $y^2=x^3+p$ (resp. $y^2=x^3+p^2$). 

We show: 

\mytheo[existence-of-quotient-bk]{Assume that the image by $\rho_{E,p}$ of an inertia subgroup at $p$ is properly contained in $C$. Let $\omega \subset \mathcal{S}_2(\Gamma_0(p))$ be a Galois orbit of newforms of size $r \geq 1$ and $\psi: C \rar \C^{\times}$ be a character of order exactly $3$. The Artin representation $\mrm{Ind}_C^N{\psi} \circ \rho_{E,p}$ is attached to a newform $g \in \mathcal{S}_1(\Gamma_1(M))$ which has rational coefficients and CM by $K$. 

There is a quotient $A$ of $J_E^{\alpha}(p)$ with dimension $2r$ which has finitely many rational points if the special Rankin--Selberg $L$-value $L(f,g,1)$ is non-zero for some (or equivalently, any) $f \in \omega$. }

In the situation that we describe, $M$ is prime to $p$, and the ``imprimitive'' Rankin--Selberg $L$-function satisfies the expected functional equation, so the non-vanishing of the $L$-value $L(f,g,1)$ can be verified numerically when the product $Mp$ is not too large.\footnote{We are not aware of an available implementation even in this case, so we wrote \emph{ad hoc} PARI-GP code with the help of P. Molin.} Using an explicit Waldspurger formula due to Cai--Shu--Tian \cite{CaiShuTian}, we can give a criterion for the non-vanishing of the $L$-value that relies only on discrete calculations. 

\mytheo[nonvan-l-function]{Let $\psi: C \rar \C^{\times}$ be a character of order exactly $3$. The Artin representation $\mrm{Ind}_C^N{\psi} \circ \rho_{E,p}$ is attached to a newform $g \in \mathcal{S}_1(\Gamma_1(M))$. Assume that the Hilbert class polynomial of discriminant $-M$ is not a cube modulo $p$, which is the case if $p > M$. Then there exists $f \in \mathcal{S}_2(\Gamma_0(p))$ such that $L(f,g,1) \neq 0$.   }

To complete Mazur's strategy, it remains to determine whether some map $f: X_E^{\alpha}(p) \rar A$ (where $A$ is the quotient appearing in Theorem \ref{existence-of-quotient-bk}) is a formal immersion in finite characteristic at the cusps. We use a slightly weaker condition than the ``formal immersion'' for technical reasons, which turns out to be a harmless simplification when attempting to apply Mazur's strategy --- we end up replacing the ``geometric'' formal immersion argument with an ``analytic'' version involving a simple case of the Chabauty--Coleman method. 

What makes this delicate is that $A/\Q$ is not the largest quotient of $J_E^{\alpha}(p)$ on which the Hecke algebra acts in a certain way. We also need to keep track of some of the $\GL{\F_p}$-action on $\coprod_{\alpha \in \F_p^{\times}}{X_E^{\alpha}(p)_{\Qbar}}$ (this is done by choosing a basis of $E[p](\Qbar)$, which lets us identify this disjoint reunion to the disconnected modular curve $X(p)^{\text{not conn.}}_{\Qbar}$). The quotient $A$ is the largest quotient of $J_E^{\alpha}(p)$ on which not only the Hecke operators act in a prescribed way, but so do certain of their products with elements of the normalizer of a suitable Cartan subgroup of $\GL{\F_p}$ (see Situation \ref{quotient-by-3C}). As a consequence, the usual $q$-expansion argument found in e.g. \cite[Proposition 4.1]{LFL} cannot be applied as is. 

\defin[exceptional-primes-intro]{Let $p \equiv 2,5 \mod{9}$ be a prime. For every coset $j$ of $\F_{p^2}^{\times}$ modulo its subgroup $\F_{p^2}^{\times 3}$ of cubes, we define \[u_j := \sum_{b \in \F_p, 1+b\sqrt{-3} \in j}{e^{\frac{2i\pi}{p}b}} \in \Z[e^{\frac{2i\pi}{p}}].\]
Let $I_p \subset \Z[e^{\frac{2i\pi}{p}}]$ be the ideal generated by the three $u_j$ for $j \in \F_{p^2}^{\times}/\F_{p^2}^{\times 3}$. We say that a prime $\ell$ is \emph{exceptional} for $p$ if $\ell$ divides the absolute norm of $I_p$.}

It turns out that the map $X_E^{\alpha}(p) \rar A$ behaves well at cusps when the characteristic is non-exceptional for $p$. Applying our analytic version of the formal immersion argument, we find: 

\mytheo[formal-immersion-applied]{Let us keep the notations of Theorem \ref{existence-of-quotient-bk} and suppose that $A$ has finitely many rational points. Then, for any elliptic curve $F/\Q$ congruent to $E$ modulo $p$, the primes $\ell$ dividing the denominator of $j(F)$ are exceptional for $p$ and congruent to $\pm 1 \mod{p}$.}

Of course, Faltings's proof of the Mordell conjecture (see e.g. \cite{Faltings}) implies that, since $X_E^{\alpha}(p)$ has genus at least $2$, there is \emph{some} set of ``exceptional primes'' such that for every elliptic curve $F/\Q$ congruent to $E$ modulo $p$, the denominator of $j(F)$ is only divisible by the exceptional primes. What makes Theorem \ref{formal-immersion-applied} interesting is the very explicit definition of this set. In fact, it is not difficult to check that when $11 \leq p \leq 10000$, there are \emph{no} exceptional primes for $p$ --- but going beyond this bound seems computationally impractical. 

At any rate, one concrete example of application of the above results for the elliptic curves $E_p$ previously defined is the following generalization of Theorem \ref{mytheo-23}:

\mytheo[applications-to-ep]{Let $p \geq 11$ be a prime congruent to $2$ modulo $9$ (resp. $5$ modulo $9$), and $E_{p}$ be a quadratic twist of the elliptic curve with Weierstrass equation $y^2=x^3+p^2$ (resp. $y^2=x^3+p$). Let $F/\Q$ be an elliptic curve congruent to $E$ modulo $p$. Then the prime numbers dividing the denominator of $j(F)$ are congruent to $\pm 1$ modulo $3p$ and exceptional for $p$. In particular, if $p \leq 10000$, then $j(F) \in \Z$ and $F$ has the same conductor as $E$. If $p \leq 59$, $F$ and $E$ are isogenous.}

Before concluding this introduction, let us briefly some questions left open by this work. 

\begin{enumerate}[itemsep=6pt]
\item \emph{Exceptional primes:} The strength of Theorem \ref{formal-immersion-applied} for a given prime $p$ is dependent on the determination of exceptional primes for $p$. As written above, we were able to show that there were no exceptional primes for $p \leq 10000$, but our computational approach does not seem to scale beyond this range, and we do not know if some prime $p \geq 10000$ congruent to $2,5$ modulo $9$ admits exceptional primes. A heuristic due to discussions with Hendrik Lenstra and Pieter Moree (see Section \ref{subsect:exceptional-primes-steinberg}), using results of Fan, Languasco, Lunia and Moree \cite{Moree} suggests that this should only occur for finitely many $p$. We tentatively conjecture that no prime $p \geq 11$ congruent to $2,5$ modulo $9$ admits any exceptional prime. 
\item \emph{Checking the non-vanishing of the $L$-value:} The assumptions of in Theorem \ref{nonvan-l-function} are much stronger than what seems empirically needed. For instance, when $-3 \geq D \geq -10000$ is a discriminant and $p \geq 3\sqrt{|D|}$ is a prime congruent to $2,5$ modulo $9$, inert in $\Q(\sqrt{D})$ and prime to $D$, the Hilbert class field polynomial of discriminant $D$ is not a cube modulo $p$ (see Section \ref{subsect:lvalue-quaternion-supersing}). Can we improve on Theorem \ref{nonvan-l-function}, for instance in the important case (cf. Proposition \ref{result-for-en}) where $\Q(\sqrt{D})=\Q(\sqrt{-3})$? 
\item \emph{The integral case:} Let $p \equiv 5\mod{9}$ be a prime and $p^{\ast}=\pm p$ be congruent to $1$ modulo $4$. Let $E$ be the elliptic curve $y^2=x^3+p^{\ast}$ (with conductor $108p^2$, the smallest one among its quadratic twists). Suppose that there are no exceptional primes for $p$. Let $F/\Q$ be an elliptic curve congruent to $E$ modulo $p$. By Theorem \ref{applications-to-ep}, $E$ and $F$ have the same conductor $108p^2$. If $p$ is not too large, we can use either the LMFDB or the algorithm of \cite{BeGhRe} to prove, by enumerating all possible curves, that $E$ and $F$ are isogenous. Could we prove that $E$ and $F$ are isogenous in some other way when $p$ is too large for such a search? 
\begin{itemize}[noitemsep,label=$-$] 
\item If $F$ has complex multiplication, this seems easy: it should not be difficult to show that $j(F)=j(E)=0$, so that $F$ has an equation of the form $y^2=x^3 \pm 54 \cdot 2^{\alpha}3^{\beta}p^{\gamma}$ (cf. e.g. \cite[\S 3.1.7]{BeGhRe}) with $0 \leq \alpha,\beta,\gamma < 6$. By the Hasse bound, one has $a_{\ell}(F)=a_{\ell}(E)$ for primes $2\leq \ell < \frac{p^2}{16}$. Using the explicit calculation of \cite[Example II.10.6]{AEC2} and the Serre weight calculation of \cite{Kraus-thesis}, it should not be too difficult to show that $E$ and $F$ are isogenous.     
\item By the Sturm bound (e.g. \cite[Lemma 5.2]{Zywina-Surj}), it suffices to show that $a_q(E)=a_q(F)$ for every prime $q \leq 9p(p+1)$ --- yet applying the Hasse bound to the congruence only yields the equality when $q \leq \frac{p^2}{16}$. Can the argument be improved in this situation?
\item For $p$ large enough, the effective Chebotarev theorem following from GRH implies the result by \cite[Th\'eor\`eme 21]{Serre-Chebotarev} (although not trivially so, since the argument requires a control on the conductor of $F$ which is given by Theorem \ref{applications-to-ep}). Checking the proof carefully and with the previous item, it turns out that a much weaker effective Chebotarev would be needed: instead of the bound ``$\pi_C(x) \geq 1$ if $x \geq c(\log{d_E})^{2}$'' (cf. \cite[Th\'eor\`eme 5]{Serre-Chebotarev}, whose notation we adapt) implied by GRH, we would need a bound of the form ``$\pi_C(x) \geq 1$ if $x \geq cd_E^{c'}$'', where $c,c' > 0$ are constants, and $c'$ is below a certain numerical bound. 
\item As previously mentioned, the other classical method for determining integral points on modular curves is \emph{Runge's method}. This method only works when the so-called \emph{Runge condition} is satisfied: it is necessary to have at least two Galois orbits of cusps. However, one can check (applying a Galois twist to the description of the cuspidal subscheme found in \cite[Theorem 10.9.1]{KM}) that $X_E^{\alpha}(p)$ has a unique Galois conjugacy class of cusps, because under our assumptions, the Galois action on $E[p](\Qbar) \backslash \{0\}$ is transitive.
\end{itemize} 
\item \emph{Other factors of the Jacobian with rank zero:} As hinted in this introduction, one can sometimes construct a different quotient $A$ of $J_E^{\alpha}(p)$ (still when $E$ has bad reduction at $p$) in which the functional equation does not force $L(A,1)$ to vanish. However, there are various technical difficulties in adapting the argument described in this paper to the quotient $A$: determining whether it even occurs\footnote{This involves the following question of independent interest: which characters $\F_{p^2}^{\times}/\F_p^{\times} \rar \C^{\times}$ arise as parameters of cuspidal representations of $\GL{\Q_p}$ attached to newforms in $\mathcal{S}_2(\Gamma_0(p^2))$? It should be possible to experiment using the results of \cite{LW}.}, whether $L(A,1)$ vanishes, or what the equivalent of what we call ``exceptional primes'' should be, are all non-trivial problems in general, although we expect them to be tractable in any given situation, at least in principle. 

\end{enumerate}

Let us now describe the structure of this paper. Section \ref{sect:XEP-general} describes the general set-up of the paper and discusses the structure of the Tate module of the Jacobian $J_E(p)$ of $X_E(p)$ for a general elliptic curve $E/\Q$. In fact, it is more convenient --- and this is what we do --- to work with finite \'etale group schemes $\Gamma$ (over suitable base rings $R$) that are \'etale-locally isomorphic to $(\Z/p\Z)^{\oplus 2}$, and consider the compactification of the curve $Y_{\Gamma}(p)$ parametrizing elliptic curves $F$ along with an isomorphism of group schemes $\Gamma \rar F[p]$. Section \ref{sect:cartan-bsd} considers the specialization of this situation where the image of $\rho_{E,p}$ is contained in the normalizer of a non-split Cartan subgroup. We recall some local results, refine the decomposition of the Tate module of $J_E(p)$ found in Section \ref{sect:XEP-general}, and apply the machinery of the Beilinson--Flach Euler system in order to prove Theorem \ref{existence-of-quotient-bk}. Section \ref{sect:fimm-ccol} proves Theorem \ref{formal-immersion-applied}. Section \ref{sect:effectivity-results} is concerned with the practical application of Theorems \ref{existence-of-quotient-bk} and \ref{formal-immersion-applied} to a specific elliptic curve $E/\Q$. We describe how to check the assumption of small inertia image at $p$, describe some properties of the newform $g$ appearing in Theorem \ref{existence-of-quotient-bk}, prove Theorem \ref{nonvan-l-function}, and deduce Theorems \ref{applications-to-ep} and \ref{mytheo-23} as consequences. \\

\textbf{Remark:} The results of this paper are in large part (but not entirely) contained in the author's PhD dissertation \cite{Studnia-thesis}. However, we were able to simplify many arguments of \emph{loc.cit.}, leading to somewhat easier proofs of Theorems \ref{existence-of-quotient-bk} and \ref{formal-immersion-applied}. Theorem \ref{nonvan-l-function} and its applications to any $p \equiv 2,5 \mod{9}$ (i.e. Theorem \ref{applications-to-ep}) do not appear in the author's PhD dissertation. \\

\textbf{Acknowledgements:} Since this text was part of my PhD dissertation, I would like to thank my former supervisor Lo\"ic Merel for his support as well as his many insightful questions and suggestions, and to Pierre Parent and Henri Darmon for their reading and commenting on the PhD dissertation this work was based upon. I am grateful to my colleagues Jan Vonk, David Lilienfeldt, Jonathan Love for many fruitful discussions regarding this work, and in particular regarding the use of quaternion algebras appearing in the proof of Theorem \ref{nonvan-l-function}. I would also like to thank Ye Tian, Jean-Loup Waldspurger and Andrew Sutherland for discussions on this very topic. I would also like to thank Pascal Molin for his help with some computational aspects of this work, Hendrik Lenstra and Pieter Moree for our discussions on the ``exceptional primes'', and David Loeffler for useful discussions regarding the twisted Bloch--Kato conjecture in the case where the two modular forms do not have coprime conductors. 

The author was partly supported by ERC Grant ``GAGARIN'' $101076941$.

\textbf{Notations:}
\begin{itemize}[noitemsep,label=\tiny$\bullet$]
\item If $X$ is a scheme and $x \in X$, $\mfk{m}_{X,x}, \kappa(x)$ denote the maximal ideal of $\OO_{X,x}$ and the residue field at $x$ respectively. 
\item $\mu_N$ denotes the finite flat group scheme $\Sp{\Z[X]/(X^N-1)}$ over $\Sp{\Z}$. 
\item $\mu_N^{\times}$ denotes the finite \'etale $\Z[1/N]$-scheme $\Sp{\Z[\zeta_N,1/N]}$, and $\zeta_N$ is a primitive $N$-th root of unity for $N \geq 1$.  
\item If $\ell$ is a prime and $n \in \Z$, $v_{\ell}(n)$ is the $\ell$-adic valuation of $n$. We extend this notation to $\Q_{\ell}$, and will in this case often omit the $\ell$ in the index. 
\item If $F$ is a field, $G_F$ denotes the absolute Galois group of $F$. If $F$ is a local field, $I_F$ denotes the inertia group and $I_F^+$ the wild inertia group. If $F$ is a number field and $v$ is a place of $F$, $G_v$ and $I_v$ denote some choice of decomposition and inertia groups of $G_F$ at the place $v$.  
\item As a consequence of the above choice, we will use $\mrm{id}_n$ to denote the identity matrix in a space of $n \times n$ matrices. 
\item For any $n \geq 1$, the cyclotomic character modulo $n$ is denoted by $\omega_n$. 
\end{itemize}

\section{The curve $X_{\Gamma}(p)$ and its Jacobian}
\label{sect:XEP-general}

\subsection{The general set-up}
\label{subsect:general-setup}

Let us recall some notations and constructions from \cite[Section 4]{Studnia-moduli}. Throughout Section \ref{subsect:general-setup}, $p \geq 3$ denotes a prime, and (unless otherwise specified) $m \geq 1$ denotes a square-free integer which is coprime to $p$. 

\defi{If $R$ is a $\Z[1/p]$-algebra, a \emph{$p$-torsion group} over $R$ is a finite \'etale commutative $R$-group scheme $\Gamma$ which is \'etale-locally isomorphic to $(\Z/p\Z)^{\oplus 2}$.}

Let $\Gamma$ be a $p$-torsion group over some $\Z[1/p]$-algebra $R$. There exists a fine moduli space $Y_{\Gamma}(p,\Gamma_0(m))$ parametrizing triples $(E/S,\iota,C)$, where $S$ is a $R$-scheme and $E$ is an elliptic curve over $S$, $\iota: \Gamma_S \rar E[p]$ is an isomorphism of $S$-group schemes, and $C \leq E$ is a finite locally free $S$-group scheme of degree $m$, admitting fppf-locally over $S$ a generator in the sense of \cite[(3.4)]{KM} (we will say for short that $C$ is cyclic). Moreover, $Y_{\Gamma}(p,\Gamma_0(m))$ is affine, flat of relative dimension one over $R$ with geometrically reduced Cohen--Macaulay fibres, smooth except at the supersingular $j$-invariants in characteristic dividing $m$, and the $j$-invariant $j: Y_{\Gamma}(p,\Gamma_0(m)) \rar \mathbb{A}^1_R$ is flat. We write $Y_{\Gamma}(p)$ when $m=1$, since any elliptic curve has a unique cyclic subgroup of degree one, the zero section. 

The $Y_{\Gamma}(p,\Gamma_0(m))$ are equipped with finite locally free \emph{degeneracy maps} defined as follows: if $r,s \geq 1$ are such that $rs \mid m$, $D_{r,s}: Y_{\Gamma}(p,\Gamma_0(m)) \rar Y_{\Gamma}(p,\Gamma_0(s))$ maps a triple $(E/S,\iota,C)$ to $((E/C[r])/S,\pi\circ \iota, \pi(C)[s])$, where $C[r]$ (resp. $\pi(C)[s]$) is the ``standard cyclic subgroup'' of $C$ (resp. of $\pi(C)$, which is cyclic of degree $m/r$ by \cite[Theorem 6.7.4]{KM}) of degree $r$ (resp. $s$) in the sense of \cite[Theorem 6.7.2]{KM}\footnote{Since $m$ is square-free, $r$ is prime to $m/r$, so, by \cite[Proposition 1.7.2]{KM} and its proof, $C[r]$ is simply the $r$-torsion subgroup of $C$, whence the notation.}, and $\pi: E \rar E/C[r]$ is the projection. 

\defi[we-pairing]{A \emph{Weil pairing} on a $p$-torsion group $\Gamma$ over $R$ is a morphism of schemes $\alpha: \Gamma \times_R \Gamma \rar (\mu_p)_R$ which is bilinear, alternating and fibrewise perfect.}

The terminology comes from the following situation, which is essentially the one we are interested in: let $E/R$ be an elliptic curve over a $\Z[1/p]$-algebra $R$, then $E[p]$ is a $p$-torsion group, and its usual Weil pairing (normalized as in e.g. \cite[(2.8)]{KM}) is a Weil pairing in the sense of Definition \ref{we-pairing}. 

\rem{Let $k$ be a field of characteristic not dividing $p$ with separable closure $k_s$. Then, by the theory of the fundamental group \cite[Exp. V, Th\'eor\`eme 4.1, \S\S 7, Proposition 8.1]{SGA1}, the datum of a $p$-torsion group $\Gamma$ over $k$ (up to isomorphism) is equivalent to that of the equivalence class of the representation $\rho: G_{k} \rar \GL{\F_p}$ given by matrix of the action of $G_k$ in some basis of $\Gamma(k_s)$. The group $\Gamma$ admits a Weil pairing over $k$ if and only if $\det{\rho}$ is the cyclotomic character modulo $p$; specifying the Weil pairing is equivalent to fixing a primitive $p$-th root of unity $\xi \in k_s$ and a $\F_p$-basis $(P,Q)$ of $\Gamma(k_s)$, up to the relation $((P,Q),\xi) \sim ((aP+bQ,cP+dQ),\xi^{ad-bc})$ for $\begin{pmatrix}a & b\\c & d\end{pmatrix} \in \GL{\F_p}$. 

Suppose now that $k$ is a number field and $v$ is a place of $k$ not dividing $p$. By specialization theory (\emph{loc.cit.} and Proposition 8.2 of \emph{op.cit.}), the group $\Gamma$ extends to a $p$-torsion group $\Gamma'$ over $\OO_{k,v}$ if and only if $\rho$ is unramified at $v$. In such a situation, $\Gamma$ uniquely determines $\Gamma'$ and any Weil pairing on $\Gamma$ uniquely extends to a Weil pairing on $\Gamma'$.}

If $\Gamma$ is a $p$-torsion group over the $\Z[1/p]$-algebra $R$, then the functor mapping a $R$-scheme $S$ to the collection of Weil pairings on $\Gamma_S$ (which is obviously a pre-sheaf over $\Sch_R$ in $\F_p^{\times}$-torsors) is represented by a finite \'etale $\F_p^{\times}$-torsor $\mrm{Pair}(\Gamma)$. 

One has a map $\mrm{We}: Y_{\Gamma}(p,\Gamma_0(m)) \rar \mrm{Pair}(\Gamma)$ sending a triple $(E/S,\iota,C)$ to 
\[\Gamma_S \times_S \Gamma_S \overset{\iota \times \iota}{\rar} E[p] \times_S E[p] \overset{\mrm{We}_{E[p]/S}}{\rar} (\mu_p)_S,\] where $\mrm{We}_{E[p]/S}$ denotes the Weil pairing for $E/S$. The morphism $\mrm{We}$ map is flat of relative dimension one, its geometric fibres are reduced, Cohen--Macaulay, and connected, and it is smooth except possibly at the supersingular $j$-invariants in residue characteristic dividing $m$ by \cite[Theorem 1]{Studnia-moduli}. Furthermore, $\mrm{We}$ is compatible with the degeneracy maps, in the sense that the following diagram commutes:

\[
\begin{tikzcd}[ampersand replacement=\&]
Y_{\Gamma}(p,\Gamma_0(m)) \arrow{d}{D_{r,s}} \arrow{r}{\mrm{We}}\& \mrm{Pair}(\Gamma) \arrow{d}{r}\\
Y_{\Gamma}(p,\Gamma_0(s))  \arrow{r}{\mrm{We}}\& \mrm{Pair}(\Gamma)
\end{tikzcd}
\]

If $s: \Gamma \rar \Gamma'$ is an isomorphism of $p$-torsion groups, it induces (contravariant) functorial isomorphisms $Y_{\Gamma'}(p,\Gamma_0(m)) \rar Y_{\Gamma}(p,\Gamma_0(m))$ which commute to all the degeneracy maps. Thus the group functor $\mrm{Aut}(G)$ acts on the \emph{right} on the collection of $Y_{\Gamma}(p,\Gamma_0(m))$ (and preserves the degeneracy maps). 

Since every $p$-torsion group is \'etale-locally isomorphic to a base change of the constant $p$-torsion group scheme $G_0$ over $R_0 := \Z[1/p]$, let us discuss the case of $G_0$ some more. In this case, $\mrm{Pair}(G_0)$ is the scheme $\mu_p^{\times}$ of primitive $p$-th roots of unity (i.e. the closed open subscheme of $(\mu_p)_{R_0}$ where its unit section has been removed, i.e. $\mu_p^{\times}=\Sp{\frac{R_0[t]}{(1+t+\ldots+t^{p-1})}}$). Moreover, the functor of points of $Y_{G_0}(p,\Gamma_0(m))$ can be identified with the functor mapping a $R$-scheme $S$ to the collection of triples $(E,(P,Q),C)$ where $E/S$ is an elliptic curve, $P,Q \in E[p](S)$ are two sections such that $(P_{\overline{s}},Q_{\overline{s}})$ forms a basis of $E[p](\overline{s})$ for any geometric point $\overline{s}$ of $S$, and $C$ is a cyclic subgroup of $E$ of degree $m$. We can then let $\GL{\F_p}$ act on $Y_{G_0}(p,\Gamma_0(m))$ (compatibly with all degeneracy maps) as follows: $\begin{pmatrix}a & b\\c & d\end{pmatrix} \cdot (E/S,(P,Q),C) = (E/S,(aP+bQ,cP+dQ),C)$. It is not difficult to check that the action of $M \in \GL{\F_p}$ defined in this manner is the right action of $M^T \in \mrm{Aut}(G_0)(R)$ described in the previous paragraph. 

Moreover, this action of $\GL{\F_p}$ is compatible with $\mrm{We}$ in that the following diagram commutes for $M \in \GL{\F_p}$: 

\[
\begin{tikzcd}[ampersand replacement=\&]
Y_G(p,\Gamma_0(m)) \arrow{d}{M} \arrow{r}{\mrm{We}}\& \mrm{Pair}(G) \arrow{d}{\det{M}}\\
Y_G(p,\Gamma_0(m))  \arrow{r}{\mrm{We}}\& \mrm{Pair}(G)
\end{tikzcd}
\]

Let us now assume that $R$ is a regular excellent Noetherian $\Z[1/p]$-algebra and $\Gamma$ is a $p$-torsion group over $R$. 

The curve $Y_{\Gamma}(p,\Gamma_0(m))$ admits a compactification $X_{\Gamma}(p,\Gamma_0(m))$ whose formation commutes with any base change of regular excellent Noetherian $\Z[1/p]$-algebras, and the maps $j, \mrm{We}$ extend to a finite locally free map $(j,\mrm{We}): X_{\Gamma}(p,\Gamma_0(m)) \rar \mathbb{P}^1_R \times \mrm{Pair}(\Gamma)$. The extended map $\mrm{We}$ has relative dimension one, connected reduced Cohen--Macaulay geometric fibres, and it is smooth outside the supersingular $j$-invariants in characteristic dividing $m$. 

The \emph{cuspidal subscheme} of $X_{\Gamma}(p,\Gamma_0(m))$ is the reduced closed subscheme corresponding to the closed subspace $\{j=\infty\}$; it is a relative effective Cartier divisor of $X_{\Gamma}(p,\Gamma_0(m))$ (over $\Sp{R}$) and it is finite \'etale over $\Sp{R}$. Its formation commutes with any base change of regular excellent Noetherian $\Z[1/p]$-algebras.

The above constructions (degeneracy maps, functoriality with respect to $\Gamma$, left action of $\GL{\F_p}$ in the case of $G_0$) all extend to the compactifications, remain finite locally free (or isomorphisms), and the above diagrams all commute when every instance of $Y_{\Gamma}$ is replaced with $X_{\Gamma}$.

We now define Hecke correspondences and the corresponding operators on the Jacobians of $X_{\Gamma}(p)$. Since these ``curves'' over $R$ do not have geometrically connected fibres, the standard definitions and results about curves and Jacobians need to be adapted. We sketch the necessary modifications in Appendix \ref{appendix:algeom}. We let $J_{\Gamma}(p)$ denote the relative Jacobian of $X_{\Gamma}(p)/R$: it represents the fppf-sheafification of the functor mapping a $R$-scheme $T$ to the collection of line bundles $\mathcal{L}$ on $X_{\Gamma}(p) \times_R T$ such that for any $t \in T$, the restriction of $\mathcal{L}$ to any connected component of $X_{\Gamma}(p) \times_R \Sp{\kappa(t)}$ has degree zero.

For every prime $\ell \neq p$, the \emph{Hecke operator} $T_{\ell}$ is the endomorphism of $J_{\Gamma}(p)$ defined by $(D_{\ell,1})_{\ast}(D_{1,\ell})^{\ast}$, where $D_{1,\ell}, D_{\ell,1}: X_{\Gamma}(p,\Gamma_0(\ell)) \rar X_{\Gamma}(p)$ are the two previously defined degeneracy maps. We let $n \in \F_p^{\times}$ act by the push-forward of the map $X_{\Gamma}(p) \overset{[n]^{\ast}}{\rar} X_{\Gamma}(p)$, where $[n]: \Gamma \rar \Gamma$ is the isomorphism of multiplication by $n$. To be completely unambiguous, $[n]$ acts on the non-cuspidal $S$-point $x = (E/S, \iota: \Gamma_S \rar E[p])$ of $X_{\Gamma}(p)$ as follows: one has $[n]^{\ast}(x) = (E/S, n\cdot \iota: \Gamma_S \rar E[p])$. We denote by $\diam{n}$ the action of $n$ on $J_{\Gamma}(p)$. 

\prop[hecke-commutes]{The sub-algebra $\mathbb{T}_{\Gamma/R}$ of $\mrm{End}(J_{\Gamma}(p))$ generated by the $T_q$ and the actions of the $n \in \F_p^{\times}$ is commutative, and furthermore, $n \in \F_p^{\times} \mapsto \diam{n}\in \Aut{J_{\Gamma}(p)}$ is a group homomorphism.}

\demo{It is clear that $[n]^{\ast}\circ [m]^{\ast}=[mn]^{\ast}$ for $n,m \in \F_p^{\times}$ already as $R$-automorphisms of $X_{\Gamma}(p)$, so we can use ordinary functoriality for Jacobians. The $[n]^{\ast}$ are defined on $X_{\Gamma}(p,\Gamma_0(m))$ and commute to all degeneracy maps, so, by ordinarity functoriality on the various Picard pre-sheaves $P^0_{X_{\Gamma}(p,\Gamma_0(m))}$ (see Appendix \ref{appendix:algeom}), the $\diam{n}$ commute to the $T_q$. 

Let $q,r \neq p$ be distinct primes, we are done if we show that $T_q,T_r$ commute as endomorphisms of the Picard pre-sheaf $P^0_{X_{\Gamma}(p)/R}$. We show in fact that $T_qT_r = (D_{qr,1})_{\ast}(D_{1,1})^{\ast}$ (where $D_{1,1}, D_{qr,1}$ are two degeneracy maps $X_{\Gamma}(p,\Gamma_0(qr)) \rar X_{\Gamma}(p)$), which is enough. Unpacking the definitions, it is enough to show that in the following diagram, 
\[
\begin{tikzcd}[ampersand replacement=\&]
X_{\Gamma}(p,\Gamma_0(qr)) \arrow{r}{D_{r,q}} \arrow{d}{D_{1,r}} \& X_{\Gamma}(p,\Gamma_0(q)) \arrow{d}{D_{1,1}}\\
X_{\Gamma}(p,\Gamma_0(r)) \arrow{r}{D_{r,1}} \& X_{\Gamma}(p)
\end{tikzcd}
\]

the maps $D_{1,1}^{\ast}(D_{r,1})_{\ast}$ and $(D_{r,q})_{\ast}D_{1,r}^{\ast}$ agree as maps $P^0_{X_{\Gamma}(p,\Gamma_0(r))/R} \rar P^0_{X_{\Gamma}(p,\Gamma_0(q))/R}$. We apply Lemma \ref{mixed-pic0-functoriality}: the diagram commutes, all maps are finite locally free, the diagram is Cartesian above $Y_G(p)$ (since the base change of this Cartesian diagram by $Y_{\Gamma}(p) \rar X_{\Gamma}(p)$ is the same diagram with all instances of $X_{\Gamma}$ replaced by $Y_{\Gamma}$), and $Y_{\Gamma}(p) \rar X_{\Gamma}(p)$ is an affine open immersion satisfying the injectivity requirement since it is a finite locally free base change (by the $j$-invariant) of the affine immersion $\mathbb{A}_1^R \rar \mathbb{P}^1_R$, which clearly satisfies the injectivity condition.}

Next, we check that the Hecke algebra (along with its Hecke operators) does not depend on the choice of $p$-torsion group. We first need a Lemma that we record here for later use. 

\lem[endo-abvar-spe]{Let $R$ be a connected ring and $R'$ be a $R$-algebra. Let $A$ be an abelian scheme over $R$. Then $\mrm{End}_R(A) \rar \mrm{End}_{R'}(A)$ is injective.} 

\demo{When $R,R'$ are both fields, this follows from \cite[Theorem 12.5]{MilAb} since $A$ and $A_{R'}$ have the same $\ell$-adic Tate module for any $\ell$ invertible in $R$. We are therefore reduced to the case where $R'$ is a residue field $k$ of $R$. In this case, this is the classical rigidity result \cite[Proposition 6.1]{GIT}.}

\prop[hecke-uniquely-defined]{If $\Gamma,\Gamma'$ are $p$-torsion groups over the regular excellent Noetherian $\Z[1/p]$-algebras $R,R'$, then there exists a unique isomorphism $\mathbb{T}_{\Gamma/R} \rar \mathbb{T}_{\Gamma'/R'}$ sending each $T_{\ell}$ to $T_{\ell}$ and each $\diam{n}$ to $\diam{n}$. If $R'$ is a $R$-algebra and $\varphi: \Gamma_{R'} \rar \Gamma'$ is an isomorphism, then the induced map $J_{\Gamma'}(p) \rar J_{\Gamma_{R'}}(p)$ (by push-forward functoriality of the isomorphism $\varphi^{\ast}: X_{\Gamma'}(p) \rar X_{\Gamma_{R'}}(p)$) is $\mathbb{T}_{\Gamma'/R'}$-linear.  
}

\demo{Let us call \emph{convenient} a morphism $\mathbb{T}_{\Gamma/R} \rar \mathbb{T}_{\Gamma'/R'}$ sending each $T_{\ell}$ to $T_{\ell}$ and each $\diam{n}$ to $\diam{n}$. If a convenient morphism exists, it is necessarily unique and surjective, and the composition of two convenient morphisms is convenient. Convenient morphisms exist in the following situations:
\begin{itemize}[noitemsep,label=\tiny$\bullet$]
\item one has convenient morphisms $\mathbb{T}_{\Gamma/R} \rar \mathbb{T}_{\Gamma_{R'}/R'}$ for any morphism $R \rar R'$ of regular excellent Noetherian $\Z[1/p]$-algebras, which are injective --- thus isomorphisms --- by Lemma \ref{endo-abvar-spe},
\item if $\varphi: \Gamma \rar \Gamma'$ is an isomorphism, the isomorphism $J_{\Gamma'}(p) \rar J_{\Gamma}(p)$ induced by push-forward of $\varphi^{\ast}$ induces a convenient morphism, which is functorial with respect to $\varphi$, and is therefore an isomorphism. In particular, $\mathbb{T}_{(G_0)_R/R}$ commutes to the action of $\GL{\F_p}$.
\end{itemize}
$\Gamma,\Gamma'$ are \'etale-locally isomorphic to $G_0$, so the conclusion follows. }

We thus write $\mathbb{T}:= \mathbb{T}_{G_0/\Z[1/p]}$ and call it the \emph{Hecke algebra}; by the above Proposition, it acts on every $J_{\Gamma}(p)$ in an obvious way that respects base change and isomorphisms of $p$-torsion groups. By \cite[Theorem 12.5]{MilAb}, $\mathbb{T} \simeq \mathbb{T}_{(G_0)_{\Q}/\Q}$ is a finite free commutative $\Z$-algebra. If $H$ is a subgroup of $\GL{\F_p}$ containing $\F_p^{\times}\mrm{id}_2$, we write $\mathbb{T}[H] := \mathbb{T} \otimes_{\Z[\F_p^{\times}]} \Z[H]$, where $n \in \F_p^{\times}$ is sent to $\diam{n} \in \mathbb{T}^{\times}$ and to $n\mrm{id}_2 \in H$; it is a finite free $\mathbb{T}$-module which also acts on $J_{G_0}(p)$ (on the left). 

We define the elements $T_n \in \mathbb{T}$ for $n \geq 1$ prime to $p$ as follows: 
\begin{itemize}[noitemsep,label=\tiny$\bullet$]
\item $T_1$ is the identity,
\item if $\ell$ is a prime and $r \geq 2$, then $T_{\ell^r} = T_{\ell}T_{\ell^{r-1}}-\ell \diam{\ell} T_{\ell^{r-2}}$, 
\item if $n,m \geq 1$ are coprime and both prime to $p$, then $T_{nm}=T_nT_m$. 
\end{itemize}

The main interest of the $T_n$ is that the $T_n\diam{m}$ generate $\mathbb{T}$ as an abelian group.

\lem[hecke-preparatory-grading]{Let $n,m \geq 1$ be integers such that $p \nmid nm$. There exist $a_1,\ldots,a_r \in \F_p^{\times}$, integers $b_1, \ldots, b_r \geq 1$ such that $a_i^2b_i \equiv 1\mod{p}$, $T_nT_m = \sum_{i=1}^r{\diam{a_i} T_{b_i}}$, and every $b_i$ is only divisible by primes which divide $nm$. }

\demo{By multiplicativity, we easily reduce to the case where $m,n$ are powers of the same prime $\ell \neq p$. The recurrence relation solves the case where $\min(n,m) \leq \ell$. If now $\ell^2 \mid n \mid m$, one has $T_mT_n = (T_{m}T_{n/\ell})T_{\ell}-\ell \diam{\ell} T_mT_{n/\ell^2}$, and the conclusion holds by induction $\min(m,n)$. }

Let now $\Gamma$ be a $p$-torsion group over some regular excellent Noetherian $\Z[1/p]$-algebra $R$. For every Weil pairing $\alpha$ for $\Gamma$, its pre-image under $\mrm{We}$ defines a closed open $R$-subscheme $X_{\Gamma}^{\alpha}(p,\Gamma_0(m))$ of $X_{\Gamma}(p,\Gamma_0(m))$, and $X_{\Gamma}^{\alpha}(p,\Gamma_0(m)) \rar \Sp{R}$ is a flat proper map of relative dimension one with reduced connected geometric fibres, which is smooth outside the supersingular $j$-invariants in characteristic dividing $m$. In particular, $X_{\Gamma}^{\alpha}(p)$ has a (classical) Jacobian $J_{\Gamma}^{\alpha}(p)$, which is an abelian scheme over $R$. 

One sees easily that the intersection $Y_{\Gamma}^{\alpha}(p,\Gamma_0(m))$ of $X_{\Gamma}^{\alpha}(p,\Gamma_0(m))$ and $Y_{\Gamma}(p,\Gamma_0(m))$ is the fine moduli space of triples $(E/S, \iota,C)$, where $E$ is an elliptic curve over the $R$-scheme $S$, $\iota: \Gamma_S \rar E[p]$ is an isomorphism such that the pairing $\Gamma_S \times_S \Gamma_S \overset{\iota \times \iota}{\rar} E[p] \times E[p] \overset{\mrm{We}_{E[p]/S}}{\rar} (\mu_p)_S$ is exactly $\alpha_S$, and $C$ is a cyclic subgroup of $E$ of degree $m$. One clearly has \[X_{\Gamma}(p,\Gamma_0(m))=\coprod_{n \in \F_p^{\times}}{X_{\Gamma}^{\alpha^n}(p,\Gamma_0(m))};\] in particular, $J_{\Gamma}(p)$ is canonically isomorphic to $\prod_{n \in \F_p^{\times}}{J_{\Gamma}^{\alpha^n}(p)}$. This decomposition lets us define a $\F_p^{\times}$-grading on $\mathbb{T}$ as follows:

\prop[hecke-grading]{For $n \geq 1$ prime to $p$ (resp. $n \in \F_p^{\times}$), $T_n$ (resp. $\diam{n}$) sends the subvariety $J_{\Gamma}^{\alpha}(p)$ of $J_{\Gamma}(p)$ into its subvariety $J_{\Gamma}^{\alpha^n}(p)$ (resp. $J_{\Gamma}^{\alpha^{n^2}}(p)$). 

Therefore, one has a decomposition of abelian groups $\mathbb{T} = \bigoplus_{n \in \F_p^{\times}}{\mathbb{T}_n}$, where $\mathbb{T}_n$ is generated as an abelian group by the $\diam{b} \cdot T_a$ with $b^2a \equiv 1\mod{p}$, and $\mathbb{T}_n\mathbb{T}_{n'} \subset \mathbb{T}_{nn'}$ for $n,n' \in \F_p^{\times}$. 
}

\demo{First, $\diam{n}$ induces an isomorphism $X_{\Gamma}^{\alpha}(p) \rar X_{\Gamma}^{\alpha^{n^2}}(p)$. By multiplicativity, it is enough to show that $T_{\ell}$ sends the Picard pre-sheaf $P^0_{X_{\Gamma}^{\alpha}(p)}$ into $P^0_{X_{\Gamma}^{\alpha^{\ell}}(p)}$ for $\ell \neq p$ prime. Now, the maps $D_{1,1}, D_{\ell,1}: X_{\Gamma}(p,\Gamma_0(\ell)) \rar X_{\Gamma}(p)$ send $X_{\Gamma}^{\alpha}(p,\Gamma_0(\ell))$ into $X_{\Gamma}^{\alpha}(p)$ and $X_{\Gamma}^{\alpha^{\ell}}(p)$ respectively, so the first sentence of the statement of Proposition \ref{hecke-grading} follows. 

To prove the second part of the claim, let $\mathbb{T}_n$ be the abelian subgroup of $\mathbb{T}$ generated by $\diam{b}T_a$ where $a \geq 1$ is an integer, $b \in \F_p^{\times}$ and $b^2a\equiv n \mod{p}$, then we know that $\mathbb{T}_n\mathbb{T}_{n'}\subset \mathbb{T}_{nn'}$ by Lemma \ref{hecke-preparatory-grading}. The first part of the proof shows that the $\mathbb{T}_n$ are in direct sum. Finally, Lemma \ref{hecke-preparatory-grading} implies that $\mathbb{T}$ is generated as an abelian group by the $\diam{b} \cdot T_a$, so that $\sum_{n \in \F_p^{\times}}{\mathbb{T}_n} = \mathbb{T}$. 
}

Let now $H$ be a subgroup of $\GL{\F_p}$ containing $\F_p^{\times}\mrm{id}_2$, we define the following $\F_p^{\times}$-grading on $\mathbb{T}[H]$: for $n \in \F_p^{\times}$, $\mathbb{T}[H]_n$ is generated by the $T_mh$ for $m \geq 1$ and $h \in H$ such that $m\det{h} \equiv n \mod{p}$. 

\lem[hecke-gl2-grading]{One has $\mathbb{T}[H] = \bigoplus_{n \in \F_p^{\times}}{\mathbb{T}[H]_n}$. Moreover, if $\Gamma$ is the constant $p$-torsion group scheme over a regular excellent $\Z[1/p]$-algebra $R$ and $\alpha: \Gamma \times_R \Gamma \rar (\mu_p)_R$ is a Weil pairing, then, for every $t \in \mathbb{T}[\GL{\F_p}]_n$, the morphism $t: J_{\Gamma}(p) \rar J_{\Gamma}(p)$ sends $J_{\Gamma}^{\alpha}(p)$ into $J_{\Gamma}^{\alpha^n}(p)$. }

\demo{Let $g \in \GL{\F_p}$. Since $g: X_{\Gamma}(p) \rar X_{\Gamma}(p)$ sends $X_{\Gamma}^{\alpha}(p)$ into $X_{\Gamma}^{\alpha^{\det{g}}}(p)$ (because of its action through $\mrm{We}$), its push-forward to $J_{\Gamma}(p)$ sends $J_{\Gamma}^{\alpha}(p)$ into $J_{\Gamma}^{\alpha^{\det{g}}}(p)$. By Proposition \ref{hecke-grading}, the endomorphism of $J_{\Gamma}(p)$ induced by any $t \in \mathbb{T}[H]_n$ sends $J_{\Gamma}^{\alpha}(p)$ into $J_{\Gamma}^{\alpha^n}(p)$. In particular, the $\mathbb{T}[H]_n$ are in direct sum, and they generate $\mathbb{T}[H]$ as an abelian group, whence the conclusion. }

Let us finally mention, for future reference, that $\mathbb{T}$ possesses the following interesting involution:

\defi[involution-of-hecke]{The rule $t \in \mathbb{T}_n \mapsto \diam{n}^{-1}t \in \mathbb{T}_{n^{-1}}$ defines an involution $w$ of the Hecke algebra. }

\rem{Let $k$ be a field of characteristic not $p$, and $\Gamma$ be a $p$-torsion group over $k$. Let $\chi: G_k \rar \{\pm 1\}$ be a character, then $\Gamma$ admits a quadratic twist $\Gamma'$ by $\chi$ (which is still a $p$-torsion group): if $\Gamma$ comes from the representation $\rho: G_k \rar \GL{\F_p}$, then $\Gamma'$ comes from $\rho \otimes \chi$. By a slight adaptation of \cite[Proposition 5.1.4, Theorem 5]{Studnia-moduli}, one constructs a canonical isomorphism $\mrm{Pair}(\Gamma) \simeq \mrm{Pair}(\Gamma')$ and $\mathbb{P}^1_k \times \mrm{Pair}(\Gamma)$-isomorphisms $X_{\Gamma}(p,\Gamma_0(m)) \simeq X_{\Gamma'}(p,\Gamma_0(m))$ for every square-free $m$ that commute to the degeneracy maps and are given on non-cuspidal points by the quadratic twists of elliptic curves. This construction induces a $\mathbb{T}$-linear isomorphism $J_{\Gamma}(p) \rar J_{\Gamma'}(p)$. In other words, it is equivalent to find non-trivial rational points on $X_{\Gamma}(p)$ and $X_{\Gamma'}(p)$. This is why we will often freely replace the elliptic curves we are working with by their quadratic twists when convenient.  
}

\bigskip

\subsection{Reminders: modular forms, automorphic representations and compatible systems}
\label{subsect:langlands-stuff}

This Section is dedicated to recalling some well-known facts about the relations between modular forms, automorphic representations, and the attached Galois representations. Our main references are \cite[\S 2]{LW}, \cite{NTB}, \cite{BH}. The interested reader should consult the books \cite{Gelbart, Bump, Getz-Hahn} for more detail about automorphic representations and their links to modular forms, although they should bear in mind that the choices of normalization are not entirely equivalent. 

First, we adelize Dirichlet characters as in \cite[Proposition 3.1.2]{Bump} and \cite[\S 2]{LW}: given a Dirichlet character $\chi: (\Z/N\Z)^{\times} \rar \C^{\times}$, $\chi_{\mathbb{A}}$ is the unique continuous character $\mathbb{A}_{\Q}^{\times}/\Q^{\times} \rar \C^{\times}$ such that for every prime $\ell \nmid N$, the restriction of $\chi_{\mathbb{A}}$ to $\Q_{\ell}^{\times}$ is the character $\Q_{\ell}^{\times} \rar \Q_{\ell}^{\times}/\Z_{\ell}^{\times} \overset{n \mapsto \chi(\ell)^{v_{\ell}(n)}}{\rar} \C^{\times}$.

By construction, $\chi_{\mathbb{A}}$ is trivial on $(\hat{\Z}^{\times}) \cap (1+N\hat{\Z})$, it only depends on the primitive Dirichlet character attached to $\chi$, and the rule $\chi \mapsto \chi_{\mathbb{A}}$ preserves multiplication of Dirichlet characters.  

Recall that a complex representation $V$ of $\GL{\Q_p}$ for a prime $p$ is \emph{admissible} (cf. e.g. \cite[2.1]{BH}) if every vector has an open stabilizer and, for every compact open subgroup $K \leq \GL{\Q_p}$, $V^K$ is a finite-dimensional complex vector space. By Schur's lemma (cf. \cite[Remark 4.15]{Gelbart}), if $V$ is an irreducible admissible representation of $\GL{\Q_p}$, $\Q_p^{\times}\mrm{id}_2$ acts on $V$ through a continuous character $\Q_p^{\times} \rar \C^{\times}$, which is the \emph{central character} of $V$. We say that $V$ is furthermore \emph{spherical} if it has a fixed vector under $\GL{\Z_p}$: such a vector is then uniquely determined up to scalar (cf e.g. \cite[Theorem 4.23]{Gelbart}). 
 
To any newform $f \in \mathcal{S}_k(\Gamma_1(N))$ with character $\chi$ and any place $v$ of $\Q$ one can attach (by e.g. \cite[\S 3.6]{Bump}, see also \cite{LW}) an infinite-dimensional irreducible admissible unitary automorphic representation $\pi_{f}$ of $\GL{\mathbb{A}_{\Q}}$, which factors as the completed restricted tensor product of irreducible admissible representations $\pi_{f,v}$ of $\GL{\Q_v}$, all but finitely many of which are spherical, with respect to the lines $\pi_{f,v}^{\GL{\Z_v}}$. This construction has the following properties: 
\begin{itemize}[noitemsep,label=\tiny$\bullet$]
\item $\pi_{f,\infty}$ is a discrete series representation (cf. e.g. \cite[p. 58]{Gelbart}) that only depends on the weight $k$ of $f$,
\item if $v$ is a non-trivial normalized valuation of $\Q$ (attached to the finite place denoted by the same letter), $\pi_{f,v}$ is spherical if and only if $v \nmid N$, 
\item the central character of $\pi_f$ is precisely $\chi_{\mathbb{A}}$,
\item if $v$ is a place of $\Q$ attached to the prime $p \nmid N$, then $\pi_{f,v}$ is the principal series representation attached to the two unramified characters $\Q_v^{\times}/\Z_v^{\times} \rar \C^{\ast}$ mapping uniformizers to the two roots of $x^2-p^{(1-k)/2}a_p(f)X+\chi(p)$ (that the representation is unitary then follows from the Ramanujan conjecture proved by Deligne \cite{Deligne-Ram,Weil-I}). 
\end{itemize}

One very important feature of this construction is its compatibility with respect to twists. 

\defi[twist-of-newform]{Let $f \in \mathcal{S}_k(\Gamma_1(N))$ be a newform and $\chi: (\Z/C\Z)^{\times} \rar \C^{\times}$ be a Dirichlet character. The \emph{twist} of $f$ by $\chi$ is the unique newform (denoted by $f \otimes \chi$) in $\mathcal{S}_k(\Gamma_1(N'))$ (for some uniquely determined integer $N' \geq 1$ dividing $NC^2$) such that, for all $n \geq 1$ coprime to $NC$, one has $a_{n}(f\otimes \chi)=a_{n}(f)\chi(n)$.}

\rem{In the above definition, 
\begin{itemize}[noitemsep,label=\tiny$\bullet$]
\item The twist of $f$ by $\chi$ only depends on the primitive Dirichlet character attached to $\chi$,
\item If $\chi_1$ and $\chi_2$ are Dirichlet characters modulo $C$, one has $f \otimes (\chi_1\chi_2)=(f \otimes \chi_1) \otimes \chi_2$.
\end{itemize}}

\prop[compatibility-with-twists]{Let $f \in \mathcal{S}_k(\Gamma_1(N))$ be a newform and $\chi: (\Z/C\Z)^{\times} \rar \C^{\times}$ be a Dirichlet character. Then $\pi_{f \otimes \chi} = \pi_f \otimes \chi_{\mathbb{A}}(\det)$. }

\demo{Both sides are automorphic representations of $\GL{\mathbb{A}_{\Q}}$ and they agree at almost every component by definition of the twist; by strong multiplicity one (e.g. \cite[Theorem 11.4.3]{Getz-Hahn}) they are equal everywhere. }

We check another form of compatibility, that with actions of $\Aut{\C}$. This is perhaps slightly more complex than it appears, since the construction of $\pi_{f,v}$ involves the upper half-plane, and is in particular not obviously algebraic. Nonetheless, the following result holds:

\prop[compatibility-with-galois]{Let $f \in \mathcal{S}_k(\Gamma_1(N))$ be a newform. Let $\sigma\in \Aut{\C}$ and $v$ be a finite place of $\Q$. Then $\sigma(f) \in \mathcal{S}_k(\Gamma_1(N))$ is a newform and the representations $\pi_{\sigma(f),v} \otimes |\det|^{(1-k)/2}$ and $\sigma(\pi_{f,v} \otimes |\det|^{(1-k)/2})$ are isomorphic (where $|\cdot|$ is the absolute value and, the action of $\sigma$ consists in applying $\sigma$ to all the coefficients of an ``infinite matrix representation'').}

\demo{The first part is a well-known consequence of the Eichler--Shimura isomorphism (i.e. the theory of modular symbols). To show the second part of the claim, it is enough to show by strong multiplicity one that both sides are automorphic (which is \cite[Theorem 1.8.1]{Waldspurger} --- after taking some care when unpacking the definitions, which account for the ``twist'' by the non-algebraic power of $|\det|$) and that the isomorphism holds for all but finitely many $v$ (because it is obviously true when $v \nmid N$). } 

We now discuss compatible systems of Galois representations. 

\defi[compatible-systems]{Let $K,F$ be number fields and $S$ be an infinite set of finite places of $F$. A \emph{compatible system of Galois representations} of $G_K$ over $F$ \emph{supported over $S$} is a collection of continuous $F_v[G_K]$-modules $V_v$ over $v \in S$ satisfying the following conditions:
\begin{itemize}[noitemsep,label=\tiny$\bullet$]
\item There is a finite set $T$ of places of $K$ such that for every $v \in S$ and every finite place $w \notin T$ of $K$ prime to the residue characteristic of $v$, $V_v$ is unramified at $v$. 
\item Let $w$ be a finite place of $K$, then there is a Frobenius-semisimple Weil--Deligne representation $V(w)$ of $K_w$ over $\overline{F}$ (see e.g. \cite[(4.1.2)]{NTB}) such that, for every $v \in S$ prime to the residue characteristic of $w$, the Frobenius-semisimplification $V_{v,w}$ of the Weil--Deligne representation attached to $(V_v)_{|G_{K_w}}$ (see e.g. \cite[(4.2)]{NTB}) --- which is defined over the field $F_v$ --- satisfies $V_{v,w} \otimes_{F_v} \overline{F_v} \simeq V(w) \otimes_{\overline{F}} \overline{F_v}$ as Weil--Deligne representations of $K_w$ over $\overline{F_v}$ for any embedding $\overline{F} \rar \overline{F_v}$.  
\end{itemize}}

\rems{\begin{itemize}[noitemsep,label=\tiny$\bullet$]
\item The formulation of the second condition is complex, but for $w \notin T$, it can be stated in a much simpler way: we simply require that the characteristic polynomial of the Frobenius $\Fr_w$ acting on $V_v$ be contained in $F[t]$ and does not depend on $v$.
\item As a comparison, if $S$ is the set of finite places of $F$ and the second condition is restricted to the $w \notin T$, we recover what is called in \cite[\S 2.3, I-11 Definition, I-13 Remark]{SerreMcGill} a strictly compatible $F$-rational $v$-adic system of Galois representations of $K$. 
\item We do not include any $p$-adic Hodge-theoretic condition, as is often the case in applications of the notion to modularity lifting theorems (see e.g. \cite[Definition 2.30]{GeeMLT}), because we can mostly ``black box'' them.
\item In the compatible systems that we are interested in, it is likely that $V(w)$ can instead be defined over $F$. The definition that we are working with sidesteps this issue since it is irrelevant. Note that it at least implies that every $V(w)$ is isomorphic to any of its conjugates by $\mrm{Gal}(\overline{F}/F)$. 
\item Note that by \cite[(3.4.5), (4.1.6), (4.2.4)]{NTB} (and working prime by prime), a compatible system of Galois representations $V$ has a \emph{conductor ideal} $\mfk{f}_V \subset \OO_K$, which is a non-zero ideal of $\OO_K$. 
\end{itemize}
}

It is clear that compatible systems with a given support are stable under the usual operations of representation theory (direct sum, tensor product, dual, induction, restriction, extension of scalars, restriction of support). It is also clear that representations $G_K \rar \mrm{GL}_n(F)$ with open kernel (resp. the cyclotomic character), yield compatible systems of Galois representations of $G_K$ over $F$ (resp. over $\Q$).

To a compatible system $(V_v)_v=V$ of Galois representations of $G_K$ over a subfield $F$ of $\C$ we attach local $L$-factors as follows: for every finite place $w$ of $K$, the local factor $L_w(V,s)$ is the local $L$-factor $L((V(w) \otimes_{\overline{F}} \C)^{\ast},s)$ attached to complex Weil--Deligne representations\footnote{Because the $\overline{F}$-isomorphism class of $V(w)$ is invariant by $\mrm{Gal}(\overline{F}/F)$, the $\C$-isomorphism class of $V(w) \otimes_{\overline{F}} \C$ does not depend on the implicit choice of $F$-embedding $\overline{F} \rar \C$.} as in \cite[(4.1.6)]{NTB} and \cite[\S 29.3, \S 31.3]{BH}. The \emph{global $L$-function} of the compatible system $V$ is $L(V,s) = \prod_{w}{L_w(V,s)}$, where the product runs over all finite places.  

\defi{Let $K,F$ be number fields. A \emph{complete compatible system of Galois representations} of $G_K$ over $F$ is the datum of:
\begin{itemize}[noitemsep,label=\tiny$\bullet$]
\item a compatible system $V=(V_v)_v$ of Galois representations of $G_K$ over $F$, where every $V_v$ has dimension $d$ over $F_v$,
\item an embedding $F \rar \C$,
\item for every infinite place $w$ of $K$, a morphism $V(w): W_{K_v} \rar \mrm{GL}_d(\C)$. 
\end{itemize}}

It is also clear that complete compatible systems of Galois representations are stable under direct sum, tensor product, dual, induction, restriction, and extension of scalars. Representations $G_K \rar \mrm{GL}_n(F)$ with open kernel (resp. the cyclotomic character) clearly yield complete compatible systems of Galois representations of $G_K$ over $F$ (resp. of $G_{\Q}$ over $\Q$), by considering the representation $W_{F_v} \rar \mrm{Gal}(\overline{F_v}/F_v) \overset{\rho}{\rar} \mrm{GL}_d(\C)$ for $v$ real (resp. the representation $\omega_1: W_{\Q} \rar \C^{\ast}$ in the notation of \cite[(2.2)]{NTB}). 

Given a complete compatible system $V$ of Galois representations of $G_K$ over $F$, we denote its \emph{completed $L$-function} by $\Lambda(V,s) = \mathbf{N}(\mfk{f}_V)^{s/2}\prod_{w}{L(V(w)^{\ast},s)}L(V,s)$, where the product is over all the infinite places of $K$, and the $L$-factor is the one given in \cite[(3.3.1)]{NTB}. 

As an example, if $V$ is the complete compatible system arising from an Artin representation $\rho: G_K \rar \mrm{GL}_d(F)$ for some number field $F \subset \C$, then $L(V,s)$ (resp. $\Lambda(V,s)$) are the standard Artin $L$-function (resp. completed Artin $L$-function as in \cite[VII.12.2]{Neukirch-ANT}) of $\rho$. Similarly, in our definition, if $V$ is any complete compatible system, then the completed $L$-function attached to the complete compatible system $V(1)$ is $\Lambda(V(1),s)=\Lambda(V,s-1)$.

\defi{Let $f \in \mathcal{S}_k(\Gamma_1(N))$ be a newform with character $\chi$ and fix a number field $F \subset \C$ containing the Fourier coefficients of $F$. Let $\lambda$ be a finite place of $F$ with residue characteristic $\ell$, we call $V_{f,\lambda}: G_{\Q} \rar \GL{F_{\lambda}}$ the Galois representation attached to $f$ (see e.g. \cite[Theorems 2.1, 2.3]{Ribet-Antwerp5} for $k \geq 2$ and \cite[Th\'eor\`eme 4.1]{Del-Ser} for $k=1$), uniquely determined by the following properties:
\begin{itemize}[noitemsep,label=\tiny$\bullet$]
\item It is continuous, irreducible, and unramified outside $N\ell$,
\item For every prime $q \nmid N\ell$, the characteristic polynomial of an arithmetic Frobenius at $q$ is given by $X^2-a_q(f)+q^{k-1}\chi(q)$.  
\end{itemize}
We will also, by a slight abuse of notation, use $V_{f,\lambda}$ to denote the total space of this representation.  
}

The following result, which we will use throughout the article, is mainly a consequence of Carayol's theorem on local-global compatibility \cite[Th\'eor\`eme A]{Carayol-lg}, at least in weight $k \geq 2$ --- it follows in weight one from \cite[Th\'eor\`eme 4.6]{Del-Ser} and the fact that the newform is attached to an Artin representation.  

\prop[carayol-local-global]{Let $f \in \mathcal{S}_k(\Gamma_1(N))$ be a newform and let $F \subset \C$ be a number field containing its coefficients. The family $(V_{f,\lambda})_{\lambda}$ forms a compatible system $V_f$ of conductor ideal $N\Z$ with $L$-function $L(f,s) := \sum_{n \geq 1}{a_n(f)n^{-s}}$. 

If one attaches to $V_f$ the representation $W_{\R} \rar \GL{\C}$ given by $\mrm{Ind}_{W_{\C}}^{W_{\R}}{z^{k-1}}$ (where $z$ denotes either isomorphism $W_{\C}^{ab} \simeq \C^{\ast}$ as in \cite[(2.2.2)]{NTB}), then $V_f$ is a complete compatible system and one has $\Lambda(V_f,s) = 2N^{s/2}(2\pi)^{-s}\Gamma(s)L(f,s)$. 
}

\rem{Let $v$ be a place of $\Q$ and $f$ be a newform of weight $k \geq 2$, then the construction and Carayol's theorem (the latter when $v$ is finite; for $v$ infinite it is essentially by construction) imply that $V_f(v)^{\ast}$ is equal to the twist of the Weil--Deligne representation $\sigma_{f,v}$ attached by the local Langlands correspondence to $\pi_{f,v}$, up to a twist, i.e. $V_f(v)^{\ast} = \sigma_{f,v} \otimes \omega_{(1-k)/2}$ in the notation of \cite[(2.2)]{NTB}. While Proposition \ref{carayol-local-global} would hold for many choices of $V_f(\infty)$, this compatibility is precisely why we defined the representation $V_f(\infty)$ in this way. 

This claim is also true when $k=1$, but it is not a consequence of Carayol's result. Using the global equality of $L$-functions \cite[Th\'eor\`eme 4.6]{Del-Ser} as well as various well-known results regarding the behavior of local constants with respect to unramified and highly ramified twists (cf. \cite{AL78} for the modular aspects, and \cite[\S 3]{NTB}, \cite[\S 29.4]{BH} for the Galois side), this can easily be proved as a consequence of the local converse theorem \cite[27.1]{BH}. A full proof is contained in \cite[Proposition 4.3.11]{Studnia-thesis}. }

Given that the local Langlands correspondence is known to be compatible with pairs (see e.g. \cite[Th\'eor\`eme 1, \S 2.4]{Carayol-LL} or \cite[Theorem 12.4.1]{Getz-Hahn}), we can in fact apply Carayol's theorem to tensor products of complete compatible systems attached to modular forms. By the theory of the Rankin--Selberg products (cf. e.g. \cite[Theorem 11.7.1]{Getz-Hahn}), one has the following result: 

\prop[rankin-selberg-galois]{Let $f \in \mathcal{S}_k(\Gamma_1(N)), g \in \mathcal{S}_l(\Gamma_1(M))$ be two newforms with characters $\chi,\psi$. 

Then the $L$-function $\Lambda(f,g,s) := \Lambda(V_f \otimes V_g,s)$ is well-defined for $\mrm{Re}(s)$ large enough, and extends to a meromorphic function satisfying the functional equation 
\[\Lambda(f,g,s) = \prod_v{\varepsilon_v((V_f(v) \otimes V_g(v))^{\ast},\frac{k+l-1}{2},\theta_v)} \Lambda(\overline{f},\overline{g},k+l-1-s),\]
where the local constants $\varepsilon_v$ are those defined in \cite[(3.6.4)]{NTB}, with the obvious adaptation of \cite[(4.1.6)]{NTB} to Weil--Deligne representations, and $\prod_v{\theta_v}: \mathbb{A}_{\Q} \rar \C^{\times}$ is the additive character that vanishes on $\Q$ such that $\theta_{\infty}(z) = e^{-2i\pi z}$ for $z \in \R$. It is entire if $f \neq \overline{g}$. Moreover, the quotient
\[\frac{L(V_f \otimes V_g,s)}{L(\chi\psi,2s+2-k-l)\sum_{n \geq 1}{a_n(f)a_n(g)n^{-s}}}\] is a finite product of functions in $\C(p^{-s})$ for primes $p \mid MN$.  

} 

\rem{A construction of the Rankin--Selberg $L$-function for the product of two newforms, purely in terms of modular forms (i.e. that does not involve the correspondence between modular forms and either automorphic representations or Galois representations), and a proof of its functional equation, can be found in \cite{Li-RS}. As remarked in \emph{loc.cit.}, the $L$-function and the determination of local constants is place-by-place what appears in the derivation of the Rankin--Selberg functional equation by automorphic means, except in one specific case where, in \emph{loc.cit.}, some local $L$-factors contain information that is ``usually'' put in the local constant when working with automorphic representations. }

It is also not difficult to deduce from \cite[\S 4]{Ribet-Antwerp5}, \cite[\S II-2.5]{SerreMcGill} and some well-known facts about class field theory (e.g. \cite[Theorems 15.10, 15.11, Remark 15.12]{Harari}) that the characters attached to CM modular forms form a compatible system. We record a precise statement here for future reference. 

\prop[compatible-cm-systems]{Let $k \geq 2$ and $f \in \mathcal{S}_k(\Gamma_1(N))$ be a newform with character $\chi$ and complex multiplication by the imaginary quadratic field $K$. There exists a number field $F \subset \C$ containing the image of $K$ in $\C$, all the Fourier coefficients of $f$ and, for every prime $q \nmid N$ split in $K$, the complex roots of the Hecke polynomial $X^2-a_q(f)X+q^{k-1}\chi(q)$. 

Let $\sigma,\tau: K \rar F$ denote the two possible embeddings. There exist two compatible systems (in the sense of \cite[Definition p. I-11, Remark p. I-13]{SerreMcGill}) $(\psi_{\lambda,\sigma})_{\lambda}$ and $(\psi_{\lambda,\tau})_{\lambda}$ of characters of $G_K$ over $F$ such that for every maximal ideal $\lambda \subset \OO_F$:
\begin{itemize}[noitemsep,label=\tiny$\bullet$]
\item one has $(V_{f,\lambda})_{|G_K} \simeq \psi_{\lambda,\sigma} \oplus \psi_{\lambda,\tau}$,
\item $\psi_{\lambda,\sigma}$ and $\psi_{\lambda,\tau}$ are conjugate by $\mrm{Gal}(K/\Q)$,
\item if $\mfk{a}: \mathbb{A}_K^{\times}/K^{\times} \rar G_K^{ab}$ denotes the class field theory reciprocity map sending a uniformizer to an arithmetic Frobenius, for each $u \in \{\sigma,\tau\}$, $\psi_{\lambda,u}\circ\mfk{a}_{|K_{u^{-1}(\lambda)}}: K_{u^{-1}(\lambda)}^{\times} \rar F_{\lambda}^{\times}$ is given in a neighborhood of $1$ by $z \mapsto u(z)^{1-k}$. 
\end{itemize}
Moreover, for any $\lambda$, the character $\psi_{\lambda,\sigma}\psi_{\lambda,\tau}^{-1}$ has infinite image. 
}

\subsection{Cohomology of $X_{G_0}(p)$ and Tate module of $J_{G_0}(p)$}
\label{subsect:cohomology-general}

From now on and until the end of Section \ref{sect:XEP-general}, $p \geq 3$ is a prime. Recall that $G_0$ is the constant group scheme equal to $(\Z/p\Z)^{\oplus 2}$ over $R_0 := \Z[1/p]$. As we saw, the smooth projective $R_0$-scheme $X_{G_0}(p)$ of relative dimension one has a Jacobian $J_{G_0}(p)$ which is an abelian scheme over $R_0$ and is endowed with an action of $\mathbb{T}[\GL{\F_p}]$.

We call $\mathcal{N}$ the collection of newforms contained in $\mathcal{S}_2(\Gamma_1(p))$ or $\mathcal{S}_2(\Gamma_0(p^2) \cap \Gamma_1(p))$. 
 Let $\mathcal{D}$ be the group of Dirichlet characters modulo $p$.  

\lem[spcm-quadrichotomy]{The map $(\chi,f) \in \mathcal{D} \times \mathcal{N} \rar \mathcal{N}$ is a group action of $\mathcal{D}$ on $\mathcal{N}$. If $\omega \subset \mathcal{N}$ is one of its orbits, it satisfies exactly one of the following conditions:
\begin{itemize}[noitemsep,label=\tiny$\bullet$]
\item $\omega$ contains exactly one $f \in \mathcal{S}_2(\Gamma_1(p))$. One has $f \in \mathcal{S}_2(\Gamma_0(p))$, and $\mathcal{D}$ acts freely on $\omega$. 
\item $\omega$ contains exactly two elements $f_1,f_2 \in \mathcal{S}_2(\Gamma_1(p))$. One has $f_2=\overline{f_1}$, the character of $f_1$ is not trivial, and $\mathcal{D}$ acts freely on $\omega$.
\item $\omega$ does not meet $\mathcal{S}_2(\Gamma_1(p))$ and $\mathcal{D}$ acts freely on $\omega$. 
\item $\omega$ does not meet $\mathcal{S}_2(\Gamma_1(p))$ and the unique quadratic character $\lambda: \F_p^{\times} \rar \{\pm 1\}$ acts trivially on $\omega$. Then $\mathcal{D}/\langle\lambda\rangle$ acts freely on $\omega$, every $f \in \omega$ has complex multiplication by $\Q(\sqrt{-p})$; in particular, one has $p \equiv -1 \mod{4}$. 
\end{itemize}}

\demo{That the twist is well-defined follows from \cite[Proposition 3.1]{AL78}; it is easily seen to be a group action. Moreover, if $f \otimes \chi=f$ for $f \in \mathcal{N}$ and a non-trivial $\chi \in \mathcal{D}$, then comparing characters yields $\chi^2=1$, thus $\chi=\lambda$. Then $f$ has complex multiplication by $\lambda$ and one has $\lambda(-1)=-1$ by \cite[Theorem 4.5]{Ribet-Antwerp5}, whence $p \equiv -1 \mod{4}$. Because $\mathcal{D}$ is abelian, this implies that either $\mathcal{D}$ acts freely on $\omega$, or $\lambda$ acts trivially on $\omega$ and $\mathcal{D}/\langle\lambda\rangle$ acts freely on $\omega$. 

Let $f \in \mathcal{S}_2(\Gamma_1(p))$ have character $\alpha \in \mathcal{D}$ and let $\beta \in \mathcal{D}$ be non-trivial. If $\beta \neq \alpha^{-1}$, then \cite[Theorem 4.1]{AL78} implies that $f \otimes \beta = \sum_{n \geq 1}{a_n(f)\beta(n)q^n}$ is a newform of level $p^2$. If $\beta=\alpha^{-1}$, then $f \otimes \beta=\overline{f}$ by \cite[Proposition 1.2, (1.1)]{AL78}. In particular, if $f \otimes \beta=f$, then $\beta=\alpha^{-1}=\lambda$, so $\alpha(-1)=\lambda(-1)=-1$ (since $p \equiv -1\mod{4}$), which is impossible since $f$ has weight $2$. The conclusion follows. }

\defi{We let 
\begin{itemize}[noitemsep,label=\tiny$\bullet$]
\item $\mathscr{S}$ be the collection of newforms in $\mathcal{S}_2(\Gamma_0(p))$,
\item $\mathscr{P}$ be the collection of newforms in $\mathcal{S}_2(\Gamma_1(p))$ with non-trivial character; the complex conjugation acts on $\mathscr{P}$ as an involution without a fixed point,
\item $\mathscr{C}$ be the collection of newforms in $\mathcal{S}_2(\Gamma_0(p^2))$ which do not have a twist of level $p$; the twist by the quadratic character $\F_p^{\times} \rar \{\pm 1\}$ defines an involution of $\mathscr{C}$,
\item $\mathscr{C}_M$ be the collection of $f \in \mathscr{C}$ with complex multiplication: they are the fixed points in $\mathscr{C}$ of the twist by the unique quadratic Dirichlet character of conductor $p$. 
\end{itemize}
}

For $f \in \mathcal{N}$, we define the left representation $R(f)$ of $\GL{\F_p}$ by $R(f) = \pi_{f,p}^{1+p\mathcal{M}_2(\Z_p)}$. It is a complex irreducible representation of $\GL{\F_p}$ (cf. e.g. \cite[\S 6]{BH} for a brief reminder of the representation theory of $\GL{\F_p}$) and satisfies the following properties: 

\begin{itemize}[noitemsep,label=\tiny$\bullet$]
\item If $f \in \mathscr{S}$, then $R(f)$ is the Steinberg representation $\mrm{St}$. 
\item If $f \in \mathscr{P}$ has character $\chi$, then $R(f)$ is the principal series representation attached to the pair of characters $\mathbf{1},\chi$ (which we denote $\pi(1,\chi)$). 
\item If $f \in \mathscr{C}$, then $R(f)$ is a cuspidal representation of $\GL{\F_p}$, and is attached to a pair $\{\phi,\phi^p\}$ of complex characters of $\F_{p^2}^{\times}$ with $\phi \neq \phi^p$.
\item If $f \in \mathcal{N}$ and $\alpha \in \mathcal{D}$, then $R(f \otimes \alpha)\simeq R(f) \otimes \alpha(\det)$ (by Proposition \ref{compatibility-with-twists}).
\item If $f \in \mathcal{N}$ and $\sigma \in \mrm{Aut}(\C)$, then $R(\sigma(f)) \simeq \sigma(R(f))$ (by Proposition \ref{compatibility-with-galois}).  
\end{itemize}

By \cite[Cor. to Theorem 24]{SerreLinReps}, $R(f)$ can be defined over any subfield of $\C$ containing the $p(p^2-1)$-th roots of unity. In fact, if $R = \Z\left[\zeta_{p(p^2-1)}\right]$, it is easy to construct (e.g. by averaging some $R$-lattice in $R(f)$) a projective $R$-module $M$ such that $M \otimes_R \C$ is isomorphic to $R(f)$.  

\medskip

We now discuss the space of differential forms of $X_{G_0}(p)$. 

\prop[tglmu-xp]{The finite free $R_0$-module $H^0(X_{G_0}(p),\Omega^1_{X_{G_0}(p)/R_0})$ carries a right action of $\mathbb{T}[\GL{\F_p}]$ and a left action of $\OO(X_{G_0}(p)) \simeq R_0[\zeta_p]$. For $n \in \F_p^{\times}$, $t \in \mathbb{T}[\GL{\F_p}]_n$, $\varphi \in \OO(X_{G_0}(p))$, $\omega \in H^0(X_{G_0(p)}, \Omega^1_{X_{G_0}(p)/R_0})$, one has $(\varphi\omega) \mid t = (\underline{n}\varphi) \cdot (\omega \mid t)$, where $\underline{n} \in \Aut{\OO(X_{G_0}(p))}$ maps $\zeta_p$ to $\zeta_p^n$. }

\demo{It is enough to prove the result when $t \in \GL{\F_p}$ and $t=gT_{\ell}$ for $\ell \neq p$ prime and $\ell\det{g} \equiv 1 \pmod{p}$. The case of $\GL{\F_p}$ is completely formal. For $t=gT_{\ell}$, this follows from the fact that $D_{1,1}: X_{G_0}(p,\Gamma_0(\ell)) \rar X_{G_0}(p)$ and $g \circ D_{\ell,1}: X_{G_0}(p,\Gamma_0(\ell)) \rar X_{G_0}(p)$ are both morphisms of $\mrm{Pair}(G_0)\simeq \Sp{\OO(X_{G_0}(p))}$-schemes.}

One convenient way to manipulate the holomorphic differentials on $X_{G_0}(p)$ (which we will need in order to apply the Chabauty--Coleman method) is the \emph{$q$-expansion}. The standard references for such a result are \cite[VII.3]{DeRa} and \cite[\S 1.6, \S 1.11]{Katz}\footnote{Since we are only interested in the $q$-expansion of holomorphic differentials, one can give a more elementary construction that avoids the Kodaira--Spencer isomorphism used by Deligne--Rapoport and Katz.}\footnote{A crucial input is \cite[Chapters 8.11, 10.9]{KM}, which describes the cuspidal subscheme and its link to the Tate curve.}

Let $R := \Z[1/p,\zeta_p] = \OO(X_{G_0}(p))$ for brevity. The triple $\tau := (\mrm{Tate}(q),q^{1/p},\zeta_p^{-1})$ is a $R((q^{1/p}))$-point of $Y_{G_0}(p)$, where $\mrm{Tate}(q)$ denotes the Tate curve over $\Z((q))$ (see \cite[(A1.2.3)]{Katz} and \cite[(8.8)]{KM}\footnote{The value of $c_6$ claimed in \emph{loc.cit.} differs from that of \cite[(8.4)]{Deligne-formulaire} by a sign $-1$, so it is either wrong, or uses a different definition than \cite[Chapter III.1]{AEC1} and \cite[(1.5)]{Deligne-formulaire}. Indeed, when specializing to the case of $\Q_r$ with $q \in r\Z_r$, \cite[(8.8)]{KM} implies that the curve has split multiplicative reduction, so in the convention of \cite{AEC1} $\frac{-c_6(q)}{c_4(q)}$ has to be a square.}). The $\zeta_p^{-1}$ comes from the normalization of \cite[(2.8.5.2)]{KM} for the Weil pairing which we follow here. The integral closure of $R[[q]]$ in $R((q^{1/p}))$ is $R[[q^{1/p}]]$, so $\tau$ extends to a $\mu_p^{\times}$-morphism $\overline{\tau}: \Sp{R[[q^{1/p}]]} \rar X_{G_0}(p)$. We call the image under $\overline{\tau}$ of $\Sp{R} \rar \Sp{R[[q^{1/p}]]}$ the \emph{cusp at infinity}, and call it $\infty \in X_{G_0}(p)(R)$.

Let us recall the properties of $q$-expansion, which can be deduced from the works cited in the two previous paragraphs. 

\prop[properties-of-q-exp]{The cuspidal subscheme of $X := X_{G_0}(p)$ is a disjoint reunion of copies of $\Sp{R}$. For every section $c$ of the cuspidal subscheme of $X_{G_0}(p)$, there exists $\gamma \in \SL{\F_p}$ such that, for every regular excellent $R_0$-algebra $R'$, there exists a homomorphism of $q$-expansion at $c$
\[\mrm{Exp}_{c,R'}: H^0(X_{R'},\Omega^1_{X_{R'}/R'}) \rar q^{1/p}R'[[q^{1/p}]] \otimes R\] satisfying the following properties:  
\begin{itemize}[noitemsep,label=\tiny$\bullet$]
\item The formation of $\mrm{Exp}_{c,R'}$ commutes with base change (in $R'$), it is injective with flat cokernel,
\item For any $\omega \in H^0(X_{R'},\Omega^1_{X_{R'}/R'})$ with $\mrm{Exp}_{c,R'}(\omega)=f(q^{1/p})$ (for $f \in tR'[[t]] \otimes_{R_0} R$), one has an isomorphism of $R' \otimes R_0$-modules 
\[\frac{H^0((\mu_p^{\times})_{R'},c^{\ast}\Omega^1_{X_{R'}/{R'}})}{(R'\otimes R_0) \cdot c^{\ast}\omega} \simeq (R' \otimes R_0)/(f'(0)).\]
\item Let $(a,b) \in \F_p\times \F_p^{\times}$ and $M_{a,b} = \gamma\begin{pmatrix}1 & a\\0 & b\end{pmatrix}\gamma^{-1}$. Let $\omega \in H^0(X_{R'},\Omega^1_{X_{R'}/R'})$ and write $\mrm{Exp}_{c,R'}(\omega)=f(q^{1/p})$, then \[\mrm{Exp}_{c,R'}(\omega \mid M_{a,b})=(\underline{b}f)(\zeta_p^{-a}q^{1/p}),\] where $\underline{b}$ acts on $q^{1/p}R'[[q^{1/p}]] \otimes R$ exactly on the second component.  
\item For $R'=\C$, the image of $\mrm{Exp}_{c,R'}$ is exactly $\mathcal{S}_2(\Gamma(p)) \otimes_{\Z} R$. Let moreover $\mrm{Exp}^1_{c,R'}$ be the composition of $\mrm{Exp}_{c,R'}$ with the map $\C \otimes R \rar \C, \zeta_p \mapsto e^{-\frac{2i\pi}{p}}$. Then for $\omega \in H^0(X_{R'},\Omega^1_{X_{R'}/R'})$ and $g \in \SL{\F_p}$, one has $\mrm{Exp}^1_{c,R'}(\omega \mid \gamma g \gamma^{-1}) = \mrm{Exp}^1_{c,R'}(\omega) \mid g$.
\item When $c=\infty$, one has $\gamma=\mrm{id}_2$. For $\omega \in H^0(X_{R'},\Omega^1_{X_{R'}/R'})$, and $\ell \neq p$ prime, let \[\mrm{Exp}_{\infty,R'}\left(\omega\right) = \sum_{n \geq 1}{a_nq^{n/p}},\qquad \mrm{Exp}_{\infty,R'}\left(\omega \mid \begin{pmatrix}\ell & 0\\0 & 1\end{pmatrix}\right) = \sum_{n \geq 1}{b_nq^{n/p}},\] then one has \[\mrm{Exp}_{\infty,R'}\left(\omega \mid T_{\ell}\right) = \sum_{n \geq 1}{(b_{n\ell} + \ell \underline{\ell} a_{n/\ell})q^{n/p})}.\] In particular, then $R'=\C$, the action of $w(T_{\ell})$ through $\mrm{Exp}_{\infty,R'}$ is given by the double coset operator $\Gamma(p)\begin{pmatrix}1 & 0\\0 & \ell\end{pmatrix}\Gamma(p) \otimes \underline{\ell^{-1}}$. 
\end{itemize}
}

\rem{In the $q$-expansion at infinity that we just defined, over $R'=\C$, the pull-back by the Hecke operator $T_{\ell}$ acts by $\Gamma(p)\begin{pmatrix}\ell & 0\\0 & 1\end{pmatrix}\Gamma(p) \otimes \underline{\ell}$, that is, the Peterson dual of the classical Hecke operator appears. This is in contrast with the calculations in \cite{Studnia-thesis} or the discussion of \cite[p. 68]{Ribet-Stein}. The reason is the following: through our choice of uniformization (which is the $\Z[\zeta_p]$-linear one and remains compatible with the exponential notation $q=e^{2i\pi\tau}$), the image of the infinity cusp through the projection $X(p)^{an}_{\mrm{We}=e^{-\frac{2i\pi}{p}b}} \rar X_1(p)^{an}$ (sending $(E,P,Q)$ to $(E,Q)$) \emph{depends on $b$}. 
}

\prop[h10-xp]{Let $F \subset \C$ be a number field containing the Fourier coefficients of every $f \in \mathcal{N}$ and the $p(p-1)^2(p+1)$-th roots of unity. Then, as a right $\mathbb{T}[\GL{\F_p}] \otimes F$-module, $H^0(X_{G_0}(p),\Omega^1_{X_{G_0}(p)/R_0}) \otimes F$ is isomorphic to the direct sum of $M_f \otimes R(f)_F^{\vee}$ over $f \in \mathcal{N}$, where:
\begin{itemize}[noitemsep,label=\tiny$\bullet$]
\item $R(f)_F$ denotes a model of $R(f)$ over $F$,
\item $R(f)_F^{\vee} := \mrm{Hom}(R(f)_F,F)$ is a right absolutely irreducible $F[\GL{\F_p}]$-module,
\item $M_f$ is a right $\mathbb{T} \otimes F$-module where $T_n$ for $n \geq 1$ prime to $p$ (resp. $\diam{m}$ for $m \in \F_p^{\times}$) acts by $a_n(f)$ (resp. by $\chi(m)$ where $\chi$ is the character of $f$).  
\end{itemize}
}

\demo{It is enough to show that this decomposition holds when $F=\C$. This seems like a classical result on the coherent cohomology of modular curves, but we are not sure about a suitable reference. This can certainly be deduced from the results of \cite[\S 2]{Langlands-Antwerp2} using the Hodge decomposition of $H^1_{dR}(X_{G_0}(p)_{\C}^{an},\C)$ and Serre duality, along with the fact that the adjoint under Poincar\'e duality of the pull-back by $t \in \mathbb{T}$ (resp. $g \in \GL{\F_p}$) is the pull-back by $w(t)$ (resp. the pull-back by $g^{-1}$). A proof can also be found in \cite[Theorem 2.B]{Studnia-thesis}. }

In the rest of the section, we fix a number field $F \subset \C$ containing all $p(p^2-1)(p-1)$-th roots of unity and the Fourier coefficients of all newforms in $\mathcal{N}$.

\cor[t-structure]{For every $f \in \mathcal{N}$, there is a ring homomorphism $\mu_f: \mathbb{T} \rar \OO_F$ sending $T_n$ for $n \geq 1$ prime to $p$ (resp. $\diam{n}$ for $n \in \F_p^{\times}$) to $a_n(f)$ (resp. $\chi(n)$, where $\chi$ is the character of $f$). Moreover, $\prod_{f \in \mathcal{N}}{\mu_f}: \mathbb{T} \otimes F \rar F^{\mathcal{N}}$ is a $F$-linear isomorphism.}

\demo{Every $\mu_f$ exists by Proposition \ref{h10-xp} (since the Fourier coefficients of newforms are algebraic integers, cf. e.g. \cite[\S 1.2]{Katz} and \cite[VII.4]{DeRa}). Since $\mathbb{T} \otimes \Q \rar \mrm{End}(H^0(X_{G_0}(p)_{\Q},\Omega^1_{X_{G_0}(p)_{\Q}/\Q}))$ is injective, $\mathbb{T} \otimes F$ acts faithfully on $H^0(X_{G_0}(p)_{F},\Omega^1_{X_{G_0}(p)_F/F})$. Therefore, by Proposition \ref{h10-xp}, $\mathbb{T} \otimes F$ is reduced and the $F$-algebra homomorphism $\prod_{f \in \mathcal{N}}{\mu_f}: \mathbb{T} \otimes F \rar F^{\mathcal{N}}$ is injective. Since the $\mu_f$ are themselves pairwise distinct, the conclusion follows by elementary commutative algebra. }

\prop[untwisted-tate-module]{For every finite prime $\lambda$ of $F$ with residue characteristic $\ell$, $\Tate{\ell}{J_{G_0}(p)} \otimes_{\Z_{\ell}} F_{\lambda}$ is isomorphic, as a $\mathbb{T}[\GL{\F_p} \times G_{\Q}] \otimes F_{\lambda}$-module, to the direct sum of the $V_{f,\lambda} \otimes_{F[\F_p^{\times}]} R(f)$ over $f \in \mathcal{N}$, where: 
\begin{itemize}[noitemsep,label=\tiny$\bullet$]
\item $f$ runs through newforms contained in $\mathcal{N}$, 
\item the action of $\mathbb{T}[G_{\Q}]$ is only carried by $V_{f,\lambda}$, the Hecke operator $T_n$ (for $n \geq 1$ prime to $p$) acts by $a_n(f)$, and $n$ acts by $\chi(n)$ where $\chi$ is the character of $f$.  
\end{itemize}
}

\demo{This is a classical rephrasing of well-known work on the cohomology of modular curves, see e.g. \cite[\S 3]{Langlands-Antwerp2}.}

\subsection{Semi-direct products and a geometric action of $\GL{\F_p}$ on $X_{\Gamma}(p)$}

Let $\Gamma$ be a $p$-torsion group over $\Q$. We want to use Proposition \ref{untwisted-tate-module} to determine the Tate module of the Jacobian of $X_{\Gamma}(p)$, and, when $\Gamma$ admits a Weil pairing, its connected components (which are the $X_{\Gamma}^{\alpha}(p)$ over $\alpha \in \mrm{Pair}(\Gamma)(\Q)$). To do this, we need some information on the relationship between $X_{\Gamma}(p)$ and $X_{G_0}(p)$. The following is a reformulation of \cite[Propositions 5.2.4, 5.2.5]{Studnia-moduli}. 

\lem[j1-p-q]{Let $k$ be a field where $p$ is invertible, $k_s$ be its separable closure, and $\Gamma$ be a $p$-torsion group scheme over $k$. Choose a basis $(P,Q)$ of $\Gamma(k_s)$. For every square-free $m \geq 1$ which is prime to $p$, there is an isomorphism $j_{m,P,Q}: X_{\Gamma}(p,\Gamma_0(m))_{k_s} \rar X_{G_0}(p,\Gamma_0(m))_{k_s}$ of $\mathbb{P}^1_{k_s}$-schemes (through the $j$-invariant) whose moduli description on non-cuspidal $k_s$-points is $(E/k_s,\alpha,C) \mapsto (E/k_s,(\alpha(P),\alpha(Q)),C)$. 
Moreover, the collection of the $j_{m,P,Q}$ satisfies the following properties:
\begin{itemize}[noitemsep,label=\tiny$\bullet$]
\item $j_{m,P,Q}$ commutes with all degeneracy maps and the action of the $\diam{n}$ for $n \in \F_p^{\times}$,
\item If one has an isomorphism $\iota: \Gamma \rar \Gamma'$ of $p$-torsion groups, and $(P',Q') = (\iota(P),\iota(Q))$, then $j_{m,P',Q'} = j_{m,P,Q} \circ \tilde{\iota} $, where $\tilde{\iota}: X_{\Gamma'}(p,\Gamma_0(m))_k \rar X_{\Gamma}(p,\Gamma_0(m))_k$ is the isomorphism attached to $\iota$,
\item For any $z \in X_{\Gamma}(p,\Gamma_0(m))(k_s)$ and any $\sigma \in G_k$, one has $j_{m,P,Q}(\sigma(z)) = M\sigma[j_{m,P,Q}(z)]$, where $M$ is the inverse transpose of the matrix of the action of $\sigma \in \Aut{\Gamma(k_s)}$ in the basis $(P,Q)$ (in other words, one has $M\begin{pmatrix}\sigma(P)\\\sigma(Q)\end{pmatrix} = \begin{pmatrix}P \\Q\end{pmatrix}$).
\item If $(P', Q')$ is another basis of $\Gamma(k_s)$, then $j_{m,P',Q'} = M \circ j_{m,P,Q}$, where $M \in \GL{\F_p}$ is such that $M\begin{pmatrix}P \\ Q\end{pmatrix} =\begin{pmatrix}P'\\Q'\end{pmatrix}$.     
\end{itemize}}

In particular, in the situation of Lemma \ref{j1-p-q}, the $k_s$-scheme $X_{\Gamma}(p,\Gamma_0(m))_{k_s}$ carries an action of $\GL{\F_p}$ to which $j_{m,P,Q}$ commutes (this uniquely defines this action). This action of $\GL{\F_p}$ clearly depends on the choice of $(P,Q)$, even though we do not include it in the notation. 

In particular, the following diagram commutes for every $M \in \GL{\F_p}$:
\[
\begin{tikzcd}[ampersand replacement=\&]
X_{\Gamma}(p,\Gamma_0(m))_{k_s} \arrow{r}{\mrm{We}} \arrow{d}{M} \& \mrm{Pair}(\Gamma)_{k_s} \arrow{d}{\det{M}}\\
X_{\Gamma}(p,\Gamma_0(m))_{k_s} \arrow{r}{\mrm{We}} \& \mrm{Pair}(\Gamma)_{k_s}\\
\end{tikzcd}
\]

Thus, the push-forward of $j_{1,P,Q}$ is an isomorphism $j_{\ast,P,Q}: J_{\Gamma}(p)_{k_s} \rar J_{G_0}(p)_{k_s}$ of abelian varieties over $k_s$, it commutes with the action of $\GL{\F_p}$ and the action of $\mathbb{T}$. 

\defi{If $N,H$ are two groups, and if we fix an action $\ast$ of $H$ on $N$ by group automorphisms, recall that the \emph{semi-direct product} of $N$ by $H$, denoted $N \rtimes H$, is the group structure on the set $N \times H$ given by $(n,h) \cdot (n', h') = (n (h \ast n'),hh')$. The maps 
\[\alpha: n \in N \mapsto (n,e_H) \in N \rtimes H,\quad \alpha': h \in H \mapsto (e_N,h) \in N \rtimes H,\quad \beta: (n,h) \in N \rtimes H \mapsto h \in H\] are group homomorphisms, the sequence 
\[1 \rar N \overset{\alpha}{\rar} N \rtimes H \overset{\beta}{\rar} H \rar 1\] is exact and $\alpha'$ is a section of this exact sequence. 

If furthermore $N,H$ are topological groups and $\ast:N \times H \rar N$ is continuous, then when endowed with the product topology and the semi-direct group law previously defined, the set $N \times H$ is a topological group.   

If now $\ast$ is given by $h \ast n =  \varphi(h)n\varphi(h)^{-1}$ (for all $(n,h) \in N \times H$) for some group homomorphism $\varphi: H \rar N$, then we also write $N \rtimes_{\varphi} H$ for the semi-direct product $N \rtimes_{\ast} H$. In this case, one checks directly that $\tilde{\varphi}: (n,h) \in N \rtimes_{\ast} H \mapsto (n\varphi(h),h) \in N \times H$ is a group isomorphism, and that it is a homeomorphism if $N,H$ are topological groups and $\varphi$ is continuous. 
}

\prop[semidirect-gl-gal]{Let $k$ be a field of characteristic distinct from $p$ with separable closure $k_s$ and $\Gamma$ be a $p$-torsion group scheme over $k$. Let $P,Q \in \Gamma(k_s)$ form a basis of $\Gamma(k_s)$. Let $\rho: G_k \rar \GL{\F_p}$ be given by the matrix of the action of $G_k$ on $\Gamma(k_s)$ in the basis $(P,Q)$. Then $\GL{\F_p} \rtimes_{\rho^{\ast}} G_k$ acts on $X_{\Gamma}(p)(k_s)$ and $J_{\Gamma}(p)(k_s)$ (in the latter case $\mathbb{T}$-linearly) in such a way that, for $M \in \GL{\F_p}$ and $\sigma \in G_k$, $(M,\mrm{id}_{k_s})$ acts by the $\GL{\F_p}$-action (attached to $P,Q$) defined in Lemma \ref{j1-p-q} and $(\mrm{id}_2,\sigma)$ acts by the natural action of $\sigma$ on $X_{\Gamma}(p)(k_s)$ (or $J_{\Gamma}(p)(k_s)$). }

\demo{If $\sigma \in \mrm{Gal}(k_s/k)$ and $M \in \GL{\F_p}$, a direct calculation shows that $\sigma M$ acts on the non-cuspidal points of $X_{\Gamma}(p)(k_s)$ by $\rho^{\ast}(\sigma)M\rho^{\ast}(\sigma)^{-1} \sigma$. In particular, the two $k_s$-automorphisms $\Sp{\sigma^{-1}} \circ M \circ \Sp{\sigma}$ and $\rho^{\ast}(\sigma)M\rho^{\ast}(\sigma)^{-1}$ of $X_{\Gamma}(p)_{k_s}$ agree on the non-cuspidal $k_s$-points. Since $X_{\Gamma}(p)$ is smooth proper over $k$, these two automorphisms agree. Therefore, by push-forward functoriality, $\Sp{\sigma^{-1}} \circ M \circ \Sp{\sigma}$ and $\rho^{\ast}(\sigma)M\rho^{\ast}(\sigma)^{-1}$ agree as $\mathbb{T}$-automorphisms of $J_{\Gamma}(p)_{k_s}$, and the conclusion follows.  }

\subsection{The Tate module of the Jacobian of $X_{\Gamma}^{\alpha}(p)$}

\lem[decomposition-rf-cm]{Let $f \in \mathscr{C}_M$. Then $p \equiv -1\mod{4}$ and $R(f)$ is the cuspidal representation of $\PGL{\F_p}$ attached to the unique pair of characters $\F_{p^2}^{\times} \rar \C^{\times}$ of order $4$.} 

\demo{One has $p \equiv -1\mod{4}$ by Lemma \ref{spcm-quadrichotomy}. Let $\lambda$ be the unique non-trivial quadratic character of $\F_p^{\times}$ and $\lambda_2 = \lambda \circ N_{\F_{p^2}/\F_p}$. One has $R(f)=R(f \otimes \lambda) = R(f) \otimes \lambda(\det)$. Hence, if $\phi: \F_{p^2}^{\times} \rar \C^{\times}$ is one of the two characters attached to $R(f)$, one has $\lambda_2\phi \in \{\phi,\phi^p\}$: $\lambda_2$ is non-trivial, so $\lambda_2\phi=\phi^p$. Taking squares, it follows that $\phi^{2p-2}=1$. Since $f \in \mathcal{S}_2(\Gamma_0(p^2))$, $\phi$ is trivial on $\F_p^{\times}$, so $\phi^{p+1}=1$, hence $\phi^4=1$. Since $\phi \neq \phi^p$, $\phi$ has order exactly $4$. }

\medskip

We now fix a number field $F \subset \C$ containing the $p(p-1)(p^2-1)$-th roots of unity, the Fourier coefficients of every $f \in \mathcal{N}$, and, for every $f \in \mathscr{C}_M$ \footnote{As we saw, this condition is empty if $p \not\equiv -1 \mod{4}$}, the complex roots of $X^2-a_q(f)X+q$ for every prime $q \neq p$ which is a square modulo $p$.

\lem[specify-r1-r2]{Assume that $p \equiv -1 \mod{4}$. Let $R$ be the cuspidal representation of $\GL{\F_p}$ attached to the pair of characters of $\F_{p^2}^{\times}$ of order $4$. Then $R_{|\F_p^{\times}\SL{\F_p}}$ has exactly two irreducible components $R_!,R_?$, which are dual to each other and conjugate under $\GL{\F_p}/\F_p^{\times}\SL{\F_p}$. Given a non-trivial unipotent $U \in \SL{\F_p}$ and $\zeta \in \mu_p^{\times}(F)$, there is a unique $R_{U,\zeta} \in \{R_!,R_?\}$ such that $\zeta$ is an eigenvalue of $R_{U,\zeta}(U)$. If $\zeta' \in \mu_p^{\times}(F)$ and $U' \in \SL{\F_p}$ is another non-trivial unipotent, $R_{U,\zeta} \simeq R_{U',\zeta'}$ if and only if there exists $\gamma \in \GL{\F_p}$ such that $(U',\zeta')=(\gamma U \gamma^{-1}, \zeta^{\det{\gamma}})$. 
 In particular, one has $R_{U,\zeta} \simeq R_{U^r,\zeta'}$ if and only if $\zeta'=\zeta^{rn^2}$ for some $n \in \F_p^{\times}$. In particular, the eigenvalues of $R_{U,\zeta}(U)$ and $R_{U,\zeta^{-1}}(U)$ are simple and form a partition of $\mu_p^{\times}(F)$. }

\demo{This is a character calculation. The main inputs are that, if $U \in \SL{\F_p}$ is a non-trivial unipotent, one has $\det(X \cdot  \mrm{id}-U \mid R) = \frac{X^p-1}{X-1}$, and that $U$ is $\SL{\F_p}$-conjugate to $U^{n^2}$ for every $n \in \F_p^{\times}$. }

\medskip

We are now ready to describe the Tate module of the Jacobian of $X_{\Gamma}^{\alpha}(p)$ as a $G_{\Q}$-module. However, we need to keep track of slightly more information in order to construct suitable quotients of $J_{\Gamma}^{\alpha}(p)$ in cases of small Galois image. 
If $\mathscr{C}_M$ is not empty, we also write $K=\Q(\sqrt{-p})$. Recall that, for every $f \in \mathcal{N}$, we defined in Corollary \ref{t-structure} a ring homomorphism $\mu_f: \mathbb{T} \rar F$.

\prop[decomp-tate-twisted-general]{Fix an algebraic closure $\Qbar$ for $\Q$. Let $\Gamma$ be a $p$-torsion group over $\Q$ endowed with a Weil pairing $\alpha: \Gamma \times \Gamma \rar \mu_p$, and let $P,Q \in \Gamma(\Qbar)$ form a $\F_p$-basis.  
Let $\rho: G_{\Q} \rar \GL{\F_p}$ be the continuous representation giving the action of $G_{\Q}$ on $\Gamma(\Qbar)$ in the basis $(P,Q)$, and $\rho^{\ast}$ be its contragredient, so that $\det{\rho^{\ast}} = \omega_p^{-1}$ (where $\omega_p$ is the cyclotomic character modulo $p$).  

Let $G'_{\Q} := \GL{\F_p} \rtimes_{\rho^{\ast}} G_{\Q}$, so that $(\tilde{\rho},\pi): (g,\sigma) \in G'_{\Q} \mapsto (g\rho^{\ast}(\sigma),\sigma) \in \GL{\F_p} \times G_{\Q}$ is an isomorphism. By Lemma \ref{semidirect-gl-gal}, $G'_{\Q}$ acts on $\coprod_{q \in \F_p^{\times}}{X_{\Gamma}^{\alpha^q}(p)_{\Qbar}}$ and $\prod_{q \in \F_p^{\times}}{J_{\Gamma}^{\alpha^q}(p)_{\Qbar}}$. The subgroup $G'_{\Q,1} = \SL{\F_p} \rtimes_{\rho^{\ast}} G_{\Q}$ stabilizes every component $X_{\Gamma}^{\alpha^q}(p)(\Qbar)$ and every summand $J_{\Gamma}^{\alpha^q}(p)(\Qbar)$. Let $\mathbb{T}[G'_{\Q}]_1$ be the subalgebra of $\mathbb{T}[G'_{\Q}] := \mathbb{T} \otimes_{\Z[\F_p^{\times}]} \Z[\GL{\F_p} \times G_{\Q}]$ generated by the $T_{n} \cdot (g,\sigma)$ where $n\det{g} \equiv 1 \mod{p}$, where $n \geq 1$ is prime to $p$. 

If $\mathscr{C}_M$ is not empty, let $\sigma: K \rar F$ be an embedding and $\tilde{\sigma}: \Qbar \rar \C$ extending $\sigma$, so that $\tilde{\sigma}(\alpha(P,Q)) \in \mu_p^{\times}(F)$. 
 The set $\sigma(\alpha(P,Q)) := \{\zeta^n \mid n \in \F_p^{\times 2}\}$ does not depend on $\tilde{\sigma}$ (thus justifying the notation), and, if $\tau: K \rar F$ is the other embedding, then $\sigma(\alpha(P,Q)),\tau(\alpha(P,Q))$ form a partition of $\mu_p^{\times}(F)$. 

Let $\lambda \subset \OO_F$ be a maximal ideal. Then, the direct sum of the following $\mathbb{T}[G'_{\Q}]_1 \otimes F_{\lambda}$-modules is isomorphic to $\Tate{\ell}{J_{\Gamma}^{\alpha}(p)} \otimes_{\Z_{\ell}} F_{\lambda}$: 

\begin{itemize}[noitemsep,label=\tiny$\bullet$]
\item for every $f \in \mathscr{S}$, a copy of $(V_{f,\lambda} \circ \pi) \otimes [\mrm{St}\circ \tilde{\rho}]$, where $\mathbb{T}$ acts through $\mu_f$, 
\item for every $f \in \mathscr{P}$ modulo complex conjugation, a copy of $(V_{f,\lambda} \circ \pi)\otimes [\pi(\chi,1) \circ \tilde{\rho}]$, where $\chi$ is the character of $f$ and $\mathbb{T}$ acts through $\mu_f$,
\item for every $f \in \mathscr{C} \backslash \mathscr{C}_M$ modulo quadratic twists, a copy of $V_{f,\lambda} \otimes [R(f) \circ \tilde{\rho}]$ where $\mathbb{T}$ acts through $\mu_f$,
\item for every $f \in \mathscr{C}_M$, a copy of \[\mrm{Ind}_{\F_p^{\times}\SL{\F_p} \rtimes_{\rho^{\ast}} G_K}^{\F_p^{\times}\SL{\F_p} \rtimes_{\rho^{\ast}} G_{\Q}}{(\psi_{\lambda,\sigma}\circ \pi) \otimes R_{U^{\varepsilon_f},\sigma(\alpha(P,Q))} \circ \tilde{\rho}},\]
 where $U=\begin{pmatrix}1 & 1\\0 & 1\end{pmatrix}$, $\varepsilon_f \in \{\pm 1\}$ is a canonical sign only depending on $f$, $\bigoplus_{n \in \F_p^{\times 2}}{\mathbb{T}_n}$ acts through $\mu_f$, while $T_ng$ acts by $0$ when $n\det{g}\equiv 1\mod{p}$ and $\det{g} \in \F_p^{\times}$ is not a square.
\end{itemize}

}

\rem{We formulated the result in this generality for three reasons:
\begin{itemize}[noitemsep,label=\tiny$\bullet$]
\item The first and main one is that, in order to prove our main result, we will need to consider quotients of the Jacobian $J_{\Gamma}^{\alpha}(p)$ not only by ideals of the Hecke algebra, as is the case in many other works using Mazur's strategy \cite{FreyMazur,Merel,Darmon-Merel,Parent-Torsion}, but also using the $\GL{\F_p}$-action. It is actually crucial that we are able to use elements of $\GL{\F_p}$ with non-square determinant\footnote{If $N \leq \GL{\F_p}$ is the normalizer of a non-split Cartan subgroup and $\mrm{St}$ is the Steinberg representation of $\GL{\F_p}$, then $\mrm{St}_{|N \cap \F_p^{\times}\SL{\F_p}}$ contains multiple copies of the same representation, but they ``correspond'' in quotients of $J_{\Gamma}^{\alpha}(p)$ to different representations of $G_{\Q}$.}. 
\item The second reason is that, while it makes the notation and proof more cumbersome, it does not introduce any new difficulty; in fact, it makes the argument simpler by making all the components irreducible, which is not always the case if we forget the action of $\GL{\F_p}$ and consider non-surjective Galois representations $\rho$. 
\item Finally, this decomposition is completely canonical and reveals the sign $\varepsilon_f$. Is there a way to tell whether it is $\pm 1$? 
\end{itemize}}

\rem{Using Part 4 of the proof, it is not difficult to adapt this result to arbitrary fields $k$ of definition of characteristic $0$, then by specialization to arbitrary fields $k$ of characteristic not $p$. The only difference in the statement is that, when $\omega_p(G_k)\subset \F_p^{\times}$ is contained in $\F_p^{\times 2}$ (i.e. if $k$ contains the algebraic integer $\frac{1+\sqrt{-p}}{2}$ when $p \equiv -1 \mod{4}$), the definition of the CM factor must be changed --- it is no longer induced, but is instead a direct sum.

Over the base field $\C$, what $\varepsilon_f$ is does not change the decomposition of the Tate module. However, a similar phenomenon can be seen at the level of the Hodge filtration of $H^1(X_{\Gamma}^{\alpha}(p)_{\C},\C)$. Indeed, for $f \in \mathscr{C}_M$, the $f$-isotypic component in the space $H^0(J_{\Gamma}^{\alpha}(p)_{\C},\Omega^1)$ is one of the two representations of $\SL{\F_p}$ constructed in Lemma \ref{decomposition-rf-cm}, and, for another choice of Weil pairing $\beta$, it will change if and only if $\beta=\alpha^n$ for some $n \in \F_p^{\times}$ which is not a square. In other words, $H^1(X_{\Gamma}^{\alpha}(p)_{\C},\C)$ as a $\mathbb{T}_1[\SL{\F_p}] \otimes \C$-module does not depend on $\alpha$, but its Hodge filtration does. }  

Before discussing the proof of Proposition \ref{decomp-tate-twisted-general}, let us state the following perhaps more evocative Corollary, which slightly refines a result of Virdol \cite[Theorem 1.1]{Virdol}. 

\cor[decomp-tate-twisted-gq]{Let $p \geq 7$ be a prime. Let $\Gamma$ be a $p$-torsion group over $\Q$ endowed with a Weil pairing $\alpha: \Gamma \times \Gamma \rar \mu_p$, and let $P,Q \in \Gamma(\Qbar)$ form a $\F_p$-basis of $\Gamma(\Qbar)$. 

Let $\rho: G_{\Q}  \rar \GL{\F_p}$ be the continuous representation such that for any $g \in G_{\Q}$, $\rho(g)$ is the matrix of the action of $g$ on $\Gamma(\Qbar)$ in the basis $(P,Q)$. Let $\rho^{\ast}$ be its contragredient, so that $\det{\rho^{\ast}} = \omega_p^{-1}$.

If $p \equiv -1 \mod{4}$, let $K = \Q(\sqrt{-p})$. Let $\lambda \subset \OO_F$ be a maximal ideal. Then, as a $\mathbb{T}_1[G_{\Q}] \otimes F_{\lambda}$-module, $\Tate{\ell}{J_{\Gamma}^{\alpha}(p)} \otimes_{\Z_{\ell}} F_{\lambda}$ is the direct sum of the following:

\begin{itemize}[noitemsep,label=\tiny$\bullet$]
\item for every $f \in \mathscr{S}$, a copy of $V_{f,\lambda} \otimes [\mrm{St}\circ \rho^{\ast}]$, 
\item for every $f \in \mathscr{P}$ with character $\chi$ modulo complex conjugation, a copy of $V_{f,\lambda} \otimes [\pi(\chi,1) \circ \rho^{\ast}]$,
\item for every $f \in \mathscr{C} \backslash \mathscr{C}_M$ modulo quadratic twists, a copy of $V_{f,\lambda} \otimes [R(f) \circ \rho^{\ast}]$,
\item for every $f \in \mathscr{C}_M$, if $R'$ denotes one of the representations of $\F_p^{\times}\SL{\F_p}$ described in Lemma \ref{specify-r1-r2}, a copy of $\mrm{Ind}_{K}^{\Q}{\psi_{\lambda,\sigma} \otimes [R' \circ \rho^{\ast}_{|G_K}]}$ for one choice of embedding $\sigma: K \rar F$. Exactly one such choice is valid if $\rho$ is surjective.   
\end{itemize}

The algebra $\mathbb{T}_1$ acts of any of the above summands through the character $\mu_f$.
} 

\rems{Let $\Gamma$ be a $p$-torsion group over $\Q$, endowed with a Weil pairing $\alpha$.

\begin{itemize}[noitemsep,label=\tiny$\bullet$]
\item It can be checked that all the $F_{\lambda}$-linear representations of $G_{\Q}$ described above are self-dual. Their Galois-theoretic root number is therefore equal to $\pm 1$, and is computed in \cite[Chapter 4]{Studnia-thesis}. Apart from the representations attached to $f \in \mathscr{C}_M$, the root number only depends on $p \mod{4}$, the restriction of $\rho$ to an inertia subgroup at $p$, and local invariants at $p$ of the modular form $f$. 
\item In the ``typical'' case where $\rho$ is surjective, if $L/\Q$ is the extension cut out by $\rho$, $L$ not solvable over its maximal totally real subfield (which is actually contained in $\Q(\zeta_p)$). Thus it seems that proving even meromorphicity and/or the functional equation for the complex $L$-functions attached to these Galois representations is beyond current techniques when $f \notin \mathscr{C}_M$. On the other hand, when $\rho$ has smaller image or in the case of the summands attached to $f \in \mathscr{C}_M$, more is known by work of Virdol \cite[\S 3]{Virdol}.  
\item Let $n \in \F_p^{\times}$ be a non-square element. The curves $X_{\Gamma}^{\alpha}(p)$ and $X_{\Gamma}^{\alpha^n}(p)$ are not \emph{a priori} isomorphic, but the only difference between the Tate modules of their Jacobians resides in their summands attached to $f \in \mathscr{C}_M$. In particular, if $p \equiv 1\mod{4}$, $J_{\Gamma}^{\alpha}(p)$ and $J_{\Gamma}^{\alpha^n}(p)$ are isogenous. This is not so surprising: if $\mathscr{C}_M = \emptyset$, it is not difficult to show that there exists $t \in \mathbb{T}_r$ (with $r \in \F_p^{\times}$ a quadratic non-residue) not contained in the kernel of any $\mu_f$, and $t: J_{\Gamma}^{\alpha}(p) \rar J_{\Gamma}^{\alpha^r}(p)$ is an isogeny. 
\end{itemize}}

\rem{Let $\Gamma$ be the $p$-torsion of some elliptic curve $E/\Q$ and let $f \in \mathscr{C}_M$ (so that $1$ and $-1$ are representatives for $\F_p^{\times}$ modulo its squares). Then $\Gamma$ has a canonical Weil pairing $W$ (the classical Weil pairing of $E$). Let $M_{f,E}^{\pm}$ be the $f$-isotypic components of the Tate module of $J_E^{\pm 1}(p) := _{\Gamma}^{W^{\pm 1}}(p)$. Using Corollary \ref{decomp-tate-twisted-gq}, one can compute in practice the unordered pair $\mathcal{L}_{unord}(E) := \{L(M_{f,E}^{\pm 1},s)\}$ (at least away from the Euler factors at $p$ and the primes of bad reduction for $E$). If $\varepsilon_f$ is known, one can actually compute the ordered pair $\mathcal{L}_{ord}(E) := (L(M_{f,E}^{+},s), L(M_{f,E}^{-},s))$. This can be used to determine $\varepsilon_f$ when $p=7,11$: given an elliptic curve $E_0/\Q$, equations for $X_{E_0}^{\pm 1}(p)$ are known (cf. \cite{Fisher-711}), so we can compute the first coefficients of $\mathcal{L}_{ord}(E_0)$ by point-counting. For suitable $E_0/\Q$, knowing a small truncation of $\mathcal{L}_{ord}(E_0)$ and $\mathcal{L}_{unord}(E_0)$ suffices to determine $\mathcal{L}_{ord}(E_0)$ and thus deduce $\varepsilon_f$. }

\rem{Let us keep the same notations as the previous remark. Assume furthermore that $\rho$ is wildly ramified at $p$. Then the $M_{f,E}^{\pm 1}$ have the same local root number at every place but $p$, but their local root numbers at $p$ are distinct, and we do not currently know a practical method to determine this root number. 
Relatedly, even in the case $p=7$ (when $\mathscr{S}=\mathscr{P}=\mathscr{C} \backslash \mathscr{C}_M=\emptyset$), we do not know of a way to determine the sign of the conjectural functional equation for the curve $X_{\Gamma}^{\mrm{We}^a}(p)$ when $\rho$ is wildly ramified at $p$.}

\medskip

We are now ready to prove Proposition \ref{decomp-tate-twisted-general}.

\demo{This argument is somewhat a sketch --- a complete proof, albeit using a slightly more complex route, can be found in the author's PhD dissertation \cite[Chapter 3]{Studnia-thesis}. 

The semi-direct product $\GL{\F_p} \rtimes_{\rho^{\ast}} G_{\Q}$ acts on $X_{\Gamma}(p)(\Qbar)$ and $J_{\Gamma}(p)(\Qbar)$ by Proposition \ref{semidirect-gl-gal}, the latter by $\mathbb{T}$-linear automorphisms. Under this action, for $n \in \F_p^{\times}$, $\begin{pmatrix}n & 0\\0 & 1\end{pmatrix}$ is an isomorphism from $X_{\Gamma}^{\alpha}(p)_{\Qbar}$ to $X_{\Gamma}^{\alpha^n}(p)_{\Qbar}$, and its push-forward $J_{\Gamma}^{\alpha}(p)(\Qbar) \rar J_{\Gamma}^{\alpha^n}(p)(\Qbar)$ is a $\mathbb{T}_1$-linear isomorphism. Therefore,
\begin{enumerate}[noitemsep,label=(\roman*)]
\item \label{dttc-6} the $\mathbb{T}_1 \otimes F_{\lambda}$-module structure of $\Tate{\ell}{J_{\Gamma}^{\alpha^n}(p)} \otimes_{\Z_{\ell}} F_{\lambda}$ does not depend on $n$. 
\end{enumerate}

\emph{Part 1: Modules, algebras and modular forms. }

Let $\mathbb{T}_{\rho^{\ast}} = \mathbb{T} \otimes_{\Z[\F_p^{\times}]} \Z[\GL{\F_p} \rtimes_{\rho^{\ast}} G_{\Q}]$ and let $\mathbb{T}_{\rho^{\ast},1}$ be its abelian subgroup generated by the $T_m (g,\gamma)$ for $m \geq 1, g \in \GL{\F_p}, \gamma \in \GL{\F_p}$ such that $m\det{g} \equiv 1 \mod{p}$. One checks easily that $\mathbb{T}_{\rho^{\ast},1}$ is a subring of $\mathbb{T}_{\rho^{\ast}}$. We also write $\mathbb{T}_{\rho^{\ast},1,1} = \mathbb{T}_1 \otimes_{\Z[\pm \mrm{id}_2]} \Z[\SL{\F_p} \rtimes_{\rho^{\ast}} G_{\Q}]$, which is a subring of $\mathbb{T}_{\rho^{\ast},1}$. 

For $f \in \mathcal{N}$, we define $M_f$ (we may add $\Gamma,P,Q,\rho,\lambda,\alpha$ as additional indices if they are required to make our claims unambiguous) to be the $\mathbb{T}_{\rho^{\ast}} \otimes F_{\lambda}$-module $(V_{f,\lambda}\circ \pi) \otimes (R(f)\circ \tilde{\rho})$, where $\mathbb{T}$ acts through $\mu_f$. This definition is valid because the central character of $R(f)$ is exactly $\chi$. 

Assuming $p \equiv -1 \mod{4}$, for $f \in \mathscr{C}_M$, $U \in \GL{\F_p}$ non-trivial unipotent, and $\sigma: K \rar F$, $\tilde{\sigma}: \Qbar \rar \C$ compatible field embeddings, let $\psi_{\lambda,\sigma}: G_K \rar F_{\lambda}^{\times}$ be the character attached to $f$. We define 

\begin{align*}
G'_K := \F_p^{\times}\SL{\F_p} \rtimes_{\rho^{\ast}} G_K, \qquad G'_{\Q,2} := \F_p^{\times}\SL{\F_p} \rtimes_{\rho^{\ast}} G_{\Q},\\ 
M_{f,U,\alpha} := \mrm{Ind}_{G'_K}^{G'_{\Q,2}}{\left\{[\psi_{\lambda,\sigma}\circ \pi]\otimes [R_{U,\tilde{\sigma}(\alpha(P,Q))} \circ \tilde{\rho}_{|G'_K}]\right\}},
\end{align*}
where $M_{f,U,\alpha}$ is a $\mathbb{T} \otimes_{\Z[\F_p^{\times}]} \Z[G'_{\Q,2}] \otimes F_{\lambda}$-module, and $\mathbb{T}$ acts through $\mu_f$. This is meaningful because $R_{U,\tilde{\sigma}(\alpha(P,Q))}$ has trivial central character. 

It is not difficult to see that $M_{f,U,\alpha}$ does not depend on $\sigma$, thus justifying the notation, and that $M_{f,U,\alpha} \simeq M_{f,U',\alpha^n}$ if $U'$ is $\SL{\F_p}$-conjugate to $U^n$. For $f \in \mathscr{C}_M$, $M_{f,U,\alpha}$ can be given a $\mathbb{T}_{\rho^{\ast},1} \otimes F_{\lambda}$-module structure compatible with the previous one since $a_n(f)=0$ for every $n \geq 1$ which is not a square modulo $p$ (i.e. for any such $n$, $T_n(g,\gamma)$ always acts by $0$).

Note the following facts; we only give sketches of their proof. 

\begin{enumerate}[noitemsep,label=(\roman*)]
\setcounter{enumi}{1}
\item \label{dttc-1} For $f \in \mathcal{N}$, the $\mathbb{T}_{\rho^{\ast},1} \otimes F_{\lambda}$-module structure of $M_f$ only depends on $f$ modulo twists by Dirichlet characters modulo $p$. Indeed, this amounts to showing that for $f \in \mathcal{N}$ and $\chi$ Dirichlet character modulo $p$, one has an isomorphism 
\[V_{f \otimes \chi,\lambda} \otimes R(f \otimes \chi) \simeq V_{f,\lambda} \otimes R(f) \otimes \chi(\det \cdot\, \omega_p)\] of $F_{\lambda}[\GL{\F_p} \times G_{\Q}]$-modules. 

\item \label{dttc-2} For $f \in \mathcal{N}$, $\mathbb{T}_1 \otimes F_{\lambda}$ acts on $M_f$ through the $F_{\lambda}$-algebra homomorphism $\mu_{f,1}: \mathbb{T}_1 \otimes F_{\lambda} \rar F_{\lambda}$; for $f,f' \in \mathcal{N}$, $\mu_{f,1} = \mu_{f',1}$ if and only if $f$ and $f'$ are twists of each other by a Dirichlet character modulo $p$. Indeed, the reverse direction follows from \ref{dttc-1}. To prove the forward direction, we need to show that if $f,f' \in \mathcal{N}$ are such that $a_{\ell}(f)=a_{\ell}(f')$ for every prime $\ell \equiv 1 \mod{p}$, then $f'$ is a twist of $f$ by a Dirichlet character modulo $p$. The condition implies that the semi-simple representations $(V_{f,\lambda})_{|\Q(\mu_p)}$ and $(V_{f',\lambda})_{|\Q(\mu_p)}$ are isomorphic. In particular, for some Dirichlet character $\chi$, the $\chi\circ\omega_p$-eigenspace of $\mrm{Hom}(V_{f,\lambda},V_{f',\lambda})$ is contains a non-zero element $\varphi$: then $\varphi$ induces a non-zero $F_{\lambda}[G_{\Q}]$-homomorphism $\varphi_1: V_{f,\lambda} \otimes (\chi\circ\omega_p) \rar V_{f',\lambda}$. Both sides are absolutely irreducible, so $\varphi_1$ is an isomorphism and $f \otimes\chi=f'$.

\item \label{dttc-3} For $f \in \mathcal{N}$ without complex multiplication, $M_f$ is irreducible as a $\mathbb{T}_{\rho^{\ast},1} \otimes F_{\lambda}$-module. Indeed, by \cite[Cor. 2 au Th. 15]{Serre-Chebotarev} and the prime number theorem for arithmetic progressions \cite[XV, Theorems 4 and 5]{Lang-ANT}, there are primes $q$ with arbitrary non-zero residue modulo $p$ such that $a_q(f) \neq 0$. As a consequence, this claim reduces to showing that $V_{f,\lambda} \otimes R(f)$ is an irreducible $F_{\lambda}[G_{\Q} \times \GL{\F_p}]$-module, which is an easy application of e.g. \cite[Lemma 2.3]{Studnia-Euler}.
 
\item \label{dttc-4} Assuming $p \equiv -1 \mod{4}$, for $f \in \mathscr{C}_M$, for any $U \in \GL{\F_p}$ unipotent and any $\alpha$, one has $M_f \simeq M_{f,U,\alpha} \oplus M_{f,U,\alpha^{-1}}$ as $\mathbb{T}_{\rho^{\ast},1} \otimes F_{\lambda}$-modules. This is a direct calculation using Proposition \ref{compatible-cm-systems} and Lemma \ref{decomposition-rf-cm}.  

\item \label{dttc-5} Let $f,f' \in \mathscr{C}_M$, $\sigma,\sigma' : K \rar F$ be embeddings, $\zeta,\zeta' \in F$ be primitive $p$-th roots of unity and $U \in \GL{\F_p}$ be a non-trivial unipotent element. Then the two $F_{\lambda}[\SL{\F_p} \rtimes_{\rho^{\ast}} G_K]$-modules $[\psi_{f,\lambda,\sigma}\circ \pi]\otimes [R_{U,\zeta} \circ \tilde{\rho}_{|G'_K}]$ and $[\psi_{f',\lambda,\sigma'}\circ \pi]\otimes [R_{U,\zeta'} \circ \tilde{\rho}_{|G'_K}]$ are isomorphic if and only if $(f',\sigma',\zeta')=(f,\sigma,\zeta^n)$ for some square $n \in \F_p^{\times}$. Indeed, we first show that $R_{U,\zeta}=R_{U,\zeta'}$ by computing the trace of $(U,\mrm{id})$, then, for each $\gamma \in G_K$, let $\delta$ be a square root of $\det{\rho^{\ast}(\gamma)}$, then taking the trace of $(\delta^{-1}\rho^{\ast}(\gamma)^{-1},\gamma)$ shows that $\psi_{f,\lambda,\sigma}(\gamma)=\psi_{f',\lambda,\sigma'}(\gamma)$. This implies by Proposition \ref{compatible-cm-systems} that $\sigma=\sigma'$, and by taking induced representations to $G_{\Q}$ that $f'=f$.    

\end{enumerate} 

Note that \ref{dttc-1} implies in particular that whether Proposition \ref{decomp-tate-twisted-general} holds does not depend on the choices of representatives of $\mathscr{P}$ modulo complex conjugation or $\mathscr{C}$ modulo quadratic twists. 

Let $\mathcal{N}'$ denote the reunion of $\mathscr{S}$, a set of representatives for $\mathscr{P}$ modulo complex conjugation, and a set of representatives for $\mathscr{C} \backslash \mathscr{C}_M$ modulo the quadratic twist. 

Furthermore, fix a non-trivial unipotent $U \in \SL{\F_p}$, then \ref{dttc-3}, \ref{dttc-5} imply that, as $\mathbb{T}_{\rho^{\ast},1} \otimes F_{\lambda}$-modules, the $M_f$ (over $f \in \mathcal{N}'$) and $M_{f,U^{\varepsilon},\alpha}$ (over $f \in \mathscr{C}_M$ and $\varepsilon \in \{\pm 1\}$) are absolutely irreducible and pairwise non-isomorphic.

\emph{Part 2: Fix $\Gamma,P,Q,\lambda$. There is a choice of $(\varepsilon_f)_{f \in \mathscr{C}_M}$ making the claim true for every $\alpha \in \mrm{Pair}(\Gamma)(\Q)$. }

Fix a Weil pairing $\alpha$ on $\Gamma$. 

Let, for every $n \in \F_p^{\times}$, $J_n = \Tate{\ell}{J_{\Gamma}^{\alpha^n}(p)} \otimes_{\Z_{\ell}} F_{\lambda}$ as a $\mathbb{T}_{\rho^{\ast},1} \otimes F_{\lambda}$-module, and we consider $J_{tot} = \bigoplus_{m \in \F_p^{\times}}{J_m}$. We fix a non-trivial unipotent $U \in \SL{\F_p}$. Note that the multiplication by $n \in \F_p^{\times}$ induces isomorphisms $J_m \rar J_{mn^2}$ of $\mathbb{T}_{\rho^{\ast},1}$-modules for each $m \in \F_p^{\times}$. Let $U=\begin{pmatrix}1 & 1\\0 & 1\end{pmatrix}$. 

As a consequence of Proposition \ref{untwisted-tate-module} and Lemma \ref{j1-p-q}, one has the following isomorphism of $\mathbb{T}[\GL{\F_p} \times G_{\Q}] \otimes F_{\lambda}$-modules: 
\[J_{tot} = \Tate{\ell}{J_{\Gamma}(p)} \otimes_{\Z_{\ell}} F_{\lambda} \simeq \bigoplus_{f \in \mathcal{N}}{M_f}.\]

Hence \ref{dttc-1}, \ref{dttc-2}, \ref{dttc-4} and Lemma \ref{spcm-quadrichotomy} imply that one has the following isomorphism of $\mathbb{T}_{\rho^{\ast},1} \otimes F_{\lambda}$-modules:

\[J_{tot} \simeq \bigoplus_{f \in \mathcal{N}'}{M_f^{\oplus (p-1)}} \oplus \bigoplus_{f \in \mathscr{C}_M}{M_{f,U,\alpha}^{\oplus \frac{p-1}{2}} \oplus M_{f,U^{-1},\alpha}^{\oplus \frac{p-1}{2}}}.\]

In particular, \ref{dttc-3},\ref{dttc-5} imply that $J_{tot}$ and therefore every $J_n$ is a semi-simple $\mathbb{T}_{\rho^{\ast},1} \otimes F_{\lambda}$-module, and their irreducible components are of the form $M_f$ for $f \in \mathcal{N}'$ or $M_{f,U^{\pm 1},\alpha}$ for $f \in \mathscr{C}_M$. 

Let now $f \in \mathcal{N}'$ and $d_n$ be the number of irreducible components of $J_n$ that are of the form $M_f$. By \ref{dttc-2}, \ref{dttc-6}, $d_n$ is constant, and $\sum_{n \in \F_p^{\times}}{d_n}$ is the number of $M_f$ appearing in $J_{tot}$, which is $p-1$, so every $d_n$ equals $1$. 

If $\mathscr{C}_M$ is empty, we are done; otherwise, let $f \in \mathscr{C}_M$ (so $p \equiv -1\mod{4}$ and $\pm 1$ are representatives for $\F_p^{\times}$ modulo squares). For every $n \in \F_p^{\times}$ and $\varepsilon \in \{\pm 1\}$, let $d_{n,\varepsilon}$ be the number of times that $M_{f,U^{\varepsilon},\alpha}$ appears as a part of $J_n$. Now, $\sum_{n \in \F_p^{\times}}{d_{n,\varepsilon}}$ is the number of times that a component of the form $M_{f,U^{\varepsilon},\alpha}$ appears in $J_{tot}$, so it is equal to $\frac{p-1}{2}$. By \ref{dttc-6}, $d_{n,1}+d_{n,-1}$ does not depend on $n$, so $d_{n,+1}+d_{n,-1}=1$ for every $n \in \F_p^{\times}$ (so $\{d_{n,\pm 1}\}=\{0,1\}$). Finally, the isomorphism class of $J_n$ as a $\mathbb{T}_{\rho^{\ast},1} \otimes F_{\lambda}$-module only depends on $n$ modulo squares, so $d_{n,\varepsilon}$ only depends on $n$ modulo squares. Therefore, there exists a sign $\varepsilon_f$ such that $d_{n,\varepsilon}=1$ if and only if $\varepsilon_f\varepsilon n$ is a square modulo $p$. Let $n \in \F_p^{\times}$ and $\epsilon \in \{\pm 1\}$ be the sign such that $n\varepsilon$ is a square modulo $p$. The subspace of $J_n$ on which $\mathbb{T}_1$ acts through $\mu_{f,1}$ is then (as a $\mathbb{T}_{\rho^{\ast},1} \otimes F_{\lambda}$-module) isomorphic to $M_{f,U^{\varepsilon_f\varepsilon},\alpha} \simeq M_{f,U^{n\varepsilon_f},\alpha} \simeq M_{f,U^{\varepsilon_f},\alpha^n}$, and the conclusion follows. 

What remains to be shown is that $(\varepsilon_f)_f$ is canonical. We will mainly exploit the following relation: for $t \in \mathbb{T}_{\rho^{\ast},1} \otimes F_{\lambda}$, one has 
\begin{align}
\sum_{f \in \mathscr{C}_M}{\mrm{Tr}(t \mid M_{f,U^{\varepsilon_f},\alpha})} = \mrm{Tr}(t \mid \Tate{\ell}{J_{\Gamma}^{\alpha}(p)} \otimes F_{\lambda}) - \sum_{f \in \mathcal{N}'}{\mrm{Tr}(t \mid M_{f,\lambda})}. \label{trace-formula}
\end{align}

\emph{Part 3: Part 2 holds for exactly one choice of $(\varepsilon_f)_{f \in \mathscr{C}_M}$, which does not depend on $\lambda$.}

Let $\lambda,\lambda'$ be two maximal ideals of $\OO_F$ of residue characteristics $\ell,\ell'$ and $(\varepsilon_f)_f,(\varepsilon'_f)_f$ be two choices of sign such that the claim holds for $(\Gamma,P,Q,\lambda,(\varepsilon_f)_f)$ and $(\Gamma,P,Q,\lambda',(\varepsilon'_f)_f)$. Let $q$ be a prime unramified for $\rho$ (so that $\Gamma,\alpha$ spread out to $\Z_{(q)}$) and coprime to $p\lambda\lambda'$ and fix a place $\mfk{q}$ of $\Qbar$ above $q$. 

Let $\gamma \in G_{\Q}$ be a Frobenius element for $\mfk{q}$ and $g \in \GL{\F_p}, n \geq 1$ are such that $n\det{g} \equiv 1 \mod{p}$, then $T_n(g,\gamma) \in \mrm{End}(\Tate{\ast}{J_{\Gamma}^{\alpha}(p)}_{\Q}) \simeq \mrm{End}(\Tate{\ast}{J_{\Gamma}^{\alpha}(p)_{\F_q}})$ (with $\ast=\ell,\ell'$) comes from the endomorphism $(T_n\gamma) \circ (\phi \times \mrm{id})$ of $J_{\Gamma}^{\alpha}(p)_{\overline{\Z}_{\mfk{q}}/\mfk{q}}$, where $\phi$ is the absolute Frobenius of $J_{\Gamma}^{\alpha}(p)_{\F_q}$ (cf. e.g. \cite[(12.1.1)]{KM}). By \cite[Proposition 12.9]{MilAb}, one has therefore
\[\mrm{Tr}(T_n(g,\gamma) \mid \Tate{\ell}{J_{\Gamma}^{\alpha}(p)}) = \mrm{Tr}(T_n(g,\gamma) \mid \Tate{\ell'}{J_{\Gamma}^{\alpha}(p)}).\] 

For $f \in \mathcal{N}'$, $\mrm{Tr}(T_n(g,\gamma) \mid M_{f,\lambda}) \in F_{\lambda},\mrm{Tr}(T_n(g,\gamma) \mid M_{f,\lambda'}) \in F_{\lambda'}$ both lie in $F$ and are equal. Similarly, $f \in \mathscr{C}_M$ and $\varepsilon=\pm 1$, $\mrm{Tr}(T_n(g,\gamma) \mid M_{f,U^{\varepsilon},\alpha,\lambda}) \in F_{\lambda},\mrm{Tr}(T_n(g,\gamma) \mid M_{f,U^{\varepsilon},\alpha,\lambda'}) \in F_{\lambda'}$) both lie in $F$ and are equal. Therefore, the result of Part 2 (applied to $\lambda,\lambda'$) and \ref{trace-formula} imply that 
\[\sum_{f \in \mathscr{C}_M}{\mrm{Tr}(T_n(g,\gamma) \mid M_{f,U^{\varepsilon_f},\alpha,\lambda})} = \sum_{f \in \mathscr{C}_M}{\mrm{Tr}(T_n(g,\gamma) \mid M_{f,U^{\varepsilon'_f},\alpha,\lambda})}.\]

By Chebotarev's theorem, this implies that for every $t \in \mathbb{T}_{\rho^{\ast},1} \otimes F_{\lambda}$, one has 
\[\sum_{f \in \mathscr{C}_M}{\mrm{Tr}(t \mid M_{f,U^{\varepsilon_f},\alpha,\lambda})} = \sum_{f \in \mathscr{C}_M}{\mrm{Tr}(t \mid M_{f,U^{\varepsilon'_f},\alpha,\lambda})}.\] By \cite[Theorem 27.8]{Curtis-Reiner} and \ref{dttc-5}, this implies that $\varepsilon_f=\varepsilon'_f$ for every $f \in \mathscr{C}_M$.

\emph{Part 4: The choice of $(\varepsilon_f)_{f \in \mathscr{C}_M}$ in Part 2 is canonical.}

We may assume that $\mathscr{C}_M$ is not empty, and thus that $p \equiv -1\mod{4}$. Let $\mathscr{J}_{\Gamma,P,Q,\alpha}$ denote the $\mathbb{T}_{\rho^{\ast},1}$-module $J_{\Gamma}^{\alpha}(p)(\Qbar)$, where the underlying $\GL{\F_p}$-action comes from the basis $(P,Q)$. Let $\Gamma'$ be another $p$-torsion group over $\Q$ equipped with a Weil pairing $\alpha'$ and $s: \Gamma_{\Qbar} \rar \Gamma'_{\Qbar}$ be a $\Qbar$-isomorphism, and assume that $\alpha' \circ (s \times s) = \alpha$; let $(P',Q')$ be a basis of $\Gamma'(\Qbar)$. Let $M \in \GL{\F_p}$ and $\begin{pmatrix}P'\\Q'\end{pmatrix}=M\begin{pmatrix}s(P)\\s(Q)\end{pmatrix}$. Let $\rho'$ be the representation attached to this basis.

Part 2 attaches signs $(\varepsilon_f)_{f \in \mathscr{C}_M}$ and $(\varepsilon'_f)_{f \in \mathscr{C}_M}$ to $(\Gamma,P,Q)$ and $(\Gamma',P',Q')$, and we are done if we show that $\varepsilon_f=\varepsilon'_f$ for every $f \in \mathscr{C}_M$ -- and Part 3 implies that we can fix $\lambda$ for the rest of the argument. 

Now, $s$ induces isomorphisms $s^{\ast}: \mathscr{J}_{\Gamma',P',Q',(\alpha')^n} \rar \mathscr{J}_{\Gamma,P,Q,\alpha^n}$ (for each $n \in \F_p^{\times}$) that commute to $\mathbb{T}$. The map $\tau: \mathbb{T}_{(\rho')^{\ast}} \rar \mathbb{T}_{\rho^{\ast}}$ sending $T_n (g,\gamma)$ to $T_n (M^{-1}g(\rho')^{\ast}(\gamma)M\rho^{\ast}(\gamma)^{-1},\gamma)$ is a ring isomorphism, and induces isomorphisms $\tau_1: \mathbb{T}_{(\rho')^{\ast},1} \rar \mathbb{T}_{\rho^{\ast},1}, \,\tau_{1,1}: \mathbb{T}_{\rho'^{\ast},1,1} \rar \mathbb{T}_{(\rho)^{\ast},1,1}$. By Lemma \ref{j1-p-q}, $s^{\ast}$ satisfies $s^{\ast}(t \cdot x) = \tau(t) \cdot s^{\ast}(x)$ for $t \in \mathbb{T}_{(\rho')^{\ast}}$ and $x \in \mathscr{J}_{\Gamma',P',Q',(\alpha')^n}$. 

Therefore, for $t \in \mathbb{T}_{(\rho')^{\ast},1} \otimes F_{\lambda}$, one has $\mrm{Tr}(t \mid \Tate{\ell}{J_{\Gamma'}^{\alpha'}(p)} \otimes F_{\lambda})=\mrm{Tr}(\tau_1(t) \mid \Tate{\ell}{J_{\Gamma}^{\alpha}(p)} \otimes F_{\lambda})$.  

On the other hand, one checks easily that 
\begin{align*}
& \forall f \in \mathcal{N}', &&\forall t \in \mathbb{T}_{(\rho')^{\ast},1} \otimes F_{\lambda},&&& \mrm{Tr}(t \mid M_{f,\rho'}) = \mrm{Tr}(\tau_1(t) \mid M_{f,\rho}), \\
&\forall (f,\varepsilon) \in \mathscr{C}_M \times \{\pm 1\}, &&\forall t \in \mathbb{T}_{(\rho')^{\ast},1} \otimes F_{\lambda}, &&& \mrm{Tr}\left(t \mid M_{f,U^{\varepsilon},\alpha',\rho'}\right) = \mrm{Tr}\left(\tau_1(t) \mid M_{f,U^{\varepsilon},\alpha,\rho}\right).
\end{align*} 
Indeed, we may assume $t=T_n(g,\gamma)$, and the trace is zero unless $\gamma \in G_K$ and $n,\det{g} \in \F_p^{\times}$ are squares, in which case we can use the decomposition appearing in \ref{dttc-5}. Therefore, for $t \in \mathbb{T}_{\rho^{\ast},1} \otimes F_{\lambda}$, one has $\sum_{f \in \mathscr{C}_M}{\mrm{Tr}\left(t \mid M_{f,U^{\varepsilon_f},\alpha',\rho'}\right)} = \sum_{f \in \mathscr{C}_M}{\mrm{Tr}\left(t \mid M_{f,U^{\varepsilon'_f},\alpha',\rho'}\right)}$. By \cite[Theorem 27.8]{Curtis-Reiner} and \ref{dttc-5}, this implies that $\varepsilon_f=\varepsilon'_f$ for every $f \in \mathscr{C}_M$ and the conclusion follows.}

\section{Structure of $J_{\Gamma}(p)$ in the non-split Cartan case}
\label{sect:cartan-bsd}

\subsection{Some local results}
\label{subsect:local-standard}
None of the results in this subsection are new, but we gathered them here for reference in other parts of the paper. 

In this section, $p \geq 11$ is a prime and $E$ is an elliptic curve over a field $F$ of characteristic zero such that the image of the morphism $\rho: G_F \rar \Aut{E[p](\overline{F})}$ is contained in the normalizer $N$ of a non-split Cartan subgroup $C$.

\lem[elliptic-curves-cartan-complex]{If $F=\R$, then $\rho^{-1}(C)$ is trivial. In particular, if $F=\Q$, the number field $K/\Q$ such that $G_K=\rho^{-1}(C)$ is imaginary quadratic. }

\demo{The only non-trivial element of $G_{\R}$ is the complex conjugation, whose action on $E[p](\C)$ has determinant $-1$ and order $2$, so $\rho(c) \notin C$. This solves the question for $F=\R$. For $F=\Q$, the group $N/C$ has cardinality two, so $K/\Q$ has degree at most $2$. Moreover, $K \not\subset \R$ by the first part of the argument, so $K$ is imaginary quadratic.}

\prop[elliptic-curves-cartan-padic]{Assume that $F=\Q_p$. Let $K/\Q_p$ be the extension such that $\rho^{-1}(C)=G_K$. Then $K/\Q_p$ is the quadratic unramified extension, so the image by $\rho$ of the inertia subgroup $I_p$ is a subgroup $I$ of $C$. The elliptic curve $E$ has potentially good supersingular reduction. 
If $I \neq C$, then $[C:I]=3$, $p$ is congruent to $2$ or $5$ modulo $9$, $v_p(j(E)) \geq 4$ and $E$ does not have good reduction.  
}

\demo{This is essentially a local version of the reunion of \cite[Proposition 1.4, Appendix B]{LFL} and \cite[Proposition 3.4]{FL}, where we removed the ``non-CM'' assumption which did not play any role.}

\prop[elliptic-curves-cartan-nram]{Suppose that $F=\Q_{\ell}$ with $\ell \neq p$. 
\begin{itemize}[noitemsep,label=\tiny$\bullet$]
\item $\rho$ is unramified if and only if $E$ has good or multiplicative reduction,  
\item if $E$ has potentially multiplicative reduction, then $j(E) \in \Q_{\ell}$ is a $p$-th power, $\ell \equiv \pm 1 \mod{p}$, and the image by $\rho$ of the inertia subgroup is contained in $\{\pm \mrm{id}\}$. 
\item if $E$ has multiplicative reduction and $\ell \equiv 1 \mod{p}$, then $\rho(\Fr) = \pm \mrm{id}$, 
\item if $E$ has multiplicative reduction and $\ell \equiv -1\mod{p}$, $\rho(\Fr) \notin C$ is an element of order $2$.
\end{itemize} 
In particular, if $\rho$ is unramified, then $E$ has good reduction if and only if $E$ has potentially good reduction, if and only if $j(E) \in \Z_{\ell}$. 
 }
 
\demo{If $E$ has potentially good reduction, then $\rho$ is ramified if and only if $E$ has bad reduction, by the N\'eron--Ogg--Shafarevich criterion (e.g. \cite[Cor. 2 (b) to Th. 2]{GoodRed}). We may thus assume that $E$ has potentially multiplicative reduction, or equivalently that $j(E) \notin \Z_{\ell}$ (by \cite[Proposition VII.5.5]{AEC1}). The rest of the claim is a direct consequence of Tate uniformization, cf. \cite[Chapter V]{AEC2}. }

\prop[elliptic-curves-cartan-potgood]{Assume that $F=\Q_{\ell}$ with $\ell \neq p$ and $E$ has potentially good reduction. Let $\Delta$ be the discriminant of a minimal model of $E$. Let $v: \Q_{\ell}^{\times} \rar \Z$ denote the normalized valuation. Then:
\begin{itemize}[noitemsep,label=\tiny$\bullet$]
\item The conductor of $E$ and the Artin conductor of $E[q](\Qbar)$ for any prime $q \nmid 6\ell$ have the same valuation.
\item Let $I_{\ell}$ be the inertia subgroup and $J_{\ell}$ be its image in $\Aut{E[p](\overline{\Q_{\ell}})}$. Then for any $n$ prime to $\ell$, the obvious morphism $I_{\ell} \rar \Aut{E[n](\overline{\Q_{\ell}})}$ factors through the projection $I_{\ell} \rar J_{\ell}$ and the injection $\operatorname{SL}(E[n](\overline{\Q_{\ell}}))\subset \Aut{E[n](\overline{\Q_{\ell}})}$; moreover the induced map $J_{\ell} \rar \operatorname{SL}(E[n](\overline{\Q_{\ell}}))$ is injective if $n \geq 3$.
\end{itemize}}

\demo{The first bullet point is classical (see e.g. \cite[II,Proposition]{Kraus-thesis}\footnote{The statement at \emph{loc.cit.} is technically global, but its proof is purely local.}). The second point is \cite[Cor. 2 to Th. 2]{GoodRed} (because the determinant of the Galois action on $\Tate{q}{E}$ is the cyclotomic character, which is unramified if $q \neq \ell$). The third and fourth bullet points come from \cite[\S 5.6]{Serre-image}. }

\subsection{Further decomposition for normalizers of non-split Cartan subgroups}
\label{subsect:cartan-decomp}

Let $p \geq 7$ be a prime and $N \leq \GL{\F_p}$ be the normalizer of a non-split Cartan subgroup $C$. Let $\varepsilon: N \rar \{\pm 1\}$ be the character with kernel $C$ and $\eta: \F_p^{\times} \rar \{\pm 1\}$ be the unique non-trivial quadratic character. Let $N[\eta]=N \cap \ker{\eta}$, $C[\eta]=C \cap \ker{\eta}$. Moreover, let $\nu_2 \subset N \backslash C$ denote the conjugacy class of elements with characteristic polynomial $X^2-1$. 

From now and until the end of Section \ref{subsect:constructing-quotients}, $F \subset \C$ is a number field which is Galois over $\Q$, contains the $p(p^2-1)$-th roots of unity, the Fourier coefficients of every newform in $\mathcal{N}$ and, for every $f \in \mathscr{C}_M$ and every prime $q$ split in $\Q(\sqrt{-p})$, the roots of $X^2-a_q(f)X+q$. Note that a character $\alpha: N \rar F^{\times}$ can be uniquely written as $\varepsilon^s \cdot (\chi \circ \det)$, where $\chi$ is a Dirichlet character modulo $p$ and $s \in \{0,1\}$. 

For $f \in \mathcal{N}$, we make the following definitions: 
\begin{itemize}[noitemsep,label=\tiny$\bullet$]
\item $\mathscr{R}_1(f)$ is the set of characters $N \rar F^{\times}$ appearing in $R(f)_{|N}$,
\item if $f$ has no complex multiplication (resp. if $f$ has CM), $\mathscr{R}_1(f)_{rep}=\mathscr{R}_1(f)$ (resp. $\mathscr{R}_1(f)_{rep}$ is the set of $\psi_{|N \cap \ker{\eta}}{\psi}$ for $\psi \in \mathscr{R}_1(f)_{rep}$),
\item $\mathscr{R}_2(f)$ is the set of characters $\psi: C \rar F^{\times}$ appearing in $R(f)_{|C}$ such that $\psi \neq \psi^p$,
\item if $f$ has no complex multiplication (resp. if $f$ has CM), $\mathscr{R}_2(f)_{rep}$ is the set of $\mrm{Ind}_C^N{\psi}$ over all $\psi \in \mathscr{R}_2(f)$ (resp. the set of $\mrm{Ind}_{C[\eta]}^{N[\eta]}{\psi_{|C[\eta]}}$). 
\end{itemize}

Note that one always has, by construction, a map $\mathscr{R}_i(f) \rar \mathscr{R}_i(f)_{rep}$. By definition, for $f \in \mathscr{N}$ without (resp. with) complex multiplication, the elements of $\mathscr{R}_i(f)_{rep}$ are representations of $N$ (resp. $N[\eta]$).

A character calculation yields that:
\begin{itemize}[noitemsep,label=\tiny$\bullet$]
\item For $f \in \mathscr{S}$, $\mathscr{R}_1(f)$ consists of the single character which is trivial on $\nu_2$ and whose restriction to $C$ is $(\eta \circ \det)$, and $\mathscr{R}_2(f)$ is the set of all $\psi: C/\F_p^{\times}\mrm{id}_2 \rar F^{\times}$ with order at least $3$. 
\item If $f \in \mathscr{P}$ and $\chi$ is its character, $\mathscr{R}_1(f)$ consists of the two characters trivial on $\nu_2$ whose restriction to $C$ is of the form $\psi\circ\det$, for each of the two Dirichlet characters $\psi$ modulo $p$ such that $\chi =\psi^2$, and $\mathscr{R}_2(f)$ is the set of all $\psi: C \rar F^{\times}$ such that $\psi_{|\F_p^{\times}\mrm{id}_2} = \chi$ and $\psi \neq \psi^p$. 
\item If $f \in \mathscr{C} \backslash \mathscr{C}_M$, let $\phi, \phi^p: C/\F_p^{\times} \rar \C^{\times}$ be the two characters attached to $R(f)$. Let $a \in C$ be a non-scalar element with scalar square. Then $\mathcal{R}_1(f)$ consists of two characters $\chi_1$ and $\chi_2$, whose respective restrictions to $C$ are trivial and $\eta\circ\det$, and such that $\chi_1(\nu_2)=-\phi(a)=-\chi_2(\nu_2)$, and $\mathscr{R}_2(f)$ is the collection of $\psi: C/\F_p^{\times}\mrm{id}_2 \rar F^{\times}$, for $\psi: C \rar F^{\times}$ of order at least $3$ and distinct from $\phi,\phi^p$,
\item If $f \in \mathscr{C}_M$, then $\mathscr{R}_1(f) = \{\varepsilon^{\frac{p-3}{4}},\varepsilon^{\frac{p-3}{4}}\eta(\det)\}$, and $\mathscr{R}_2(f)$ consists of the $\psi: C/\F_p^{\times}\mrm{id}_2 \rar F^{\times}$ such that $\psi^4 \neq 1$. 
\end{itemize}

\lem[ri-twist]{Let $f \in \mathcal{N}$. If $\sigma \in \Aut{\C}$ and $\chi$ is a Dirichlet character modulo $p$, then $\psi \mapsto (\sigma \circ \psi) \cdot (\chi \circ \det)$ defines a bijection $\mathscr{R}_i(f) \rar \mathscr{R}_i(\sigma(f) \otimes \chi)$ that respects the fibres of $\mathscr{R}_i(\bullet) \rar \mathscr{R}_i(\bullet)_{rep}$.}

\lem[ri-no-cm]{Let $f \in \mathcal{N}$ without complex multiplication. Then 
\begin{itemize}[noitemsep,label=\tiny$\bullet$]
\item $R(f)_{|N}$ is the direct sum of the representations of $\mathcal{R}_1(f)$ and those in $\mathcal{R}_2(f)_{rep}$, which are irreducible,
\item $\psi \mapsto \psi^p$ is an involution of $\mathscr{R}_2(f)$ without fixed point, and the quotient of $\mathscr{R}_2(f)$ by this involution is precisely the map $ \mathscr{R}_2(f) \rar \mathscr{R}_2(f)_{rep}$. 
\end{itemize}}

\lem[ri-with-cm]{Let $f \in \mathcal{N}$ with complex multiplication. Then
\begin{itemize}[noitemsep,label=\tiny$\bullet$]
\item $R(f)_{|N[\eta]}$ is the direct sum of the representations in $\mathscr{R}_1(f)_{rep}$ and $\mathscr{R}_2(f)_{rep}$, which are irreducible, 
\item The multiplication by $\eta(\det)$ is an involution of $\mathscr{R}_1(f)$ without fixed points, and the quotient of $\mathscr{R}_1(f)$ by this involution is exactly the map $\mathscr{R}_1(f) \rar \mathscr{R}_1(f)_{rep}$. 
\item The group $\F_2^{\oplus 2}$ acts on $\mathscr{R}_2(f)$ by $(a,b) \cdot \psi \mapsto \psi^{p^a} \cdot \eta(\det)^b$ without fixed points, and the orbits are exactly the fibres of the map $\mathscr{R}_2(f) \rar \mathscr{R}_2(f)_{rep}$.
\end{itemize}}

\prop[tate-module-twisted-cartan]{Let $p \geq 7$ and $\Gamma$ be a $p$-torsion group over $\Q$ with Weil pairing $\alpha$. Let $P,Q \in \Gamma(\Qbar)$ form a basis of $\Gamma(\Qbar)$. Let $\rho: G_{\Q} \rar \GL{\F_p}$ be the matrix of the action of $G_{\Q}$ in the basis $(P,Q)$. Assume that the image of its contragredient $\rho^{\ast}$ is contained in $N$. Let $\tilde{\rho}: (g,\gamma) \in \GL{\F_p} \rtimes_{\rho^{\ast}} G_{\Q} \mapsto g\rho^{\ast}(\gamma)\in \GL{\F_p}$. Recall that the choice of $P,Q$ defines an action of $\GL{\F_p}$ on $\coprod_{n \in \F_p^{\times}}{X_{\Gamma}^{\alpha^n}(p)_{\Qbar}}$.

Let $\lambda \subset \OO_F$ be a maximal ideal. Then, as a $\mathbb{T}[N\rtimes_{\rho^{\ast}} G_{\Q}]_1 \otimes F_{\lambda}$-module, $\Tate{\ell}{J_{\Gamma}^{\alpha}(p)} \otimes_{\Z_{\ell}} F_{\lambda}$ is isomorphic to the following direct sum: 

\begin{itemize}[noitemsep,label=\tiny$\bullet$]
\item for every $f \in \mathscr{S}$ and every $R \in \mathscr{R}_1(f)\cup\mathscr{R}_2(f)_{ind}$, a copy of $V_{f,\lambda} \otimes (R\circ \tilde{\rho})$,
\item for every $f \in \mathscr{P}$ modulo complex conjugation and every $R \in \mathscr{R}_1(f)\cup\mathscr{R}_2(f)_{ind}$, a copy of $V_{f,\lambda} \otimes (R\circ \tilde{\rho})$, 
\item for every $f \in\mathscr{C} \backslash \mathscr{C}_M$ modulo quadratic twists and every $R \in \mathscr{R}_1(f)\cup\mathscr{R}_2(f)_{ind}$, a copy of $V_{f,\lambda} \otimes (R\circ \tilde{\rho})$, 
\item for every $f \in \mathscr{C}_M$, a copy of $V_{f,\lambda} \otimes (\varepsilon^{\frac{p-3}{4}}\circ \tilde{\rho})$,
\item for every $f \in \mathscr{C}_M$ and every fibre $\Psi$ of $\mathscr{R}_2(f) \rar \mathscr{R}_2(f)_{rep}$, a copy of $V_{f,\lambda} \otimes ((\mrm{Ind}_C^N{\psi})\circ \tilde{\rho})$ for some $\psi \in \Psi$ (the representation does not depend on the choice of $\psi$),
\end{itemize}}

\demo{The claim for the non-CM summands follows from Proposition \ref{decomp-tate-twisted-general} and the decomposition of the representation $R(f)_{|N}$ by Lemma \ref{ri-no-cm}. So all that remains is the discussion of the summands attached to $f \in \mathscr{C}_M$, so we may assume that $p \equiv -1\mod{4}$, and let $K=\Q(\sqrt{-p})$.

Let $\Sigma_f$ be the direct sum of all summands attached to $f$. Let $\sigma: K \rar F$ be an embedding, $\psi_{\lambda,\sigma}$ be attached to $f,\sigma$ by Proposition \ref{compatible-cm-systems}, and $U=\begin{pmatrix}1 & 1\\0 & 1\end{pmatrix}$, and let $\zeta \in F$ be a primitive $p$-th root of unity. By Proposition \ref{decomp-tate-twisted-general}, we need to show that $\Sigma_f$ is isomorphic as a $\mathbb{T}[N[\eta] \rtimes_{\rho^{\ast}} G_{\Q}]_1 \otimes F_{\lambda}$-module to 
\[\Sigma'_f := \mrm{Ind}_{N[\eta] \rtimes_{\rho^{\ast}} G_{K}}^{N[\eta] \rtimes_{\rho^{\ast}} G_{\Q}}{\left\{[\psi_{\lambda,\sigma}\circ \pi]\otimes [R_{U,\zeta} \circ \tilde{\rho}_{|N[\eta] \rtimes_{\rho^{\ast}} G_K}]\right\}}.\]

Indeed, for $g \in N \backslash \ker{\eta}$, $n \geq 1$ such that $n\det{g}\equiv 1\mod{p}$, $\gamma \in G_{\Q}$, $n$ is not a square modulo $p$, so $T_n(g,\gamma)$ kills both $\Sigma_f$ and the summand of $\Tate{\ell}{J_{\Gamma}^{\alpha}(p)}$ attached to $f$ (by Proposition \ref{decomp-tate-twisted-general}). As $\mathbb{T}[N \rtimes_{\rho^{\ast}} G_{\Q}]_1$-modules, one has $V_{f,\lambda} \simeq V_{f,\lambda} \otimes [\eta(\det) \circ \tilde{\rho}]$, so $\Sigma_f$ is independent from the choice of elements in the fibres of $\mathscr{R}_i(f) \rar \mathscr{R}_i(f)_{rep}$. \footnote{There is one such fibre when $i=1$, hence the explicit choice of character --- the point is, in this case, that we could have chosen the other one as well.}

Note that the image of $\tilde{\rho}_{|N[\eta])\rtimes_{\rho^{\ast}} G_K}$ is $N[\eta]$ and that it (tautologically) extends to a group homomorphism $\tilde{\rho}: N[\eta] \rtimes_{\rho^{\ast}} G_{\Q} \rar N$. Let $R'(f) = \varepsilon^{\frac{p-3}{4}} \oplus \bigoplus_{\psi}{\mrm{Ind}_C^N{\psi}}$, where $\psi$ runs through characters in $\mathscr{R}_2(f)$ and we pick exactly one of them in each fibre of $\mathscr{R}_2(f) \rar \mathscr{R}_2(f)_{rep}$. Thus $\Sigma_f \simeq [V_{f,\lambda} \circ \pi] \otimes [R'(f)_{|N \rtimes_{\rho^{\ast}} G_{\Q}}]$. By a character calculation, $R_{U,\zeta}$ is the restriction $R'(f)_{|N[\eta]}$, which implies that $\Sigma'_f$ is isomorphic to $\left[\mrm{Ind}_K^{\Q}{\psi_{\lambda,\sigma}}\right]\circ \pi \otimes [R'(f) \circ \tilde{\rho}] \simeq \Sigma_f$. }

\cor[tate-module-twisted-cartan-gq]{Keep the notations of Proposition \ref{tate-module-twisted-cartan}. Then, as a $\mathbb{T}_1[G_{\Q}] \otimes F_{\lambda}$-module, $\Tate{\ell}{J_{\Gamma}^{\alpha}(p)} \otimes_{\Z_{\ell}} F_{\lambda}$ is isomorphic to the following direct sum: 

\begin{itemize}[noitemsep,label=\tiny$\bullet$]
\item for every $f \in \mathscr{S}$ and every $R \in \mathscr{R}_1(f)\cup\mathscr{R}_2(f)_{ind}$, a copy of $V_{f,\lambda} \otimes (R\circ \rho^{\ast})$,
\item for every $f \in \mathscr{P}$ modulo complex conjugation and every $R \in \mathscr{R}_1(f)\cup\mathscr{R}_2(f)_{ind}$, a copy of $V_{f,\lambda} \otimes (R\circ \rho^{\ast})$, 
\item for every $f \in\mathscr{C} \backslash \mathscr{C}_M$ modulo quadratic twist and every $R \in \mathscr{R}_1(f)\cup\mathscr{R}_2(f)_{ind}$, a copy of $V_{f,\lambda} \otimes (R\circ \rho^{\ast})$, 
\item for every $f \in \mathscr{C}_M$ and every $u \in \mathscr{R}_1(f)$, a copy of $V_{f,\lambda} \otimes (u\circ \rho^{\ast})$,
\item for every $f \in \mathscr{C}_M$ and every $\psi \in \mathscr{R}_2(f)_{ind}$, a copy of $V_{f,\lambda} \otimes ((\mrm{Ind}_C^N{\psi})\circ \rho^{\ast})$,
\end{itemize}
All of these $F_{\lambda}[G_{\Q}]$-modules form compatible systems. 
 }

\demo{The decomposition is an easy consequence from the previous Proposition. For the second part, it is known that the $(V_{f,\lambda})_{\lambda}$ and Artin representations define compatible systems, and that compatible systems are stable under tensor products.

}

\cor[properties-tate-module-twisted-cartan-gq]{Keep the notations of Corollary \ref{tate-module-twisted-cartan-gq}. All the $F_{\lambda}[G_{\Q}]$-modules described in this Corollary are self-dual up to Tate twist and irreducible if $\rho^{\ast}(G_{\Q})=N$. Moreover, the completed complex $L$-function attached to the corresponding complete compatible systems (as in Proposition \ref{carayol-local-global}) possesses a holomorphic continuation and satisfies a functional equation $\Lambda(s) = \varepsilon \Lambda(2-s)$, where $\varepsilon$ is a sign.  }

\demo{The self-duality is a direct calculation based on the fact that $\det{\rho^{\ast}}$ is the inverse of the cyclotomic character modulo $p$. Let $K$ be the number field such that $\rho^{-1}(C)=G_K$, then $K$ is imaginary quadratic by Lemma \ref{elliptic-curves-cartan-complex}. 

Let us discuss the irreducibility. It is a direct consequence of e.g. \cite[(4.4)]{Ribet-Antwerp5} and \cite[Lemma 2.3]{Studnia-Euler} when $f$ has no complex multiplication. Since $V_{f,\lambda}$ is absolutely irreducible, it remains irreducible when the representation is tensored by a character, so all that remains is the case of $V_{f,\lambda} \otimes ((\mrm{Ind}_C^N{\psi})\circ \rho^{\ast})$ for $f \in \mathscr{C}_M$ and $\psi \in \mathscr{R}_2(f)$. We may thus assume that $p \equiv -1 \mod{4}$, and we let $R'$ denote the $F_{\lambda}[N]$-module $\mrm{Ind}_C^N{\psi}$.  

Since $\det{\rho}$ is the mod $p$ cyclotomic character, it sends the complex conjugation to $-1$, so $\rho^{\ast}$ sends the complex conjugation to an element with characteristic polynomial $X^2-1$: such an element is not contained in $C$, whence $K$ is quadratic imaginary and distinct from $\Q(\sqrt{-p})$. Let $\psi_f$ denote one of the two $\lambda$-adic characters attached to $f$ by Corollary \ref{compatible-cm-systems}, then one has 

\[V_{f,\lambda} \otimes (R'\circ \rho^{\ast}) \simeq \mrm{Ind}_{G_{\Q(\sqrt{-p})}}^{\Q}{\psi_f \otimes [R'\circ \rho^{\ast}_{|G_{\Q(\sqrt{-p})}}]}.\] 

By Proposition \ref{compatible-cm-systems}, $\psi_f$ does not agree with its $\mrm{Gal}(\Q(\sqrt{-p})/\Q)$-conjugate over any open subgroup of $G_{\Q(\sqrt{-p})}$, so by Mackey's criterion \cite[Cor. to Prop. 23]{SerreLinReps}, it is enough to show that $\psi_f \otimes [R'\circ \rho^{\ast}_{|G_{\Q(\sqrt{-p})}}]$ is an absolutely irreducible $F_{\lambda}[G_{\Q(\sqrt{-p})}]$-module. Because $\rho^{\ast}(G_{\Q})=N$, $\rho^{\ast}(G_{\Q(\sqrt{-p})}) = N [\eta]$, so it is enough to check that $R'_{|N[\eta]}$ is absolutely irreducible. By Mackey's criterion, this amounts to showing that $\psi_{|C[\eta]}$ does not agree with its $p$-th power. If one had $\psi_{|C[\eta]}=\psi^p_{|C[\eta]}$, then $\psi^{2(p-1)}$ would be trivial; since $\psi^{p+1}$ is trivial, $\psi^4$ would be trivial, which is impossible because $\mathscr{R}_2(f)$ does not contain a character of order dividing $4$.

Finally, we discuss the functional equation. In the case of $V_{f,\lambda} \otimes (R\circ \rho^{\ast})$ with $R \in \mathscr{R}_1(f)$, $R \circ \rho^{\ast}$ is a one-dimensional representation of $G_{\Q}$, so is given by a Dirichlet character, and this is the classical theory for modular forms (Proposition \ref{carayol-local-global}). 

Now, if $\psi \in \mathscr{R}_2(f)$, let $R_{\psi} :=\mrm{Ind}_C^N{\psi}$, then $\mrm{Tr}(\nu_2 \mid R_{\psi})=0$, so $R_{\psi} \circ \rho^{\ast}$ is an odd two-dimensional Artin representation. If it is reducible, then it is the sum of two characters of $G_{\Q}$ and we are brought back to the previous paragraph. If it is irreducible, it follows from the well-known dihedral case of Serre's conjecture (see e.g. the discussion of \cite[pp. 517--518]{Del-Ser}) that $R_{\psi} \circ \rho^{\ast}$ is attached to a newform of weight one, and the result follows from the theory of the Rankin--Selberg product and local-global compatibility (i.e. Proposition \ref{rankin-selberg-galois}).   
}

One can in particular compute the functional equation signs for the various summands appearing in Corollary \ref{tate-module-twisted-cartan-gq}. As discussed in the Introduction, at least under the assumption that $\rho^{\ast}: G_{\Q} \rar N$ is surjective and that $\Gamma_{\Q_v}$ comes from an elliptic curve over $\Q_v$ for every place $v$ of $\Q$, this is carried out in the author's PhD dissertation (cf. \cite[Theorem 10]{Studnia-thesis}). In most cases, and in particular whenever $\Gamma_{\Q_p}$ comes (up to quadratic twist) from an elliptic curve with good reduction at $p$, the functional equation sign is equal to $-1$, and the central value of the $L$-function is then zero. However, there are two situations under which the functional equation sign can be equal to $+1$ (and therefore the central value of the $L$-function is not forced to vanish). We discuss only the first one, which is the simpler of the two.  

\prop[sign-plus]{Assume that $\rho^{\ast}: G_{\Q} \rar N$ is surjective and that $\Gamma_{\Q_p}$ is the $p$-torsion subgroup scheme of some elliptic curve $E_p/\Q_p$. Assume furthermore that, if $I_p$ denotes an inertia subgroup at $p$, $\rho^{\ast}(I_p)$ is properly contained in $C$.

Let $f \in \mathscr{S}$ and $\psi \in \mathscr{R}_2(f)$ be one of the two characters of order $3$. Let $g$ be the weight one newform attached to the odd irreducible Artin representation $\left(\mrm{Ind}_C^N{\psi}\right)\circ \rho$. The functional equation sign of the completed Rankin--Selberg $L$-function $\Lambda(f,g,s)$ is $+1$.}

\demo{As we previously saw, the field $K/\Q$ such that $G_K=(\rho^{\ast})^{-1}(C)$ is imaginary quadratic; by Proposition \ref{elliptic-curves-cartan-padic}, $K$ is inert at $p$, one has $\rho^{\ast}(I_p)=3C$, and $p \equiv 2,5 \mod{9}$. As a consequence, $\mathscr{R}_2(f)$ does contain a character of order $3$. Let $\varepsilon_K$ denote the character of $K$, and $D$ be the conductor of $\varepsilon_K$. 

By Proposition \ref{properties-tate-module-twisted-cartan-gq}, $\Psi := \mrm{Ind}_C^N{\psi} \circ \rho^{\ast}$ is an irreducible two-dimensional odd Artin representation of $G_{\Q}$, and it is attached to a weight one newform $g$. Now, the character of the representation $\mrm{Ind}_C^N{\psi}$ of $N$ is real-valued and its determinant is exactly the quadratic character $\varepsilon$. Because $\rho^{\ast}(I_p)=3C \subset \ker{\psi}$, $\Psi$ is unramified at $p$. By \cite[\S 4]{Del-Ser}, this implies that $g$ has real coefficients, complex multiplication by $K$, character $\varepsilon_K$, and level $M$ prime to $p$. In particular $g,f$ have coprime levels. The conclusion follows from \cite[Theorem 2.2, Example 2]{Li-RS} \footnote{Technically, \emph{loc.cit.} requires some extra assumptions, named (A), (B), (C) in that paper. While (A) and (C) always hold for newforms with coprime levels, it is not clear that this condition is sufficient to obtain (B). The formulation in \emph{loc.cit.} makes it apparent that the formula in Example 2 of \emph{op.cit.} holds even when (B) is not satisfied. At any rate, it is not a difficult exercise to check that (B) is, in this case, satisfied.} (where the squares of the Atkin--Lehner pseudo-eigenvalues are given by \cite[Proposition 1.1]{AL78}).    

}

\subsection{Central characters}
\label{subsect:central-characters}

For $f \in \mathcal{N}$ of character $\chi$, we denote as in Corollary \ref{t-structure} by $\mu_{f}: \mathbb{T} \rar F$ the character sending $T_n$ (resp. $\diam{m}$) to $a_n(f)$ (resp. $\chi(m)$) if $n\geq 1$ is coprime to $p$ (resp. $m \in \F_p^{\times}$). It is clear that $\mu_f$ determines $f$ and that its co-restriction to $\Q(f)$ is a ring epimorphism. Let now $\mu_{f,1}$ denote the restriction of $\mu_f$ to $\mathbb{T}_1$; the item \ref{dttc-2} in the proof of Proposition \ref{decomp-tate-twisted-general} shows that for $f,f' \in \mathcal{N}$, one has $\mu_{f,1} = \mu_{f',1}$ if and only if $f$ and $f'$ are twists by a Dirichlet character modulo $p$. 

Since $\mathbb{T}_1 \subset \mathbb{T}$, it is reduced, finite free over $\Z$ and $\mathbb{T}_1 \otimes \Q$ is a product of fields. Moreover, the intersection of the $\ker{\mu_{f,1}}$ over a set of representatives of $\mathcal{N}$ modulo torsion by Dirichlet characters modulo $p$ is trivial, because it is contained in the intersection of the kernels of the $\mu_f$ over $f \in \mathcal{N}$.  

For $f \in \mathcal{N}$, let $I_{G_{\Q} \cdot f,1} = \ker{\mu_{f,1}}$; the notation is obviously well-defined, i.e. $\ker{\mu_{f,1}}$ only depends on the $G_{\Q}$-orbit of $f$. 

\lem[igqf-determines-f-galtwist]{Let $f,g \in \mathcal{N}$. Then $I_{G_{\Q}\cdot f,1}=I_{G_{\Q}\cdot g,1}$ if and only if there is a Dirichlet character $\chi$ such that $f$ and $g\otimes \chi$ are Galois-conjugates. }

\demo{The ``if'' statement is a direct consequence of the preceding discussion. Assume conversely that $I_{G_{\Q}\cdot f,1}=I_{G_{\Q}\cdot g,1}$. Then let $\Q(f)_1\subset \Q(f)$, $\Q(g)_1 \subset \Q(g)$ be the fields generated by the images of $\mu_{f,1}$ and $\mu_{g,1}$, there is a homomorphism $\kappa: \Q(f)_1 \rar \Q(g)_1$ such that $\kappa \circ \mu_{f,1}=\mu_{g,1}$, i.e. $\kappa(a_n(f))=a_n(g)$ for $n \equiv 1\mod{p}$. Now, $\kappa$ lifts to an automorphism $\hat{\kappa}$ of $\C$, so the newform $f_1=\hat{\kappa}(f)$ is Galois-conjugate to $f$ (hence in $\mathcal{N}$) and such that $a_n(f_1)=a_n(g)$ if $n \equiv 1\mod{p}$. In particular, $\mu_{f_1,1}=\mu_{g,1}$, so $f_1$ and $g$ are twists of each other by a Dirichlet character modulo $p$.
}

\defi{Let $\mathcal{Z}_N \subset \mathbb{T}[N]_1$ be the subring of elements that commute to every $g \in N$. For $f \in \mathcal{N}$ and $\psi \in \mathscr{R}_1(f)$ (resp. $\psi \in \mathscr{R}_2(f)$), we can consider the $\mathbb{T}[N] \otimes F$-module $M_{f,\psi}$, where $T_n$ acts by $a_n(f)$ (for $n \geq 1$ prime to $p$), and $N$ acts by $\psi$ (resp. as the representation $\mrm{Ind}_C^N{\psi}$). By assumption, $\mathcal{Z}_N$ acts on $M_{f,\psi}$ by $F[N]$-automorphisms, so by a ring homomorphism $\mu_{f,\psi}: \mathcal{Z}_N \rar F$ which we call the \emph{central character} of $M_{f,\psi}$. }

The main result of this paragraph is that the \emph{kernel} of the central character determines $(f,\psi)$ up to some equivalence relation that we now describe. 

\defi{Let $f' \in \mathcal{N}$ and $\psi' \in \mathscr{R}_1(f') \cup \mathscr{R}_2(f')$. 

We say that the two couples $(f,\psi)$ and $(f',\psi')$ are \emph{equivalent} if the following conditions are satisfied:
\begin{itemize}[noitemsep,label=\tiny$\bullet$]
\item $\psi' \in \mathscr{R}_1(f')$ if and only if $\psi \in \mathscr{R}_1(f)$,
\item there exists a Dirichlet character $\chi$ modulo $p$ and an automorphism $\sigma$ of $\C$ such that $\sigma(f)\otimes \chi=f'$, and $(\sigma\circ\psi)\cdot (\chi\circ\det)$ has the same image in $\mathscr{R}_1(f')_{rep} \cup \mathscr{R}_2(f')_{rep}$ as $\psi'$. 
\end{itemize}}

Note that this is indeed an equivalence relation.  

\prop[central-character-determines-all]{Let $f,f' \in \mathcal{N}$, and $\psi \in \mathscr{R}_i(f), \psi' \in \mathscr{R}_j(f')$. The kernels of $\mu_{f,\psi}$ and $\mu_{f',\psi'}$ are equal if and only if $(f,\psi)$ and $(f',\psi')$ are equivalent. Moreover, if $\ker{\mu_{f,\psi}}\neq \ker{\mu_{f',\psi'}}$, their sum has finite index in $\mathcal{Z}_N$. }

\demo{Because $\mathcal{Z}_N$ is a finite free commutative $\Z$-algebra, the kernels of its ring homomorphisms to $\C$ are exactly its minimal prime ideals, which identify with the maximal ideals of $\mathcal{Z}_N \otimes \Q$, which implies the final claim. 

Our first step is to show the ``if'' of the equivalence.

If $\sigma \in \mrm{Gal}(F/\Q)$, then one easily sees that $\mu_{\sigma(f),\sigma \circ \psi} = \sigma \circ \mu_{f,\psi}$, so $\mu_{\sigma(f),\sigma\circ \psi}$ and $\mu_{f,\psi}$ have the same kernel. It is also not difficult to see that if $\psi \in \mathscr{R}_i(f)$, the kernel of $\mu_{f,\psi}$ only depends on the image of $\psi$ in $\mathscr{R}_i(f)_{rep}$ (the only non-trivial case being when $f$ has complex multiplication). Finally, if $\chi$ is a Dirichlet character, then there exists an $F$-isomorphism $\iota: M_{f,\psi} \rar M_{f\otimes \chi,\psi\cdot( \chi\circ\det)}$ such that $T_n\iota(x) = \chi(n)\iota(T_n x)$ (for $n\geq 1$ prime to $p$) and $g \cdot \iota(x) = \chi(\det{g})\iota(g \cdot x)$, so $\iota$ is a $\mathbb{T}[N]_1$-isomorphism, whence $\mu_{f\otimes\chi,\psi\cdot (\chi\circ\det)} = \mu_{f,\psi}$, and we are done. 

We now prove the ``only if'' direction. 

\emph{Step 1: We may assume that $f=f' \in \mathscr{S}\cup\mathscr{P}\cup\mathscr{C}$, $\mu_{f,\psi}=\mu_{f',\psi'}$ and $i=j$.}

Indeed, if $\ker{\mu_{f,\psi}}=\ker{\mu_{f',\psi'}}$, since $F/\Q$ is finite Galois, there exists $\sigma \in \mrm{Gal}(F/\Q)$ such that $\mu_{f',\psi'} = \sigma \circ \mu_{f,\psi}$, so, after replacing $(f,\psi)$ with $(\sigma(f),\sigma \circ \psi)$, we may assume that $\mu_{f,\psi}=\mu_{f',\psi'}$. Then $\mu_{f,1}=\mu_{f',1}$, so $f'=f \otimes \chi$, where $\chi$ is a Dirichlet character modulo $p$. Now $\mu_{f',\psi'} = \mu_{f,\psi} = \mu_{f\otimes \chi,\psi \cdot (\chi \circ \det)}=\mu_{f',\psi \cdot (\chi \circ \det)}$, so we are brought back to the case $f=f'$ and $\mu_{f,\psi}=\mu_{f,\psi'}$. After twisting $f$ one final time, we may assume $f \in \mathscr{S}\cup\mathscr{P}\cup\mathscr{C}$ by Lemma \ref{spcm-quadrichotomy}. 

Let $g$ be a generator of the cyclic group $C \cap \SL{\F_p}$, then $g+g^p \in \mathcal{Z}_N$ and one checks easily that $\mu_{f,\psi}(g+g^p)$ is $2$ if $i=1$ and $\omega+\omega^{-1}$ (where $\omega \neq 1$ is a root of unity) if $i=2$. Therefore one has $i=j$.

\emph{Step 2: Proof of the result if $f$ does not have complex multiplication.}

Let $g \in C$. By \cite[Cor. 2 to Th. 15]{Serre-Chebotarev} and the prime number theorem for arithmetic progressions \cite[XV, Theorems 4 and 5]{Lang-ANT}, there exists a prime $q$ such that $a_q(f) \neq 0$ and $q\det{g} \equiv 1\mod{p}$. It is then easy to see that $\frac{2}{i}\mu_{f,\psi}(T_q(g+g^p)) = a_q(f)\mrm{Tr}(g \mid M_{f,\psi})$, whence $\mrm{Tr}(g \mid M_{f,\psi})=\mrm{Tr}(g \mid M_{f,\psi'})$, so $\psi=\psi'$ if $i=1$. For $i=2$, by comparing traces, we see that $\mrm{Ind}_C^N{\psi}\simeq \mrm{Ind}_C^N{\psi'}$, which implies that $\psi,\psi'$ have the same image in $\mathscr{R}_2(f)_{rep}$. 

\emph{Step 3: Proof of the result if $f$ has complex multiplication.}

Since $\mathscr{R}_1(f)_{rep}$ has a single element, we may assume $i=2$. It is enough to show that $M_{f,\psi'}$ and $M_{f,\psi}$ have the same traces when restricted to $N \cap \F_p^{\times}\SL{\F_p}$, so, as in Step 2, we are done if we can show that, for any square $u \in \F_p^{\times}$, there is a prime $q \equiv u \mod{p}$ such that $a_q(f) \neq 0$. Now, by Dirichlet's theorem, there are infinitely many primes $q \equiv u \mod{p}$; if $q$ is such a prime which is unramified in $F$, then $q$ splits in $K$ and the two complex roots of $X^2-a_q(f)X+q$ lie in $F$. Therefore $\sqrt{-q}$ is not a root of this polynomial, so $a_q(f) \neq 0$ and the conclusion follows.} 

\lem[mfkzn-behaves-well]{The endomorphism of $\prod_{n \in \F_p^{\times}}{J_{\Gamma}^{\alpha^n}(p)_{\overline{K}}}$ induced by any $t \in \mathcal{Z}_N$ is contained in $\prod_{n \in \F_p^{\times}}{\mrm{End}_K(J_{\Gamma}^{\alpha^n}(p))}$.
}

\demo{Because $t \in \mathbb{T}[\GL{\F_p}]_1$, the image of $t$ lands in $\prod_{n \in \F_p^{\times}}{\mrm{End}_{\overline{K}}(J_{\Gamma}^{\alpha^n}(p)_{\overline{K}})}$ by Lemma \ref{hecke-gl2-grading}. It is therefore enough to check by descent that $t$ commutes to the action of $G_K$, which is a direct calculation using Proposition \ref{semidirect-gl-gal}. }

\subsection{Construction of quotients}
\label{subsect:constructing-quotients}

Recall the following standard result, which can be deduced for instance from \cite[Exp. $\mrm{VI}_\mrm{A}$, $\mrm{VI}_\mrm{B}$]{SGA3} and \cite[Lemma 2.3.3.2]{Anantharaman} for quotients, and \cite[\S 1]{LLR} for the claims about differentials. 

\lem[quotient-of-abelian-variety]{Let $A$ be an abelian variety over a field $K$ endowed with a ring homomorphism $R \rar \mrm{End}(A)$, and let $I \subset R$ be a left ideal, which is a finitely generated abelian group. The quotient fppf-sheaf $A/IA$ is representable by an abelian variety $A'$. Moreover,
\begin{itemize}[noitemsep,label=\tiny$\bullet$]
\item $A'$ is endowed with an action of $R/I$,
\item if $f: A \rar A'$ is the obvious projection map, then it is smooth, surjective, commutes to $R$, and $(\ker{f})(\overline{K}) = IA(\overline{K})$,
\item $A'$ only depends on the saturation $I'$ of $IR'$, where $R'$ is the saturation of $R$ in $\mrm{End}(A)$, i.e. the inverse image of $I \otimes \Q$ under the inclusion $R \rar \mrm{End}(A) \otimes \Q$,
\item the formation of $A'$ commutes with base change,
\item if $K$ has characteristic zero, $f^{\ast}H^0(A',\Omega^1) \subset H^0(A,\Omega^1)$ is exactly the collection of differentials killed by $I$,
\item one has an $R \otimes \Z_{\ell}$-morphism $\Tate{\ell}{A}/I\Tate{\ell}{A} \rar \Tate{\ell}{A'}$, if $\ell$ is a prime invertible in $K$, which is surjective with finite kernel. 
\end{itemize}
}

Let now $K$ be a field of characteristic zero, $\Gamma$ be a $p$-torsion group over $K$ equipped with a Weil pairing $\alpha$. Let $P,Q \in \Gamma(\overline{K})$ form a basis of $\Gamma(\overline{K})$. Let $\rho: G_K \rar \GL{\F_p}$ be the action of Galois on $\Gamma(\overline{K})$ in the basis $(P,Q)$ and let $\rho^{\ast}$ be the contragredient.

This defines an action of $\mathbb{T}[\GL{\F_p}]$ on $\prod_{n \in \F_p^{\times}}{J_{\Gamma}^{\alpha^n}(p)_{\overline{K}}} = J_{\Gamma}(p)_{\overline{K}}$. We will use Lemma \ref{quotient-of-abelian-variety} to construct quotients of $J_{\Gamma}^{\alpha}(p)$.

The quotients that we will consider are the following: 

\prop[quotient-by-T1]{Let $\omega \subset \mathcal{N}$ be a collection of newforms which is stable under twists by Dirichlet characters modulo $p$ and the action of $G_{\Q}$ on Fourier coefficients. Let $I_{\omega,1}$ be the saturation of $\prod_{f \in \omega}{\ker{\mu_{f,1}}}$ in $\mathbb{T}_1$. Then 
\begin{itemize}[noitemsep,label=\tiny$\bullet$]
\item the quotient abelian variety $J_{\Gamma,\omega}^{\alpha} = J_{\Gamma}^{\alpha}(p)/I_{\omega,1}J_{\Gamma}^{\alpha}(p)$ is well-defined, has an action of $\mathbb{T}[\GL{\F_p}]_1$, and its formation commutes with base change,
\item the inverse image of $H^0(J_{\Gamma,\omega}^{\alpha},\Omega^1)$ in $H^0(J_{\Gamma}^{\alpha}(p),\Omega^1)$ is exactly the collection of global differentials annihilated by $I_{\omega,1}$,
\item if $K=\Q$, for any maximal ideal $\mfk{q}\subset \OO_F$ with residue characteristic $q$, through the decomposition of Proposition \ref{decomp-tate-twisted-general}, $\Tate{q}{J_{\Gamma,\omega}^{\alpha}} \otimes_{\Z_{q}} F_{\mfk{q}}$ identifies as a $\mathbb{T}[G'_{\Q}] \otimes F_{\mfk{q}}$-module with the quotient of $\Tate{q}{J_{\Gamma}^{\alpha}(p)} \otimes_{\Z_q} F_{\mfk{q}}$ by the sum of the factors corresponding to $f \notin \omega$.
\end{itemize}}

\demo{This is a consequence of Propositions \ref{quotient-of-abelian-variety} and \ref{decomp-tate-twisted-general} because for $f,f' \in \mathcal{N}$, either $\mu_{f,1}$ and $\mu_{f',1}$ have the same kernel, and $f'$ is the twist by a Dirichlet character modulo $p$ of some newform in the $G_{\Q}$-orbit of $f$, or $\ker{\mu_{f,1}}+\ker{\mu_{f',1}}$ is an ideal of $\mathbb{T}_1$ with finite index. }

We also consider quotients by ideals of $\mathcal{Z}_N$. 

\prop[quotient-by-ZN]{Assume that the image of $\rho^{\ast}$ is contained in the normalizer $N$ of some non-split Cartan subgroup $C$. Let $\omega \subset \bigcup_{f \in \mathcal{N}}{\{f\} \times \left(\mathscr{R}_1(f) \cup \mathscr{R}_2(f)_{ind}\right)}$ be saturated under the equivalence relation defined above Proposition \ref{central-character-determines-all}. Let $I_{\omega}$ be the saturation in $\mathcal{Z}_N$ of $\prod_{(f,\psi) \in \omega}{\ker{\mu_{f,\psi}}}$. Then
\begin{itemize}[noitemsep,label=$-$]
\item the quotient abelian variety $J_{\Gamma,\omega}^{\alpha} := J_{\Gamma}^{\alpha}(p)/I_{\omega}J_{\Gamma}^{\alpha}(p)$ is well-defined, has an action of $\mathbb{T}[N]_1$, and its formation commutes with base change and isomorphisms of $p$-torsion groups,
\item the inverse image of $H^0(J_{\Gamma,\omega}^{\alpha},\Omega^1)$ in $H^0(J_{\Gamma}^{\alpha}(p),\Omega^1)$ is exactly the collection of global differentials annihilated by $I_{\omega}$,
\item if $K=\Q$, for any maximal ideal $\mfk{q}\subset \OO_F$ with residue characteristic $q$, through the decomposition of Proposition \ref{tate-module-twisted-cartan}, $\Tate{q}{J_{\Gamma,\omega}^{\alpha}} \otimes_{\Z_{q}} F_{\mfk{q}}$ is the quotient of $\Tate{q}{J_{\Gamma}^{\alpha}(p)} \otimes_{\Z_q} F_{\mfk{q}}$ by all the factors attached to couples $(f,\psi) \notin \omega$. 
\end{itemize}}

\demo{This follows again from Lemma \ref{quotient-of-abelian-variety} and Proposition \ref{central-character-determines-all}. }

\bigskip

We will mainly apply Proposition \ref{quotient-by-T1} in the following situation: \\
\textbf{Situation}\refstepcounter{propo} \label{quotient-by-T1S}(\thepropo):\; $\omega$ is the set of elements of $\mathcal{N}$ which are twists by Dirichlet characters of newforms in a collection $\omega' \subset \mathscr{S}$ of newforms which is stable by the action of $G_{\Q}$ on Fourier expansions.\\

We will mainly apply Proposition \ref{quotient-by-ZN} in the following situation:\\
\textbf{Situation}\refstepcounter{propo} \label{quotient-by-3C}(\thepropo):\; One has $p \equiv -1 \mod{3}$, so that a given character $\psi: C \rar F^{\times}$ of order exactly $3$ is in $\mathscr{R}_2(f)$ for every $f \in \mathscr{S}$. We select a collection $\omega_0 \subset \mathscr{S}$ of newforms which is stable under the action of $G_{\Q}$ on Fourier expansions, and take $\omega$ be the saturation of $\omega_0 \times \{\psi\}$. We denote the attached abelian variety by $B_{\omega_0,3}$.

\subsection{Bloch--Kato conjecture with coefficients and finiteness of Selmer group}
\label{subsect:bloch-kato-rank-zero}

The goal of the remainder of Section \ref{sect:cartan-bsd} is (mainly) to prove the finiteness of the Mordell--Weil group of the abelian variety $B_{\omega,3}$ constructed in Situation \ref{quotient-by-3C} provided that its $L$-function does not vanish. Doing so involves using known cases of the Bloch--Kato conjecture in rank zero to show that a certain $\ell$-adic \Nekovar cohomology group $\tilde{H}^2$ (with suitable local conditions) is finite. Applying duality, it follows that a suitable localization of the divisible Selmer group for $B_{\omega,3}$ is finite, which implies finiteness of the Mordell--Weil group by classical Kummer theory. 

Since we use the version of the Bloch--Kato conjecture proved in \cite{KLZ15}, we need certain ordinarity assumptions, which are not ``rational'': if $f \in \mathcal{N}$ is ordinary at a certain prime $\mfk{p}$ of $\Qbar \subset \C$, its Galois conjugates need not be ordinary at $\mfk{p}$. Yet, if $f$ ``contributes'' to $B_{\omega,3}$, so do its Galois conjugates. This is why we need to work by adding extra coefficients. Since we also work with more general abelian varieties instead of elliptic curves, we also need to involve a polarization of $B_{\omega,3}$. These modifications seem standard and dealing with them does not involve much extra work; however, since we are unaware of a reference discussing them, we spell out the complete argument here, which should be essentially viewed as a slight generalization of e.g. \cite[Theorem 11.7.4]{KLZ15}. It is a little bit more general than is warranted by our set-up, but we found that this level of generality (which can serve for other quotients of $J_{\Gamma}^{\alpha}(p)$) was easier to follow. 

We use the Selmer complexes first defined in \cite{Nekovar-Selmer}, and refer the reader to \cite[\S 11.2]{KLZ15} for their definition and important properties. We need to make the following modifications: 
\begin{enumerate}[noitemsep,label=$(\arabic*)$]
\item Our ring of coefficients $R$ is the ring of integers of some finite extension $F/\Q_p$, where $p$ is an odd prime. The maximal ideal is $\mfk{m}$. If $M$ is a $R$-module, and $d_F \in F$ is the element with the lowest valuation such that $\mrm{Tr}_{F/\Q_p}(d_FR) \subset \Z_p$, then 
\[u \in \mrm{Hom}(M,F/R) \mapsto \mrm{Tr}_{F/\Q_p}\circ (d_F u) \in \mrm{Hom}_{\Z_p}(M,\Q_p/\Z_p)\] is an isomorphism of $R$-modules, so it is equivalent to consider the Pontryagin dual of $M$ as a $R$-module and as a $\Z_p$-module. 
\item We need to work with \emph{admissible} $R$-modules\footnote{This was probably also needed in \emph{op.cit.} as well to formulate \Nekovar duality.} as defined in \cite[Chapter 3.2]{Nekovar-Selmer}: if $G$ is a topological group acting linearly on a $R$-module $M$, $M$ is admissible if the image $I$ of $G$ in $\mrm{End}_R(M)$ is finitely generated and the map $G \rar I$ is continuous. The notion of admissibility behaves well with respect to the usual operations (sub-objects, quotients, tensor products, homomorphism modules). When $R$ is finitely generated (resp. co-finitely generated), $M$ is admissible if and only if $G$ acts continuously in the $\mfk{m}$-adic (resp. discrete) topology. 
\item We slightly change the definition of simple local condition in \cite[Definition 11.2.2]{KLZ15}: if $v$ is a place of $\Q$ and $M$ is an admissible $R[G_{\Q}]$-module, the \emph{simple} local condition at $v$ attached to a given submodule $N \subset H^1(\Q_v,M)$ is the subcomplex 
\[[0 \rar \underbrace{C^0(G_{\Q_v},M)}_{\deg=0} \rar \mrm{pr}^{-1}(N) \rar 0]\] of $C^{\bullet}(G_{\Q_v},M)$, where $\mrm{pr}^{-1}(N)$ is the set of cocycles in $C^{1}(G_{\Q_v},M)$ whose attached cohomology class lies in $N$.  
\item With this definition, the \Nekovar duality theorem \cite[Theorem 11.2.7]{KLZ15} for simple local conditions is functorial in the following sense: if $f: M \rar M'$ is a morphism of admissible $R[G_{\Q}]$-modules, that respects the simple local conditions $\Delta_{M},\Delta_{M'}$, then the maps 
\begin{align*}
&f: \tilde{H}^i(G_{\Q,S},M;\Delta_M) &&\rar \tilde{H}^i(G_{\Q,S},M';\Delta_{M'}),\\ 
&f^{\vee}: \tilde{H}^{3-i}(G_{\Q,S},(M')^{\vee}(1);\Delta_{M'}^{\perp}) &&\rar \tilde{H}^{3-i}(G_{\Q,S},M^{\vee}(1);\Delta_{M}^{\perp})
\end{align*} are adjoint of each other for the \Nekovar duality pairing (this is a consequence of \cite[Theorem 6.4.3, Proposition 6.4.2, \S 6.4.4]{Nekovar-Selmer}).    
\end{enumerate}

We consider, in the rest of this section, the following simple local conditions:
\begin{itemize}[noitemsep,label=\tiny$\bullet$]
\item the Bloch--Kato condition when $v \nmid p$ and $M$ is finite free over $R$: it is attached to the inverse image $H^1_f(G_{K_v},M)$ of $H^1(G_{K_v}/I_{K_v},(M \otimes \Q_p)^{I_{K_v}})$ under $H^1(G_{K_v},M) \rar H^1(G_{K_v},M \otimes \Q_p)$,
\item the Bloch--Kato condition when $v \nmid p$ and $M=M' \otimes \Q_p/\Z_p$ where $M'$ is finite free over $R$: it is attached to the image $H^1_f(G_{K_v},M)$ of $H^1(G_{K_v}/I_v,(M' \otimes \Q_p)^{I_{K_v}})$ in $H^1(G_{K_v},M)$,
\item the Bloch--Kato condition when $v \mid p$, $M$ is finite free over $R$, and $M \otimes \Q_p$ is de Rham: it is attached to the inverse image $H^1_f(G_{K_v},M)$ under $H^1(G_{K_v},M) \rar H^1(G_{K_v},M \otimes \Q_p)$ of the subspace of crystalline classes $H^1_f(G_{K_v},M\otimes \Q_p) \subset H^1(G_{K_v},M\otimes \Q_p)$ defined in \cite[(3.7.2)]{BKTama},  
\item the Bloch--Kato condition when $v \mid p$ and $M=M' \otimes \Q_p/\Z_p$ where $M'$ is finite free over $R$ and $M' \otimes \Q_p$ is de Rham at $p$: it is attached to the image $H^1_f(G_{K_v},M)$ by $H^1(G_{K_v},M' \otimes \Q_p) \rar H^1(G_{K_v},M)$ of $H^1_f(G_{K_v},M' \otimes \Q_p)$. 
\item the Kummer local condition at a place $v$, and $M$ is the $p^{\infty}$-torsion of some abelian variety $A/K$: it is the image of the connecting map $A(K) \otimes \Q_p/\Z_p \rar H^1(G_{K_v},A(K)[p^{\infty}])$ in the Kummer exact sequence (see e.g. the introduction of \cite{Greenberg-Coates}).
\end{itemize}

The Bloch--Kato condition is clearly functorial; the Kummer local condition is functorial with respect to morphisms of abelian varieties. 

\prop[local-tate-duality-bloch-kato]{Let $M$ be a finite free $\Z_p$-module with a continuous action of $G_K$ for some finite extension $K/\Q_p$. The subspaces $H^1_f(G_K,M)$ and $H^1_f(G_K,M^{\vee}(1))$ are exact duals of each other under the perfect local Tate duality pairing $H^1(G_K,M) \times H^1(G_K,M^{\vee}(1)) \rar \Q_p/\Z_p$. 
}

\demo{This is \cite[Proposition 3.8]{BKTama}. }

\prop[nekovar-is-selmer]{Let $A/\Q$ be an abelian variety and $q$ be an odd prime. Then one has an exact sequence $0 \rar \Tate{q}{A} \rar \Tate{q}{A} \otimes_{\Z_q} \Q_q \rar A[q^{\infty}] \rar 0$ of $G_{\Q}$-modules, and $\Tate{q}{A} \otimes_{\Z_q} \Q_q$ is a de Rham representation of $G_{\Q_q}$. Let $S$ be a finite set of places of $\Q$ containing $q$, $\infty$, and the places of bad reduction for $A$, so that all our representations are admissible $\Z_q[G_{\Q,S}]$-modules. Then, if $\Delta^{BK}$ denotes the collection of Bloch--Kato conditions on $A[q^{\infty}]$ previously defined, $\tilde{H}^1(\Z[1/qS],A[q^{\infty}];\Delta^{BK})$ identifies with the classical $q^{\infty}$-Selmer group of $A$.     }

\demo{Let $A^{\vee}/\Q$ be the dual abelian variety to $A$, so that $A$ is dual to $A^{\vee}$ by \cite[Proposition 9.5]{MilAb}. Because $\Qbar^{\times}=H^0_{\text{\'et}}(A^{\vee}_{\Qbar},\mathbb{G}_m)$ is a divisible group, and since the N\'eron--Severi group of $A^{\vee}$ is torsion-free by \cite[Corollary 12.8]{MilAb}, the Kummer exact sequence yields an isomorphism of $\Z_q[G_{\Q}]$-modules $H^1_{\text{\'et}}(A^{\vee}_{\Qbar},\Z_q)(1) \simeq \Tate{q}{A}$. By classical results (e.g. \cite[Example 3.11]{BKTama}), $\Tate{q}{A} \otimes \Q_q$ is a de Rham $\Q_q[G_{\Q_q}]$-module, so the Bloch--Kato local condition is well-defined at all places. By the definition of the \Nekovar cohomology groups, it is then enough to show that: 
\begin{itemize}[noitemsep,label=$-$]
\item the images in $H^1(G_{\Q_q},A[q^{\infty}])$ of $H^1_f(G_{\Q_q},\Tate{q}{A}\otimes \Q_q)$ and $A(\Q_q) \otimes \Q_q/\Z_q$ (the former under the obvious map, the latter under the Kummer map) are equal,
\item for every finite place $v \neq q$, the Kummer map $A(\Q_v) \otimes \Q_q/\Z_q \rar H^1(G_{\Q_v},A[q^{\infty}])$ is zero,
\item for every finite place $v \neq q$, the map $H^1(G_{\Q_v}/I_{\Q_v},(\Tate{q}{A} \otimes \Q_q)^{I_{\Q_v}}) \rar H^1(G_{\Q_v},A[q^{\infty}])$ is zero.
\end{itemize}

All three items are quite classical. The first point is the most subtle one, and is discussed in \cite[Example 3.11]{BKTama}. The second point is e.g. \cite[Proposition 4.1]{Greenberg-Coates}. We were not able to find a reference for the third point, so we include the argument here. 

The group $H^1(G_{\Q_v}/I_{\Q_v},(\Tate{q}{A} \otimes \Q_q)^{I_{\Q_v}})$ is $q$-divisible, so it is enough to show that the $q^{\infty}$-torsion group $H^1(\mrm{Gal}(\Q_v^{nr}/\Q_v),A[q^{\infty}](\Q_v^{nr}))$ has finite exponent, i.e. that the exponent of $H^1(\mrm{Gal}(\Q_v^{nr}/\Q_v),A[q^r](\Q_v^{nr}))$ is bounded as $r \rar \infty$. By \cite[Proposition 10.16]{Harari}, the cardinality of this group is exactly $A[q^r](\Q_v)$, and $A[q^{\infty}](\Q_v)$ is finite by \cite[Lemma I.3.3]{MilneADT}.}

\lem[localize-duality-finiteness]{Let $A$ be a finite free $\Z_q$-algebra and $M$ be a finitely generated $A$-module. Then $\mrm{Hom}_{\Z_q}(M,\Q_q/\Z_q)$ can be made into an $A$-module by $a \cdot f := f(a \cdot)$. Let $\mfk{q}$ be a maximal ideal of $A$. Then $M_{\mfk{q}}$ is finite iff $\mrm{Hom}_{\Z_q}(M,\Q_q/\Z_q)_{\mfk{q}}$ is finite. }

\demo{By definition, $A \rar \mrm{End}(\mrm{Hom}(M,\Q_q/\Z_q))$ always factors through $A \rar \mrm{End}(M)$, so, if $A \rar \mrm{End}(M)$ factors through some localization $A \rar A_{\mfk{q}}$ at a maximal ideal, so does $A \rar \mrm{End}(\mrm{Hom}(M,\Q_q/\Z_q))$. The rule $M \mapsto \mrm{Hom}_{\Z_q}(M,\Q_q/\Z_q)$ is (finitely) additive with respect to the $A$-module structure; $\Z_q$ is a complete Noetherian local ring, so $A$ is a product of its localizations at prime ideals, which are finite free local $\Z_q$-algebras: hence we are brought to the case where $A$ is local. 

In this case, we need to show that $M$ is finite if and only if $\mrm{Hom}(M,\Q_q/\Z_q)$ is finite, which is immediate (e.g. by the structure of finitely generated $\Z_q$-modules).}

\prop[simple-H2N-to-BSD]{Let $A$ be a simple abelian variety over $\Q$ and $N\geq 1$ be an integer divisible by all its primes of bad reduction. Let $Z$ be the center of $\mrm{End}(A)$, $\lambda: A \rar A^{\vee}$ be a polarization of degree $d$. Let $F$ be a number field and $\mfk{q} \subset \OO_F$ be a maximal ideal of odd residue characteristic $q \nmid dN$. Assume that $\Tate{q}{A} \otimes_{\Z_q} F_{\mfk{q}} \simeq \bigoplus_{i=1}^r{V_i^{\oplus n_i}}$, where every $V_i$ is an absolutely irreducible $F_{\mfk{q}}[G_{\Q}]$-module and the $V_i$ are pairwise non-isomorphic, and fix a $\OO_{F_{\mfk{q}}}[G_{\Q}]$-lattice $T_1 \subset V_1$.  

There exists a ring map $Z \rar \OO_{F_{\mfk{q}}}$ such that $\Tate{q}{A} \otimes_{Z_{\mfk{p}}} F_{\mfk{q}} \simeq V_1^{\oplus n_1}$, let $\mfk{p} \subset Z$ denote the inverse image of $\mfk{q}$ (which is a maximal ideal of $Z$), and $\mfk{p}'$ be the dual of $\mfk{p}$ by the Rosati involution. Let $\Delta^{BK}$ denotes the Bloch--Kato Selmer condition, and assume that $\tilde{H}^2(\Z[1/Nq], T_1; \Delta^{BK})$ is finite. Then $\mrm{Sel}_{q^{\infty}}(A/\Q)\otimes_Z Z_{\mfk{p}'}$ is finite. As a consequence, $A(\Q)$ is finite. 
}

\demo{\emph{Step 1: Preliminaries on abelian varieties.}

By \cite[Theorem 12.5]{MilAb}, $\mrm{End}(A)$ is a finite free $\Z$-module, hence its center $Z$ is a finite free $\Z$-module as well and $\mrm{End}(A)/Z$ is torsion-free. Moreover, $A$ is simple, so for any nonzero $u \in \mrm{End}(A)$, the connected component of unity in $\ker{u}$ is a smooth proper connected $\Q$-group scheme by Cartier's theorem, i.e. an abelian subvariety of $A$. By \cite[Proposition 8.1]{MilAb}, $u$ is an isogeny, so $u$ is surjective, hence the multiplication by $u$ is injective in $\mrm{End}(A)$. Therefore $\mrm{End}(A)$ has no zero divisors. By definition of $A^{\vee}$ (see \cite[\S 9--11]{MilAb}), $u \in \mrm{End}(A)\mapsto u^{\vee} \in \mrm{End}(A^{\vee})$ is an anti-isomorphism of rings. 

Because $q$ is prime to $d$, the Rosati involution $u \mapsto u^{\dagger}$ attached to $\lambda$ defines a ring involution of $Y := Z \otimes \Z_q$. For every $n \geq 1$, $Z/q^nZ$ acts faithfully on $A[q^n]$ (by \cite[Theorem 12.5]{MilAb}), so $\Tate{q}{A}, A[q^{\infty}]$ are both $Y[G_{\Q}]$-modules. By \cite[Theorem 11.1, \S 16]{MilAb}, $\lambda$ produces, for every $n \geq 1$, a perfect duality pairing $\langle \cdot ,\,\cdot \rangle: A[q^n](\Qbar) \times A[q^n](\Qbar) \rar \mu_{q^n}(\Qbar)$ which is naturally compatible with change of $n$ and such that for every $x,y \in A[q^n](\Qbar)$ and every $t \in Y$, one has $\langle t\cdot x,\,y\rangle=\langle x,\,t^{\dagger}\cdot y\rangle$.

\emph{Step 2: Applying Faltings's theorem}

Since the $V_i$ are absolutely irreducible and pairwise non-isomorphic, $\mrm{End}_{F_{\mfk{q}}[G_{\Q}]}(\Tate{q}{A} \otimes_{\Z_q} F_{\mfk{q}})$ identifies with the product of $\prod_{i=1}^r{\mathcal{M}_{n_i}(F_{\mfk{q}})}$ by Schur's lemma. By Faltings's theorem \cite[Theorem 4]{Faltings}, the following obvious map is an isomorphism \[\gamma: \mrm{End}(A) \otimes F_{\mfk{q}}\rar \mrm{End}_{F_{\mfk{q}}[G_{\Q}]}(\Tate{q}{A} \otimes_{\Z_q} F_{\mfk{q}}). \]

Therefore, $\gamma$ is an isomorphism of $Z \otimes F_{\mfk{q}}$ to the center of $\prod_{i=1}^r{\mathcal{M}_{n_i}(F_{\mfk{q}})}$, i.e. $\prod_{i=1}^r{F_{\mfk{q}}}$, and the $i$-th projection $\pi_i$ is by construction the action of $Y \otimes_{\Z_q} F_{\mfk{q}}$ on $V_i$. It follows that the $\ker{\pi_i}$ are pairwise distinct and are all the maximal ideals of $Y \otimes_{\Z_q} F_{\mfk{q}}$, and as a consequence, $\Tate{q}{A} \otimes_{Y,\pi_1} F_{\mfk{q}} \simeq V_1^{\oplus n_1}$. 

Since $Z$ is a finite free $\Z$-algebra which is a domain, we take $\mfk{p}$ as the inverse image of $\mfk{q}\OO_{F_{\mfk{q}}}$ under $\pi_1$, and let $\mfk{p}'$ denote its image under the Rosati involution. Note that $Y_{\mfk{p}}$ identifies with the completion of $Z$ at $\mfk{p}$, so it is a finite free local $\Z_q$-algebra and a domain.

\emph{Step 3: Finiteness of the Selmer group. }

We wish to prove here that $\mrm{Sel}_{q^{\infty}}(A/\Q) \otimes_Y Y_{\mfk{p}'}$ is finite (where the action is the obvious one). 

By Proposition \ref{local-tate-duality-bloch-kato}, the \Nekovar duality theorem, Proposition \ref{nekovar-is-selmer} and Step 1, we have an isomorphism \[\mrm{Sel}_{q^{\infty}}(A/\Q) \simeq \mrm{Hom}_{\Z_q}(\tilde{H}^2(\Z[1/Nq],\Tate{q}{A}; \Delta^{BK}),\Q_q/\Z_q),\] which is $Y$-linear if $t \in Y$ acts by $t^{\dagger}$ on the left-hand side (contrarily to the above paragraph) and pre-composition by $t$ on the right-hand side. Therefore, by Lemma \ref{localize-duality-finiteness}, it suffices to show that $\tilde{H}^2(\Z[1/Nq],\Tate{q}{A}; \Delta^{BK})\otimes_Y Y_{\mfk{p}}$ is finite. 

Note that any $V_1^{\oplus g}$ (for any $g \geq 1$) is, as a $\Q_q[G_{\Q}]$-module, a direct summand of a sum of copies of $\Tate{q}{A} \otimes \Q_q$, so by e.g. \cite[Th\'eor\`eme 3.8]{FontainePeriodes} it is a de Rham representation, so the Bloch--Kato local condition for any of its lattices (in particular, $T_1$) is well-defined.  

Since $Y$ is a finite free $\Z_q$-algebra, $Y_{\mfk{p}}$ is a direct summand of $Y$, so it is equivalent to show that $\tilde{H}^2(\Z[1/Nq],\Tate{q}{A} \otimes_Y Y_{\mfk{p}}; \Delta^{BK})$ is finite. By Step 2, there is an injective map of $\Z_q[G_{\Q}]$-modules $\iota: (\Tate{q}{A} \otimes_Y Y_{\mfk{p}})^{\oplus s} \rar T_1^{\oplus n_1}$ with finite cokernel (where $s$ is the generic rank of $\OO_{F_{\mfk{q}}}$ over the domain $Y_{\mfk{p}}$) killed by $q^{\alpha}$ for some $\alpha \geq 1$. After carefully unpacking the definitions of the Selmer complex, \cite[1.3.11]{Nekovar-Selmer} implies that one has an exact triangle \[\tilde{R\Gamma}(\Z[1/Nq], (\Tate{q}{A} \otimes_Y Y_{\mfk{p}})^{\oplus s};\Delta^{BK}) \overset{\iota}{\rar} \tilde{R\Gamma}(\Z[1/Nq],T_1^{\oplus n_1};\Delta^{BK}) \rar C \rar \ldots\] where $C$ is a complex of $\Z_q/q^{\alpha}\Z_q$-modules. Hence the maps on Selmer cohomology groups induced by $\iota$ have finite kernels and cokernels. 

Since $\tilde{H}^2(\Z[1/Nq],T_1;\Delta^{BK})$ is finite by assumption, so is $\tilde{H}^2(\Z[1/Nq],T_1^{\oplus n_1};\Delta^{BK})$ (because the Bloch--Kato condition is additive, hence so is the Selmer cohomology group), and therefore $\tilde{H}^2(\Z[1/Nq],\Tate{q}{A} \otimes_Y Y_{\mfk{p}}; \Delta^{BK})$ is finite as well.  

\emph{Step 4: Conclusion.}

By Kummer theory, we have an injection of $Y$-modules $A(\Q) \otimes \Q_q/\Z_q \rar \mrm{Sel}_{q^{\infty}}(A/\Q)$. By Step 3 (and since localization is exact), $\left(A(\Q) \otimes \Q_q/\Z_q\right)_{\mfk{p}'}$ is finite, i.e. $\left(A(\Q)\otimes_Z Y_{\mfk{p}'}\right) \otimes_{\Z_q} \Q_q/\Z_q$ is finite. Since $A(\Q)$ is a finitely generated abelian group, $A(\Q)\otimes_Z Y_{\mfk{p}'}$ is a finitely generated $\Z_q$-module, so it is finite. 

As we saw in Step 1, $Z$ is a finite free $\Z$-algebra and is a domain, so there is an injective homomorphism $\pi: Z^{\oplus \rho} \rar A(\Q)$ with finite cokernel. Since $Y_{\mfk{p}'}$ is flat over $Z$, $\pi$ induces an injective morphism $Y_{\mfk{p}'}^{\oplus \rho} \rar A(\Q) \otimes_{Z} Y_{\mfk{p}'}$, whence $\rho=0$ and we are done. }

\prop[H2N-to-BSD]{Let $A/\Q$ be an abelian variety and $F$ be a number field. Let $S$ be an infinite set of finite places of $F$, and $(V_{i,v})_{v \in S}$ be a compatible system of representations of $G_{\Q}$ over $F$ supported in $S$ for every $1 \leq i \leq r$. Suppose furthermore that for every $v \in S$ with residue characteristic $p$, 
\begin{itemize}[noitemsep,label=\tiny$\bullet$]
\item $\Tate{p}{A} \otimes_{\Z_p} F_v \simeq \bigoplus_{i=1}^r{V_{i,v}^{\oplus n_{i,v}}}$ for some positive integers $n_{1,v},\ldots,n_{r,v}$,
\item for every $1 \leq i \leq r$, $V_{i,v}$ is absolutely irreducible,
\item for every $1 \leq i < j \leq r$, $V_{i,v}$ is not isomorphic to $V_{j,v}$, 
\end{itemize}
Then, for every $1 \leq i \leq r$, the integer $n_{i,v}$ does not depend on $v$. 
Let $N \geq 1$ be an integer which is divisible by all the primes of bad reduction of $A$. Assume that, for every $1 \leq i \leq r$, there are infinitely many $v \in S$ such that for some lattice $T_{i,v}\subset V_{i,v}$, $\tilde{H}^2(\Z[1/Nq],T_{i,v};\Delta^{BK})$ is finite. Then $A(\Q)$ is finite. 
}

\demo{For any $v \in S$ with residue characteristic $p$ and any $1 \leq i \leq r$, $V_{i,v}$ is, as a $\Q_p$-vector space, a direct factor of a sum of copies of $\Tate{p}{A}\otimes_{\Z_p}\Q_p$. Therefore, by the N\'eron--Ogg--Shafarevich criterion \cite[Theorem 1]{GoodRed}, $V_{i,v}$ is unramified at every place of good reduction for $A$ not dividing $p$. Let $N$ be some integer that is divisible by every prime of bad reduction for $A$. 

Since every $V_{i,v}$ is irreducible and the $V_{i,v}$ and $V_{j,v}$ are pairwise non-isomorphic when $j \neq i$, the $F_v$-vector space $\mrm{Hom}_{F_v[G_{\Q}]}(V_{i,v},V_{j,v})$ is zero if $i \neq j$ and $F_v$ if $i=j$. 

By Proposition \ref{simple-H2N-to-BSD}, and since $A$ is isogenous to a product of powers of certain of its simple abelian subvarieties, it is enough to show that the $n_{i,v}$ do not depend on $v$ and that, for every simple abelian subvariety $B/\Q$ of $A$ and every $v \in S$ with residue characteristic $p$, there are integers $m_{1,v},\ldots,m_{r,v}$ such that \[\Tate{p}{B} \otimes_{\Z_p} F_v \simeq \bigoplus_{i=1}^r{V_{i,v}^{\oplus m_{i,v}}}.\]

The second claim is elementary representation theory: if $B$ is an abelian subvariety of $A$ and $v \in S$ has residue characteristic $p$, then $\Tate{p}{B}\otimes_{\Z_p} F_v$ is isomorphic to a $F_v[G_{\Q}]$-submodule of $\bigoplus_{i=1}^r{V_{i,v}^{\oplus n_{i,v}}}$. In particular, it is semi-simple by \cite[(15.2)]{Curtis-Reiner} and its irreducible components are contained in the $V_{i,v}$, so we are done.

We now discuss the first claim. Let $v,v' \in S$ with residue characteristics $p,p'$. There is a finite set of primes $T$ containing $p,p'$ and all the bad reduction primes of $A$ such that for every prime $\ell \notin T$, $V_{i,v}$ and $V_{i,v'}$ are unramified at $\ell$ and one has
\begin{align*}
\sum_{i=1}^r{n_{i,v'}\mrm{Tr}(\Fr_{\ell}\mid V_{i,v})} &= \sum_{i=1}^r{n_{i,v'}\mrm{Tr}(\Fr_{\ell}\mid V_{i,v'})} = \mrm{Tr}(\Fr_{\ell} \mid \Tate{p'}{A} \otimes_{\Z_{p'}} F_{v'}) \\
&= \mrm{Tr}(\Fr_{\ell} \mid \Tate{p}{A}\otimes_{\Z_p} F_v) =  \sum_{i=1}^r{n_{i,v}\mrm{Tr}(\Fr_{\ell}\mid V_{i,v})}.
\end{align*} 

By Chebotarev's theorem, the $\Fr_{\ell}$ for $\ell \notin T$ form a dense subset of $G_{\Q,Np}$, so the continuous representations $\bigoplus_{i=1}^r{V_{i,v}^{\oplus n_{i,v}}}$ and $\bigoplus_{i=1}^r{V_{i,v}^{\oplus n_{i,v'}}}$ of $G_{\Q,Np}$ over $F_v$ have the same traces. By \cite[Corollary 27.13]{Curtis-Reiner}, since the $V_{i,v}$ are absolutely irreducible and pairwise non-isomorphic, one has $n_{i,v}=n_{i,v'}$ for every $i$.  }

\subsection{Quotients of $J_{\Gamma}^{\alpha}(p)$ of rank zero}
\label{subsect:rank-zero-bsd}

As previously, $F \subset \C$ is a number field Galois over $\Q$, containing the $p(p^2-1)$-th roots of unity, the Fourier coefficients of every $f \in \mathcal{N}$, and, for every $f \in \mathscr{C}_M$ and every prime $q$ split in $\Q(\sqrt{-p})$, the roots of $X^2-a_q(f)X+q$.  

We now assume that $\Gamma$ is a $p$-torsion group over $\Q$ endowed with a Weil pairing $\alpha$. Let $(P,Q)$ denote a basis of $\Gamma(\Qbar)$, let $\rho: G_{\Q} \rar \GL{\F_p}$ denote the matrix of the Galois action on $\Gamma(\Qbar)$ in the basis $(P,Q)$ and let $\rho^{\ast}$ be its contragredient. We assume that the image of $\rho^{\ast}$ is contained in the normalizer $N$ of some Cartan subgroup $C$, and furthermore that:
\begin{itemize}[noitemsep,label=\tiny$\bullet$]
\item $\rho^{\ast}: G_{\Q} \rar N$ is surjective,
\item $\Gamma_{\Q_p}$ is the $p$-torsion subgroup scheme of some elliptic curve over $\Q_p$ (i.e. $Y_{\Gamma}(p)(\Q_p) \neq \emptyset$).
\end{itemize}

\prop[BSD-that-interests-us]{Let $\omega$ be as in Proposition \ref{quotient-by-ZN} where furthermore
\begin{itemize}[noitemsep,label=\tiny$\bullet$]
\item $\omega$ does not contain any couple $(f,\psi)$ where $f \in \mathscr{C}_M$, or $\psi \in \mathscr{R}_1(f)$, or $f \in \mathscr{C}$ and $\psi^4=1$,
\item there is a subset $\omega' \subset \omega$ such that every element of $\omega$ is equivalent to some element of $\omega'$ and, for every $(f,\psi) \in \omega'$, if $g \in \mathcal{S}_1(\Gamma_1(D))$ denotes the normalized newform attached to the odd irreducible dihedral Artin representation $\mrm{Ind}_K^{\Q}{\psi\circ\rho^{\ast}_{|G_K}}$, the central Rankin--Selberg $L$-value $\Lambda(f,g,1)$ (cf. Proposition \ref{rankin-selberg-galois}) does not vanish.  
\end{itemize}    

Let $J_{\Gamma,\omega}^{\alpha}/\Q$ be the Abelian variety constructed in Proposition \ref{quotient-by-ZN}. Then $J_{\Gamma,\omega}^{\alpha}(\Q)$ is finite. }

\rem{We stress that the abelian variety $J_{\Gamma,\omega}^{\alpha}$ appearing in Proposition \ref{BSD-that-interests-us} is not of $\GL$-type. }

\rem{It is apparent from the proof that the conditions that $\rho^{\ast}$ is surjective or that $\psi \notin \mathscr{R}_1(f)$ are for convenience. The only difference is that we would be forced to deal with factors that are representations attached to classical modular forms, and we would need to use Kato's Euler system \cite{KatoBSD} rather than the Euler system of Beilinson--Flach elements. We do not do so here because, under our assumptions, the sign of the functional equation is always $-1$ (cf. \cite[Chapters 5.3, 5.4]{Studnia-thesis}), so the claim would be vacuous. }

\demo{For newforms $f \in \mathcal{S}_2(\Gamma_1(N))$, $g \in \mathcal{S}_1(\Gamma_1(M))$, the \emph{imprimitive} Rankin--Selberg $L$-function $L(f,g,s)$ is defined as in \cite[\S 2.7]{KLZ15} as the meromorphic continuation of the series 
\[L^{NM}(\chi\theta,2s-1)\sum_{n \geq 1}{a_n(f)a_n(g)n^{-s}},\] where $\chi,\theta$ are the primitive characters attached to $f,g$ respectively, and the exponent $NM$ means that we are removing the Euler factors of the $L$-function at primes dividing $NM$. 
We claim that $L(f,g,1)$ is always defined and that it equals zero if and only if $\Lambda(f,g,s) = 0$ and that whether $L(f,g,1)$ vanishes does not change if we replace $(f,g)$ with $(\sigma(f) \otimes \beta,\sigma(g) \otimes \beta^{-1})$ for some $\sigma \in \Aut{\C}$ and some Dirichlet character $\beta$. This is a consequence of the following remarks:
\begin{itemize}[noitemsep,label=\tiny$\bullet$]
\item $\Lambda(f\otimes \beta,g \otimes \overline{\beta},s)=\Lambda(f,g,s)$ (because of the interpretation as Galois representations),
\item whether $L(f,g,1)$ vanishes does not change if $(f,g)$ is replaced with a Galois conjugate by \cite[Theorem 3]{Shimura76},
\item $\frac{L(f \otimes \beta,g \otimes \overline{\beta},s)}{L(f,g,s)}$ is a finite product of $F_p(p^{-s}) \in \C(p^{-s})$, where $p$ runs over primes dividing $MN$ or the conductor of $\beta$, and one checks explicitly that $F_p(p^{-1})$ is always defined and non-zero by the Ramanujan conjecture and \cite[Theorem 4.6.17]{Miyake},
\item when $(g,\overline{f})$ satisfies the conditions A), B), C) of \cite[\S 2]{Li-RS}, \cite[Theorem 2.2, Remark 2]{Li-RS} implies using the Ramanujan conjecture and \cite[Theorem 4.6.17]{Miyake} that $L(f,g,1)$ is non-zero if and only if $\Lambda(f,g,1) \neq 0$,
\item every pair of newforms has a simultaneous twist satisfying the conditions A), B), C) of \cite[\S 2]{Li-RS} by Theorem 2.1 of \emph{op.cit.}.
\end{itemize}

By this discussion, we may replace the special values of $\Lambda$ by the special values of the imprimitive $L$-series, and since $\det{\rho^{\ast}} = \omega_p^{-1}$, we are reduced to the case where $\omega=\omega'$. 

Let $N$ be an integer divisible by $30p$ and the primes of ramification of $\rho$. By Proposition \ref{quotient-by-ZN}, Proposition \ref{properties-tate-module-twisted-cartan-gq} and Proposition \ref{H2N-to-BSD}, we are done if we can show that, for $f \in \mathscr{S} \cup \mathscr{P} \cup [\mathscr{C} \backslash \mathscr{C}_M]$ and $g \in \mathcal{S}_1(\Gamma_1(D))$ attached to the irreducible odd Artin representation $\mrm{Ind}_C^N{\psi} \circ \rho^{\ast}$ for $\psi \in \mathscr{R}_2(f)$, if $L(f,g,1) \neq 0$, then the set of places $v$ of $F$ of residue characteristic $q \geq 7$ such that $\tilde{H}^2(\Z[1/Nq],T_{f,v} \otimes_{\OO_{F_v}} T_{g,v}; \Delta^{BK})$ is finite is infinite (where $T_{f,v}, T_{g,v}$ are lattices in $V_{f,v}, V_{g,v}$ respectively).

We want to apply \cite[Theorem 11.7.3]{KLZ15} to the couple $(f,g)$. To do so, we need to check that all the conditions below are simultaneously verified for infinitely many $v$ ($q$ denotes the residue characteristic and $\mfk{q} \subset \OO_F$ the valuation ideal for $v$):
\begin{enumerate}[noitemsep,label=(\alph*)]
\item \label{bsdtiu1} $f$ is ordinary at $v$,
\item \label{bsdtiu2} $T_{f,v}/\mfk{q}$ and $T_{g,v}/\mfk{q}$ are absolutely irreducible, 
\item \label{bsdtiu7} the semi-simplification of $(T_{f,v}/\mfk{q})_{|G_{\Q_q}}$ is the sum of two distinct characters,
\item \label{bsdtiu8} the restriction of $T_{g,v}/\mfk{q}$ to $G_{\Q_q}$ is the sum of two distinct characters,
\item \label{bsdtiu4} there exists $\sigma \in G_{\Q(\mu_{q^{\infty}})}$ acting on $T_{v} := T_{f,v} \otimes_{\OO_{F_v}} T_{g,v}$ by $-\mrm{id}$, 
\item \label{bsdtiu5} the restriction to $G_{\Q(\mu_{q^{\infty}})}$ of $T_v/\mfk{q}$ is absolutely irreducible,
\item \label{bsdtiu6} there exists $\sigma \in G_{\Q(\mu_{q^{\infty}})}$ such that $T_v/(\sigma-1)T_v$ is free over $\OO_{F_v}$ of rank one. 
\end{enumerate}

\ref{bsdtiu4} and \ref{bsdtiu5} hold for all but finitely many $v$ by \cite[Proposition 2.4]{Studnia-Euler}, therefore \ref{bsdtiu2} holds for all but finitely many $v$. \ref{bsdtiu7} always holds when $q \geq 7$ by \cite[Remark 7.2.7]{KLZ15}. 

Because the image of $G_{\Q}$ in $\mrm{Aut}(T_{g,v})$ is a finite group $M_g$, \ref{bsdtiu8} holds if the image of $\Fr_q$ in $M_g$ is not scalar. Let $K_g$ be the number field such that $G_{K_g}$ is the kernel of the projective Galois representation attached to $g$. Then \ref{bsdtiu8} holds for any $q$ which is not split in $K_g$. Note that $T_{g,v}$ is absolutely irreducible, so $M_g$ is not abelian. Hence $\mrm{Gal}(K_g/\Q)$ is not cyclic, so $[K_g:\Q] \geq 4$, and by Chebotarev, \ref{bsdtiu8} holds for all $v$ such that their residue characteristic lies outside a set of primes of density at most $1/4$. 

For every field homomorphism $\sigma: F \rar \C$ and every prime $q \nmid N$, one has $|\sigma(a_q(f))| \leq 2\sqrt{q}$. Therefore, if $q$ is a prime such that if $a_q(f) \neq 0$, then either $q \mid 6N\Delta_{F/\Q}$ or $a_q(f) \notin \sqrt{q\OO_F}$.   

Let $q$ be a prime such that $a_q(f) \notin \sqrt{q\OO_F}$: then there is a maximal ideal $\mfk{q}$ of $\OO_F$ with residue characteristic $q$ not containing $a_q(f)$, and $f$ is ordinary at $\mfk{q}$. By \cite[Th\'eor\`eme 15]{Serre-Chebotarev} and the prime number theorem, one has $a_q(f) \neq 0$ for a density one set of primes $q$. Therefore, for a density one set of primes $q$, there is a place $v$ of $F$ above $\mfk{q}$ at which \ref{bsdtiu1} holds. 

Finally, we show that \ref{bsdtiu6} holds for every $v$ but finitely many. If $f \in \mathscr{S}$, then by Lemma \ref{spcm-quadrichotomy} every one of its twists by a Dirichlet character of conductor $p$ is a newform of level $p^2$. Hence $f$ has no inner twist and the conclusion follows from the first condition in \cite[Theorem A]{Studnia-Euler}. If $f \in \mathscr{P}$, then by Lemma \ref{spcm-quadrichotomy} its only inner twist is by the inverse of its character, which is even: the conclusion follows by the second condition in \cite[Theorem A]{Studnia-Euler}. These arguments also apply if $f \in \mathscr{C} \backslash \mathscr{C}_M$ has no inner twist. 

Thus, we may assume that $f \in \mathscr{C} \backslash \mathscr{C}_M$ and that $f$ has an inner twist by some Dirichlet character $\chi$. By \cite[Lemma 1.3]{Studnia-Euler}, $\chi$ is quadratic of conductor dividing $p^2$, so it is the unique quadratic character of $\F_p^{\times}$, and the field $H$ attached to $f$ by Section 1 of \emph{op.cit.} is either $\Q$ or $\Q(\sqrt{p^{\ast}})$, with $p^{\ast}=(-1)^{\frac{p-1}{2}}p$. Since $\psi^p=\psi^{-1}$, the character of $g$ is exactly the quadratic character attached to the imaginary quadratic field $L$ such that $(\rho^{\ast})^{-1}(C)=G_L$. The extension $L/\Q$ is inert at $p$ by Proposition \ref{elliptic-curves-cartan-padic}, so it is linearly disjoint from $H$. The conclusion follows from \cite[Corollary 2.8 (v)]{Studnia-Euler}. 

}

While Proposition \ref{BSD-that-interests-us} suggests that we could prove the finiteness of the Mordell--Weil groups for many quotients of $J_{\Gamma}^{\alpha}(p)$, it turns out that, in most cases (determined in \cite[Theorem 10]{Studnia-thesis}), the functional equation sign is $-1$, so that the $L$-function has a forced vanishing at the center of the functional equation. However, Proposition \ref{sign-plus} gives us one case where it is possible that Proposition \ref{BSD-that-interests-us} is not vacuous. 

\cor[BSD-for-3C]{Let $\psi: C/3C \rar F^{\times}$ be an injective character and $g \in \mathcal{S}_1(\Gamma_1(M))$ be the newform attached to the irreducible odd Artin representation $\mrm{Ind}_C^N{\psi} \circ \rho^{\ast}$. Let $\omega \subset \mathscr{S}$ be stable under the Galois action and such that every $f \in \omega$ is Galois-conjugate to some $f'\in \omega$ such that $L(f',g,1) \neq 0$. Then the abelian variety $B_{\omega,3}$ described in Situation \ref{quotient-by-3C} has finite Mordell--Weil group.}

\section{Bounding the conductor of a rational point}
\label{sect:fimm-ccol}

\subsection{Exceptional primes and representations of $\GL{\F_p}$}
\label{subsect:exceptional-primes-steinberg}

\defi[exceptional-primes]{Let $p \geq 5$ be a prime congruent to $-1$ mod $3$. Fix some quadratic nonresidue $z \in \F_p^{\times}$, and let $Q_3$ denote the quotient $\F_{p^2}^{\times}/\F_{p^2}^{\times 3}$ (i.e. the quotient of $\F_{p^2}^{\times}$ by its subgroup of cubes). For each $t \in Q_3$, let 
\[u_{z,t} := \sum_{\substack{b \in \F_p\\1+b\sqrt{z} \in t}}{e^{\frac{2i\pi}{p}b}} \in \Z\left[e^{\frac{2i\pi}{p}}\right].\]

We say that a prime number $\ell$ is \emph{exceptional for $p$} if $\ell$ divides the norm of the ideal $I_{z,p} \subset \Z\left[e^{\frac{2i\pi}{p}}\right]$ generated by the $u_{z,t}$.

For $\sigma \in \mrm{Gal}(\Q\left(e^{\frac{2i\pi}{p}}\right)/\Q)$, one has $\sigma(I_{z,p})=I_{z\omega_p^{-2}(\sigma),p}$, so the notion of exceptional prime does not depend on the choice of $z$. 
}

\rem{For any prime $p$ as in Definition \ref{exceptional-primes}, there are at most finitely many exceptional primes, and we can check using MAGMA that when $p \leq 10000$, there are no exceptional primes. We are not able to prove that this claim holds generally. However, discussing the question with Hendrik Lenstra and Pieter Moree suggested the following heuristic: the ideal $I_{z,p}$ is generated by three elements whose sum is zero, so it is generated by any two elements among the $u_j$. A possible model is that these two $u_j$ to be random independent elements of $\Z[e^{\frac{2i\pi}{p}}]$ (in some informal sense), and $p$ has no exceptional elements if and only if these two elements are coprime. Now, it is well-known that the ``probability'' for two ``random'' independent elements of $\Z[\zeta_p]$ to be coprime is $\zeta_{\Q(\zeta_p)}(2)^{-1}$. 

Hence, suppose that we are given a collection $\mathcal{P}$ of primes, and choose ``random'' independent elements $\alpha_p,\beta_p \in \Z[\zeta_p]$ for each $p \in \mathcal{P}$. The Borel--Cantelli lemma suggests the following dichotomy: either $\prod_{p \in \mathcal{P}}{\zeta_{\Q(\zeta_p)}(2)^{-1}}$ is absolutely convergent, and then $\alpha_p$ and $\beta_p$ should almost surely be coprime for all but finitely many $p$, or $\prod_{p \in \mathcal{P}}{\zeta_{\Q(\zeta_p)}(2)^{-1}}$ diverges to zero, and then $\alpha_p$ and $\beta_p$ should almost surely have a common divisor for infinitely many $p$. 

Now, Fan, Languasco, Lunia and Moree were able to show \cite{Moree} that $\prod_{p}{\zeta_{\Q(\zeta_p)}(2)^{-1}}$ was indeed convergent, and their numerical data suggests that $\prod_{p \geq 10^4}{\zeta_{\Q(\zeta_p)}(2)^{-1}}$ is close to $1$. Hence, the Borel--Cantelli lemma suggests that exceptional primes exist for only finitely many $p$, and that it is ``likely'' that there are none. This is why we tentatively conjecture that exceptional primes do not exist.     
 }

Until the end of this section, we fix a prime $p \geq 5$ congruent to $-1$ modulo $3$. A brief overview of the well-known theory of representations of $\GL{\F_p}$, which we will use freely, is given in \cite[Chapter 6]{BH}.

\lem[choice-of-unip-cartan]{Let $\mathcal{U}$ be the collection of subgroups of order $p$ of $\GL{\F_p}$ and $\mathcal{C}$ be the collection of its non-split Cartan subgroups. Then the action by conjugation of $\GL{\F_p}$ on $\mathcal{U} \times \mathcal{C}$ is transitive.}

\demo{This is standard. The key fact is that if $B,C \leq \GL{\F_p}$ are respectively conjugate to the subgroup of upper-triangular matrices and a non-split Cartan subgroup, then one has $B \cap C=\F_p^{\times}\mrm{id}_2$, so $|B||C|=|B \cap C||\GL{\F_p}|$ and $\GL{\F_p}=BC=CB$.  }

\lem[steinberg-sum-model]{Let $\ell \nmid p(p^2-1)$ be a prime and $R=\Z_{(\ell)}$. Let $S$ be a non-zero finite free $R$-module endowed with a right action of $G := \GL{\F_p}$. Assume that for some field $F$ of characteristic zero, $S \otimes F$ embeds as a $F[G]$-module in $\mrm{Div}^0(\mathbb{P}^1(\F_p))^{\oplus d} \otimes F$ for some $d \geq 1$ (where the vectors of $\mathbb{P}^1(\F_p)$ are seen as equivalence classes of row vectors). Then $S \simeq S_0^{\oplus d}$ as $R[G]$-modules.}

\demo{Note that $S_0:= \mrm{Div}^0(\mathbb{P}^1(\F_p)) \otimes_{\Z} R$ is a finite free $R$-module and that $|G|$ is invertible in $R$. By \cite[Lemma 20]{SerreLinReps}, $S_0$ is a projective $R[G]$-module. Therefore, it is enough by induction to show that there exists a non-trivial surjective map $S \rar S_0$ of $R[G]$-modules. 

Let $\phi: S \rar S_0$ be a non-zero map. Since $S_0 \otimes \Q$ is an irreducible $\Q[G]$-module, $\phi \otimes \Q$ is surjective, so $\phi$ has finite cokernel. Suppose that $\varphi$ is not surjective, then by Nakayama $\phi \otimes \F_{\ell}: S \otimes \F_{\ell} \rar S_0 \otimes\F_{\ell}$ is not surjective. Since $S_0 \otimes \F_{\ell}$ is an irreducible $\F_{\ell}[G]$-module \cite[Proposition 43]{SerreLinReps}, $\phi \otimes \F_{\ell}$ is zero, hence $\varphi/\ell: S \rar S_0$ is a well-defined non-trivial map. Hence, if there exists a non-trivial map $\phi: S \rar S_0$ of $R[G]$-modules, the map $\phi$ with the smallest possible cokernel is surjective. 

To sum up, we are done if there exists a non-trivial map $S \rar S_0$ of $R[G]$-modules; it is enough to show that a non-trivial map $S \otimes \Q \rar S_0 \otimes\Q$ of $\Q[G]$-modules exists. Now, one has an injection of $\Q[G]$-modules $S \otimes \Q \subset S \otimes F \simeq S_0^{\oplus d} \otimes F \simeq S_0^{\oplus d[F:\Q]}$, whence the conclusion. 
}

\lem[explicit-calculation-steinberg]{Let $\ell \nmid p(p^2-1)$ be a prime, $R=\Z_{(\ell)}$, $G := \GL{\F_p}$ and $S_0$ be the right $R[G]$-module $\mrm{Div}^0(\mathbb{P}^1(\F_p))$ as in Lemma \ref{steinberg-sum-model}. Fix a quadratic non-residue $\varepsilon \in \F_p^{\times}$, let $U=\begin{pmatrix}1 & 1\\0 & 1\end{pmatrix}$ and let \[C=\{\begin{pmatrix}a & b\varepsilon\\b & a\end{pmatrix},\, (a,b) \in \F_p^2 \backslash \{(0,0)\}\}.\]
Let $\mathbf{1},j,j' \subset \F_{p^2}^{\times}$ denote the three cosets modulo its cubes (where $\mathbf{1}$ is the trivial coset), and, for $w \in \{\mathbf{1},j,j'\}$, let 
\[\gamma_w = \sum_{\substack{b \in \F_p\\1+b\sqrt{\varepsilon}} \in w}{U^{b\varepsilon}} \in \Z[U].\]  
Then the right $R[U]$-module generated by $S_0^{3C}$ is $S_0 \cdot \left(\gamma_j-\gamma_1,\gamma_{j'}-\gamma_1\right)$. In particular, if $\ell$ is not exceptional for $p$, $S_0^{3C}$ generates $S_0$ as a $R[U]$-module. 
}

\demo{For $w \in \{\mathbf{1},j,j'\}$, let $I_w$ be the set of $b \in \F_p$ such that $1+b\sqrt{\varepsilon} \in w$: the three $I_w$ form a partition of $\F_p$. Note that $|I_j|=|I_{j'}|=\frac{p+1}{3}=1+|I_{\mathbf{1}}|$. 

Let $t=[1:0]-[0:1] \in S_0$: then $t \cdot U^i = [1:i]-[0:1]$ for every $i \in \F_p$, so that $z \in R[U] \mapsto t \cdot z \in S_0$ is an isomorphism of right $R[U]$-modules. By a character computation, $(S_0 \otimes \Q)^{3C}$ has dimension two, so $S_0^{3C}$ is a free $R$-module of rank $2$ which is saturated in $S_0$.

For $w \in \{j,j'\}$, let 
\[z_w = \sum_{b \in I_w}{[1:\varepsilon b]}-\sum_{b \in I_{\mathbf{1}}}{[1: \varepsilon b]}-[0:1] \in S_0.\] 
Because $\begin{pmatrix} a & \varepsilon b\\\ast & \ast\end{pmatrix} \in C\mapsto a+b\sqrt{\varepsilon} \in \F_{p^2}^{\times}$ is a group isomorphism, $z_w$ is fixed by $3C$. Moreover, $(a,b) \in R^{\oplus 2} \mapsto az_j+bz_{j'} \in S_0$ is clearly injective with $R$-saturated image, so that $Rz_j+Rz_{j'} = S_0^{3C}$.
Now, one computes directly that $z_w=t \cdot (\gamma_w-\gamma_{\mathbf{1}})$ for $w \in \{j,j'\}$. 

Finally, we show that if $\ell$ is not exceptional for $p$, then $(\gamma_j-\gamma_1,\gamma_{j'}-\gamma_1)R[U]=R[U]$. Indeed, one has an isomorphism $R[U] \rar R[e^{2i\pi/p}] \times R$, where the first projection is $\pi_1: U \mapsto e^{2i\pi/p}$ and the second projection is $\pi_2: U \mapsto 1$. One has $\pi_2(\gamma_j)=\pi_2(\gamma_{j'})=\pi_2(\gamma_{\mathbf{1}})+1=\frac{p+1}{3}$, so it is enough to show that $\pi_1(\gamma_j-\gamma_{\mathbf{1}}), \pi_1(\gamma_{j'}-\gamma_{\mathbf{1}}) \in R[e^{2i\pi/p}]$ are coprime. Now, for $w \in \{\mathbf{1},j,j'\}$, $\pi_1(\gamma_w)$ is the element $u_{\varepsilon^{-1},w}$ in the notation of Definition \ref{exceptional-primes}. Hence, we need to show that $I_{\varepsilon^{-1},p}R[e^{2i\pi/p}]=(u_{\varepsilon^{-1},j}-u_{\varepsilon^{-1},\mathbf{1}},u_{\varepsilon^{-1},j'}-u_{\varepsilon^{-1},\mathbf{1}})R[e^{2i\pi/p}]$. It is enough to show the forward inclusion: since $\ell \neq 3$, one has 
\begin{align*}
u_{\varepsilon^{-1},\mathbf{1}} &= -\frac{1}{3}\sum_{w \in \{j,j'\}}{(u_{\varepsilon^{-1},w}-u_{\varepsilon^{-1},\mathbf{1}})},\\
u_{\varepsilon^{-1},j} &= u_{\varepsilon^{-1},\mathbf{1}} + (u_{\varepsilon^{-1},j}-u_{\varepsilon^{-1},\mathbf{1}})
\end{align*}
and the conclusion follows. 
}

\prop[exceptional-and-reptheo]{Let $\ell \nmid p(p^2-1)$ be a prime, $R=\Z_{(\ell)}$, and $R_p = R \otimes_{\Z} \Z[\zeta_p]$. Let $S$ be a non-zero finite free $R$-module endowed with a right action of $G=\GL{\F_p}$. Let $N \leq G$ be the normalizer of a non-split Cartan subgroup $C$ and $U \in G$ be a unipotent. Let $a_1: S \rar R_p$ be a $R$-linear map such that for every $x \in S$, one has $a_1(x \mid U) = \zeta_p^{-1} a_1(x)$. Assume furthermore that for some field $F$ of characteristic zero, $S \otimes F$ is isomorphic to $\mrm{Div}^0(\mathbb{P}^1(\F_p))^{\oplus d} \otimes F$ for some $d \geq 1$. 
If $a_1(S)$ generates the unit ideal of $R_p$ and $\ell$ is not exceptional for $p$, then $a_1(S^{3C})$ generates the unit ideal of $R_p$.}

\demo{After conjugating by Lemma \ref{choice-of-unip-cartan} and applying the Galois action, we may assume that $U=\begin{pmatrix}1 & 1\\0 & 1\end{pmatrix}$ and $C=\{\begin{pmatrix}a & b\varepsilon\\b & a\end{pmatrix},\, (a,b) \in \F_p^2 \backslash \{(0,0)\}\}$, for some $\varepsilon \in \F_p^{\times}$ which is not a square. We may also assume by Lemma \ref{steinberg-sum-model} that $S = S_0^{\oplus d}$, where $S_0$ is as in \emph{loc.cit.}. 

Because $a_1(x \mid U)=\zeta_p^{-1}a_1(x)$ for every $x \in S$, $a_1(S)$ is an ideal of $R[\zeta_p]$, and the ideal of $R[\zeta_p]$ generated by $a_1(S^{3C})$ is exactly the image under $a_1$ of the $R[U]$-submodule of $S$ generated by $S^{3C}$. By Lemma \ref{explicit-calculation-steinberg}, since $\ell$ is not exceptional for $p$, the right $R[U]$-module generated by $S^{3C}$ is $S$: hence the ideal generated by $a_1(S^{3C})$ is exactly $a_1(S)$, and the conclusion follows.   }

\rem{One can show in a very similar way that if $a_1(S^{3C})$ generates the unit ideal of $R[\zeta_p]$, then $\ell$ is not exceptional for $p$.  }

\subsection{The Chabauty--Coleman argument}
\label{subsect:chabauty-coleman}

\lem[t1-formal-immersion]{Let $f \in \mathcal{N}$ with character $\chi$. Let $\mu_f: \mathbb{T} \rar \C$ be the character sending $T_n$ for $n \geq 1$ prime to $p$ (resp. $\diam{n}$ for $n \in \F_p^{\times}$) to $a_n(f)$ (resp. $\chi(n)$). Let $c \in X_{G_0}(p)(\Z[1/p,\zeta_p])$ be a section of the cuspidal subscheme, and $a_1: H^0(X_{G_0}(p),\Omega^1_{X(p)/\Z[1/p]}) \rar \Z[1/p,\zeta_p]$ be the first coefficient in the $q$-expansion at $c$ attached by Corollary \ref{properties-of-q-exp}. Then the restriction 
\[a_1: H^0(X_{G_0}(p),\Omega^1_{X_{G_0}(p)/\Z[1/p]})[\ker{\mu_f}] \rar \Z[1/p,\zeta_p]\] is surjective.}

\demo{Note that for some non-trivial unipotent $U \in \SL{\F_p}$, one has $a_1(\omega \mid U) = \zeta_p a_1(\omega)$. Since $H^0(X_{G_0}(p),\Omega^1_{X_{G_0}(p)/\Z[1/p]})[\ker{\mu_f}]$ is stable under $\GL{\F_p}$, the image of 
\[a_1: H^0(X_{G_0}(p),\Omega^1_{X_{G_0}(p)/\Z[1/p]})[\ker{\mu_f}] \rar \Z[1/p,\zeta_p]\]
is therefore an ideal of $\Z[1/p,\zeta_p]$. Hence, it is enough to find $\omega \in H^0(X(p),\Omega^1_{X_{G_0}(p)/\Z[1/p]})[\ker{\mu_f}]$ such that $a_1(\omega) \in \Z[1/p,\zeta_p]^{\times}$. By Corollary \ref{properties-of-q-exp}, it is enough to find a faithfully flat $\Z[1/p]$-algebra $R \subset \C$ and $\omega \in \mathcal{S}_2(\Gamma(p)) \otimes \Z[\zeta_p]$ such that
\begin{itemize}[noitemsep,label=\tiny$\bullet$]
\item The $q$-expansion of $\omega$ is contained in $R[[q^{1/p}]] \otimes \Z[\zeta_p]$ with invertible first coefficient,
\item for every prime $\ell \neq p$, one has $\omega \mid \left(\Gamma(p)\begin{pmatrix}\ell & 0\\0 & 1\end{pmatrix}\Gamma(p) \otimes \underline{\ell}\right)= a_{\ell}(f)\omega$,
\item for every $n \in \F_p^{\times}$, one has $\omega \mid \left(\begin{pmatrix}n & 0\\0 & n^{-1}\end{pmatrix} \otimes \underline{n^2}\right)= \chi(f)\omega$.
\end{itemize}

We take $R=\OO_F[1/p]$, where $F \subset \C$ is a Galois number field over $\Q$ and contains all the Fourier coefficients of $f$, and it is an elementary verification that $\omega = \sum_{n \geq 1}{\overline{a_n(f)}q^{n/p}} \otimes 1$ satisfies all the requirements. 
}

\prop[nonexc-implies-fi]{Let $\mfk{f} \subset \mathscr{S}$ be a collection of newforms for $\mathcal{S}_2(\Gamma_0(p))$ stable under the action of Galois. Fix a non-split Cartan subgroup $C \subset \GL{\F_p}$ and let $N$ be its normalizer, and let $I_{\mfk{f},3}$ be the ideal appearing in Situation \ref{quotient-by-3C}.  

Let $c \in X_{G_0}(p)(\Z[1/p,\zeta_p])$ be a section of the cuspidal subscheme and $\ell \nmid p(p^2-1)$ be a prime. If $\ell$ is not exceptional for $p$, then 
the natural map 
\[\nu: H^0(X_{G_0}(p)_{\Z_{(\ell)}},\Omega^1_{X_{G_0}(p)_{\Z_{(\ell)}}/\Z_{(\ell)}})[I_{\mfk{f},3}] \rar H^0(\Sp{\Z_{(\ell)}[\zeta_p]},c^{\ast}\Omega^1_{X_{G_0}(p)_{\Z_{(\ell)}}/\Z_{(\ell)}})\] is surjective.  
}

\demo{By construction, $H^0(X_{G_0}(p)_{\Z_{(\ell)}},\Omega^1_{X_{G_0}(p)_{\Z_{(\ell)}}/\Z_{(\ell)}})[I_{\mfk{f},3}]$ is a $\OO(X_{G_0}(p))$-module, so that $\nu$ is $\Z_{(\ell)}[\zeta_p]$-linear. 

Fix $f \in \mfk{f}$. The domain of $\nu$ clearly contains $H^0(X_{G_0}(p)_{\Z_{(\ell)}},\Omega^1_{X_{G_0}(p)_{\Z_{(\ell)}}/\Z_{(\ell)}})[\ker{\mu_f}]^{3C}$, so it is enough to show that the image of 
\[a_1: H^0(X_{G_0}(p)_{\Z_{(\ell)}},\Omega^1_{X_{G_0}(p)_{\Z_{(\ell)}}/\Z_{(\ell)}})[\ker{\mu_f}]^{3C} \rar \Z_{(\ell)} \otimes \Z[\zeta_p]\] generates the unit ideal.

By Lemma \ref{t1-formal-immersion}, $a_1: H^0(X_{G_0}(p)_{\Z_{(\ell)}},\Omega^1_{X_{G_0}(p)_{\Z_{(\ell)}}/\Z_{(\ell)}})[\ker{\mu_f}] \rar \Z_{(\ell)} \otimes \Z[\zeta_p]$ is surjective. By the properties of $q$-expansion (Corollary \ref{properties-of-q-exp}), the claim follows from Proposition \ref{exceptional-and-reptheo} if we can show that for some field $F$ of characteristic zero, $H^0(X_{G_0}(p)_{\Z_{(\ell)}},\Omega^1_{X_{G_0}(p)_{\Z_{(\ell)}}/\Z_{(\ell)}})[\ker{\mu_f}] \otimes F$ (which is stable under $\GL{\F_p}$) is a sum of copies of the (right) Steinberg representation $\mrm{St}$ of $\GL{\F_p}$. By Proposition \ref{h10-xp}, it is enough to show that if $f' \in \mathcal{N}$ is such that $\ker{\mu_f}\subset \ker{\mu_{f'}}$, then $f' \in \mathscr{S}$. Since both ideals are saturated in $\mathbb{T}$ and correspond to maximal ideals of $\mathbb{T} \otimes \Q$, one has $\ker{\mu_f}=\ker{\mu_{f'}}$. Then there is a morphism $\gamma: \Q(f) \rar \Q(f')$ sending $a_n(f) = \mu_f(T_n)$ to $\mu_{f'}(T_n)=a_n(f')$ for every $n \geq 1$ prime to $p$. Then $\gamma$ extends to an automorphism of $\C$, so $f,f'$ are Galois-conjugates and $f' \in \mathscr{S}$, whence the conclusion.}

\cor[nonexc-implies-fi-relative]{Let $\ell \nmid p(p^2-1)$ be a prime. Let $\Gamma$ be a $p$-torsion group over the ring of integers $\OO$ of some finite extension $K/\Q_{\ell}$ and let $\alpha: \Gamma \times \Gamma \rar (\mu_p)_{\OO}$ be a Weil pairing. Let $P,Q \in \Gamma(\overline{K})$ form a basis and assume that the image of the morphism $\rho: G_{K} \rar \GL{\F_p}$ giving the Galois action on $\Gamma(\overline{K})$ in the basis $(P,Q)$ lands in the transpose of the normalizer of some non-split Cartan subgroup $C$. 
Let $\mfk{f}$ and $I_{\mfk{f},3}$ be as in Corollary \ref{nonexc-implies-fi}. Let $\iota: Z \rar X_{\Gamma}(p)$ be the inclusion corresponding to the cuspidal subscheme. If $\ell$ is not exceptional for $p$, then the canonical map
\[H^0(X_{\Gamma}(p),\Omega^1_{X_{\Gamma}(p)/\OO})[I_{\mfk{f},3}] \otimes_{\OO(X_{\Gamma}(p))} \OO(Z) \rar H^0(Z,\iota^{\ast}\Omega^1_{X_{\Gamma}(p)/\OO})\]
is surjective.
}

\demo{The formation of this map, and whether it is surjective, commutes with faithfully flat base change. Therefore we may replace $(\Gamma,\OO)$ with $(G_0,\Z_{(\ell)})$, and this case follows from Proposition \ref{nonexc-implies-fi} because $Z$ is a disjoint reunion of copies of $\Sp{\Z_{(\ell)}[\zeta_p]}$.  }

\cor[nonexc-implies-fi-twisted]{Let $\ell \nmid p(p^2-1)$ be a prime. Let $\Gamma/\OO,P,Q,C,\mfk{f},I_{\mfk{f},3}$ be as in Corollary \ref{nonexc-implies-fi-twisted}, and let $\alpha$ be a Weil pairing for $\Gamma$. 

Let $\mfk{f} \subset \mathscr{S}$ be a collection of newforms for $\mathcal{S}_2(\Gamma_0(p))$ stable under the action of Galois and, for every $n \in \F_p^{\times}$, $B_{\mfk{f},n}/\Z_{\ell}$ be the quotient of $J_{\Gamma}^{\alpha^n}(p)$ constructed in Situation \ref{quotient-by-3C} (i.e. on which $I_{\mfk{f},3}$ acts trivially). 

Let $c \in X_{\Gamma}^{\alpha^n}(p)(\OO)$ be a section of the cuspidal subscheme, $\varphi: (B_{\mfk{f},n})_{K} \rar (B_{\mfk{f},n})_K$ an isogeny, and $F: X_{\Gamma}^{\alpha^n}(p)_K \rar (B_{f,n})_K$ sending $x$ to $\varphi([x-c])$. 

If $\ell$ is not exceptional for $p$, then there exists \[\omega \in H^0(X_{\Gamma}^{\alpha^n}(p),\Omega^1_{X_{\Gamma}^{\alpha^n}(p)/\OO}) \cap F^{\ast} H^0((B_{\mfk{f},n})_K,\Omega^1_{(B_{\mfk{f},n})_K/K})\] such that $c^{\ast}\omega$ generates the free $\OO$-module $H^0(\Sp{\OO},c^{\ast}\Omega^1_{X_{\Gamma}^{\alpha^n}(p)/\Sp{\OO}})$ of rank one.  
}

\demo{If $\varphi$ is an isogeny, then pulling back by $\varphi$ is an automorphism of $H^0(B_{\mfk{f},n},\Omega^1_{B_{\mfk{f},n}/K})$: as a consequence, $H^0(X_{\Gamma}^{\alpha^n}(p),\Omega^1_{X_{\Gamma}^{\alpha^n}(p)/\OO}) \cap F^{\ast} H^0((B_{\mfk{f},n})_K,\Omega^1_{(B_{\mfk{f},n})_K/K})$ does not depend on what $\varphi$ is, and it is precisely the submodule of $H^0(X_{\Gamma}^{\alpha^n}(p),\Omega^1_{X_{\Gamma}^{\alpha^n}(p)/\OO})$ killed by $I_{\mfk{f},3}$ by Proposition \ref{quotient-by-3C}. The result is then a consequence of Corollary \ref{nonexc-implies-fi-twisted}.}

\prop[coleman-integration-rules]{Let $\ell > 2$ be a prime and let $X$ be a smooth projective $\Z_{\ell}$-scheme of relative dimension one. Let $A$ be an abelian scheme over $\Z_{\ell}$, and let $f: X \rar A$ be a $\Z_{\ell}$-morphism. Let $c \in X(\F_{\ell})$ be a point and let $P,Q \in X(\Z_{\ell})$ be two points whose reduction modulo $\ell$ is $c$. Assume that there exists $\omega \in H^0(A,\Omega^1_{A/\Q_{\ell}})$ such that $\omega' := f^{\ast}\omega$ lies in $H^0(X,\Omega^1_{X/\Z_{\ell}})$ and $\omega'_c$ generates $\Omega^1_{X/\Z_{\ell},c}$. Then, either $P=Q$ or $f(P)-f(Q)$ has infinite order. 
}

\demo{We use Coleman's integration theory as defined in \cite{ColmInt}; this is essentially a retelling of a special case of the main argument of \cite{ColChab}. Note that we can multiply $f$ by any non-zero integer: it is enough to show that if $f(P)=f(Q)$, then $P=Q$.

Because $X$ is smooth over $\Z_{\ell}$, the ring $\OO_{X,c}$ (resp. $\OO_{X,c}/(\ell)$) is regular local of dimension $2$ (resp. $1$), so a system of parameters is given by $(\ell,t)$ for some $t \in \OO_{X,c}$ such that $\OO_{X,c}/(\ell)$ is a discrete valuation ring with uniformizer $t$, and in particular the obvious map $\Z_{\ell}[[t]] \rar \widehat{\OO_{X,c}}$ is an isomorphism. Therefore, we can write $\omega'_c = \sum_{n \geq 0}{a_nt^n}\mrm{d}t$ in the space of continuous differentials of $\widehat{\OO_{X,c}}$ over $\Z_{\ell}$ and the assumption on $\omega'_c$ implies that $a_0 \in \Z_{\ell}^{\times}$. 

Since $\Z_{\ell}[[t]] \rar \widehat{\OO_{X,c}}$ is an isomorphism, $z \mapsto t(z)$ defines a bijection from the set of $\Z_{\ell}$-points of $X$ that reduce to $c$ modulo $\ell$ into $\ell\Z_{\ell}$. 

Because $f(P)=f(Q)$, the properties of Coleman integration then imply that 
\begin{align*}
0&=\int_{f(P)}^{f(Q)}{\omega}=\int_P^Q{\omega'} = \sum_{n \geq 0}{a_n\frac{t(P)^{n+1}-t(Q)^{n+1}}{n+1}}\\
&=(t(P)-t(Q))\left(a_0+\sum_{n \geq 1}{\frac{a_n}{n+1}\sum_{k=0}^n{t(P)^kt(Q)^{n-k}}}\right).
\end{align*}

Since $\ell$ is odd, one has $n-v_{\ell}(n+1) \geq 1$ for every $n \geq 1$, hence $\frac{a_n}{n+1}\sum_{k=0}^n{t(P)^kt(Q)^{n-k}} \in \ell\Z_{\ell}$. Since $a_0 \in \Z_{\ell}^{\times}$, this implies that $t(P)=t(Q)$, thus $P=Q$. }

\cor[nicest-chabauty-coleman]{Let $f: X \rar A$ be a morphism from a smooth projective curve to an abelian variety over $\Q$. Let $\ell$ be an odd prime at which $X$ and $A$ have good reduction, and let $\mathcal{X}$ be a smooth proper model of $X$ over $\Z_{\ell}$. Let $c \in \mathcal{X}(\Z_{\ell})$ and $c_0 \in \mathcal{X}$ be the image by $c$ of the maximal point of $\Sp{\Z_{\ell}}$. Assume that 
\begin{itemize}[noitemsep,label=\tiny$\bullet$]
\item $f(c_{\Q_{\ell}})=0$,
\item $A(\Q)$ is finite,
\item there exists $\omega \in H^0(A,\Omega^1_{A/\Q_{\ell}})$ such that $\omega' := f^{\ast}\omega$ lies in $H^0(\mathcal{X},\Omega^1_{\mathcal{X}/\Z_{\ell}})$ and $\omega'_{c_0}$ generates $\Omega^1_{\mathcal{X}/\Z_{\ell},c_0}$.
\end{itemize}
Let $P \in X(\Q)$ be a point whose reduction modulo $\ell$ (in the model $\mathcal{X}$) is $c_{\F_{\ell}}$. Then $P_{\Q_{\ell}}=c_{\Q_{\ell}}$.}  

\demo{The N\'eron model $\mathcal{A}$ of $A_{\Q_{\ell}}$ over $\Z_{\ell}$ is an abelian scheme by \cite[1.4/3]{BLR} and $f$ extends to $f_1: \mathcal{X} \rar \mathcal{A}$. The group $A(\Q)$ is finite, so $f(P) \in A(\Q)$ is a torsion point; by the valuative criterion of properness, $P$ lifts to $P_1 \in \mathcal{X}(\Z_{\ell})$. The restriction to $\Q_{\ell}$ of $f_1(P_1)-f_1(c) \in \mathcal{A}_f(\Z_{\ell})$ is the torsion point $f(P)$, so, by the N\'eron model property, $f_1(P_1)-f_1(c)$ is a torsion section of $\mathcal{A}_f$. Moreover, the $\Z_{\ell}$-points $P_1$ and $c$ of $\mathcal{X}$ have the same reduction modulo $\ell$ by assumption. Finally, our last assumption means that we can apply Proposition \ref{coleman-integration-rules}, and thus $P_1=c$, whence $P_{\Q_{\ell}}=c_{\Q_{\ell}}$. }

\prop[fimm-main]{Let $p\geq 11$ be a prime and $E/\Q$ be an elliptic curve such that the image of the attached Galois representation $\rho: G_{\Q} \rar \GL{\F_p}$ is contained in the normalizer of a non-split Cartan subgroup $C$. Let $K/\Q$ be the imaginary quadratic extension defined by $\rho^{-1}(C)$, $\psi: C \rar \C^{\times}$ be a character of $C$ with order exactly $3$ and $\rho_3 = \mrm{Ind}_{K}^{\Q}{\left[\psi\circ\rho_{|G_K}\right]}$ be an Artin representation. Assume furthermore that:
\begin{itemize}[noitemsep,label=\tiny$\bullet$]
\item the image under $\rho$ of an inertia group at $p$ is properly contained in $C$,
\item the image of $\rho$ is equal to $N$, and therefore $\rho_3$ is attached to a newform $g \in \mathcal{S}_1(\Gamma_1(M))$,   
\item there exists $f \in \mathcal{S}_2(\Gamma_0(p))$ such that the central Rankin--Selberg $L$-value $L(f,g,1)$ is nonzero. 
\end{itemize}
Let $\ell$ be a prime at which $E$ does not have potentially good reduction. Then $\ell \equiv \varepsilon(\ell)\mod{p}$ and $\ell$ is exceptional for $p$, where $\varepsilon$ is the quadratic character attached to $K/\Q$.  }

\rem{It follows from the functional equation sign computation (e.g. \cite[Theorem 10]{Studnia-thesis} that the third assumption actually implies the first one, because otherwise the sign of the functional equation for $L(f,g,s)$ is $-1$. }

\demo{By Lemma \ref{elliptic-curves-cartan-complex} and Proposition \ref{elliptic-curves-cartan-padic}, $K$ is indeed an imaginary quadratic field and inert at $p$, $\rho$ maps any inertia subgroup of $G_{\Q}$ at $p$ to the unique subgroup of $C$ of index $3$, one has $p \equiv 2,5 \pmod{9}$, and $E$ has potentially good reduction at $p$ (so $p \neq \ell$). After replacing $E$ with a quadratic twist (which does not change $\rho_3$), we may assume by \cite[Lemma V.5.2, Theorem V.5.3]{AEC2} that $E$ has multiplicative reduction at $\ell$. By Proposition \ref{elliptic-curves-cartan-nram}, $\rho$ is then unramified at $\ell$ and one has $\ell \equiv \varepsilon(\ell)\pmod{p}$. The group scheme $E[p]_{\Q}$ with its Weil pairing extends to a $p$-torsion group $\Gamma$ over $\Z_{(\ell)}$ endowed with a Weil pairing $\alpha$, and $E$ defines a rational point $P$ of $X_{\Gamma}^{\alpha}(p)$, whose reduction modulo $\ell$ lies in the cuspidal subscheme. The cuspidal subscheme is finite \'etale over $\Z_{(\ell)}$, so it has a $\Z_{\ell}$-point $c$ such that $c \equiv P \mod{\ell}$. 

Let $q \equiv 1 \mod{p}$ be a prime distinct from $\ell$ and $\mfk{f} \subset \mathscr{S}$ be the $G_{\Q}$-orbit of the $f$ given by the assumptions. Let $B_{\mfk{f},\alpha}$ be the quotient of $J_{\Gamma}^{\alpha}(p)_{\Q}$ on which $I_{\mfk{f},3}$ acts trivially (considered in Situation \ref{quotient-by-3C}). Since $T_q-q-1$ kills the cuspidal subscheme of $X_{\Gamma}^{\alpha}(p)_{\Q}$, the map $\pi: x \in X_{\Gamma}^{\alpha}(p)_{\Q} \mapsto (T_q-q-1)(x-c) \in B_{\mfk{f},\alpha}$ is well-defined. 

Since $\rho_3$ is isomorphic to $\left[\mrm{Ind}_{C^T}^{N^T}{\psi}\right]\circ \rho^{\ast}$, Corollary \ref{BSD-that-interests-us} implies that $B_{\mfk{f},\alpha}(\Q)$ is finite. Moreover, $B_{\mfk{f},\alpha}$ is a quotient of $J_{\Gamma}^{\alpha}(p)$ (which has good reduction at $\ell$), so $B_{\mfk{f},\alpha}$ has good reduction at $\ell$ by e.g. \cite[Proposition 1.2]{FreyMazur}. Finally, $T_q-q-1$ is an isogeny of $B_{\mfk{f},\alpha}$, so by Corollary \ref{nonexc-implies-fi-twisted}, the assumptions of Corollary \ref{nicest-chabauty-coleman} are satisfied. Therefore $P_{\Q_{\ell}}=c_{\Q_{\ell}}$, which is impossible since $j(P) \in \Q$ and $j(c) =\infty$.} 

\cor[fimm-invariants]{Under the assumptions of Proposition \ref{fimm-main}, the conductor of $E$ is equal to $p^2N_{\rho}L$, where $N_{\rho}$ is the Artin conductor of $\rho$, and $L$ is a product of pairwise distinct primes $\ell$ satisfying the following properties:
\begin{itemize}[noitemsep,label=\tiny$\bullet$]
\item $\ell \equiv \varepsilon(\ell) \mod{p}$ and $\ell$ is exceptional for $p$,
\item $\ell$ is prime to $pN_{\rho}$,
\item $j(E)$ is a $p$-th power in $\Q_{\ell}$.
\end{itemize}
Moreover, the denominator of $j(E)$ is the $p$-th power of some integer $L_1$ divisible by $L$.  
}

\demo{If $\ell \neq p$ is a prime of potentially good reduction for $E$, then the conductor exponent of $E$ at $\ell$ is the same as the Artin conductor exponent for the representation $E[p](\Qbar)$ by Proposition \ref{elliptic-curves-cartan-potgood}. By Proposition \ref{elliptic-curves-cartan-padic}, $E$ has bad and potentially good reduction at $p$. 

Let us now prove that if $\ell$ is a prime of additive and potentially multiplicative reduction for $E$, then $E$ and $\rho$ have the same conductor exponents at $\ell$. By Proposition \ref{fimm-main}, $\ell \equiv \pm 1 \mod{p}$, so $\ell > 3$, so it is enough to show that $\rho$ has no invariants under the inertia at $\ell$. By Proposition \ref{elliptic-curves-cartan-nram}, the image by $\rho$ of the inertia at $\ell$ is contained in $\pm \mrm{id}_2$ and non-trivial, so this image is exactly $\{\pm \mrm{id}_2\}$, which implies the conclusion. 

The claim for primes of multiplicative reduction follows similarly from Proposition \ref{elliptic-curves-cartan-nram} and Proposition \ref{fimm-main}.}

\section{Applying the result to an elliptic curve}
\label{sect:effectivity-results}

The goal of this section is to give explicit criteria that tell us whether Proposition \ref{fimm-main} applies to a given elliptic curve $E$. It is divided in three parts. The first one is dedicated to the determination of the newform $g$ appearing in this Proposition and provides criteria in order to determine whether the Galois image of inertia at $p$ is smaller than expected, as required in \emph{loc.cit.}. The second part uses an explicit Waldspurger formula by Cai--Shu--Tian \cite{CaiShuTian} to give a ``finite'' criterion letting one check whether the final assumption of Proposition \ref{fimm-main} is satisfied. The third part applies these results to certain elliptic curves.

\subsection{The representation $\rho'$ and the small inertia image condition}
\label{subsect:explicit-rhoprime}

In this section, $E/\Q$ denotes an elliptic curve and $p \geq 11$ is a prime congruent to $2$ or $5$ modulo $9$. Let $P,Q \in E[p](\Qbar)$ be a basis of $E[p](\Qbar)$, and let $\rho: G_{\Q} \rar \GL{\F_p}$ be the matrix of the Galois action in the basis $(P,Q)$. We assume that the image of $\rho$ is contained in the normalizer $N$ of some non-split Cartan subgroup $C$. Let $K/\Q$ be the number field such that $\rho^{-1}(C) =G_K$, it is imaginary quadratic and inert at $p$ by Lemma \ref{elliptic-curves-cartan-complex} and Proposition \ref{elliptic-curves-cartan-padic}. 

We know by Proposition \ref{elliptic-curves-cartan-padic} and Lemma \ref{elliptic-curves-cartan-complex} that $3C \subset \rho(G_{\Q}) \not\subset C$. Moreover, because $p \equiv -1 \mod{3}$ and $p \not\equiv -1\mod{9}$, $C$ is the direct sum of $3C$ and its $3$-torsion group, the former having cardinality prime to $3$ and the latter being cyclic of order $3$. Since $p \equiv -1\mod{3}$, $3C$ also contains the scalar matrices and has surjective determinant.

We let $D$ denote the absolute value of the discriminant of $K$ and $\Delta$ the discriminant of a minimal model of $E$. Let $\psi: C \rar \C^{\times}$ be a character of order $3$ and $\rho' := \left(\mrm{Ind}_C^N{\psi}\right)\circ \rho^{\ast}$: it is a representation $G_{\Q} \rar \GL{\C}$. 

\rem{
\begin{itemize}[noitemsep,label=\tiny$\bullet$]
\item $\rho'$ is also isomorphic to $\left(\mrm{Ind}_C^N{\psi}\right)\circ \rho$,
\item Since $N/C \simeq \mfk{S}_3$, $\mrm{Ind}_{C/3C}^{N/3C}{\psi}$ is irreducible and can be realized over $\Q$,  
\item $\det{\rho'}$ is the quadratic character attached to $K$ --- so that $\rho'$ is odd,
\item the composition $G_{\Q}^{ab} \overset{\mrm{Ver}}{\longrightarrow} G_K^{ab} \overset{\psi \circ \rho^{\ast}_{|G_K}}{\longrightarrow} \C^{\times}$ is trivial.
\end{itemize}
}

We are thus interested to find out for which primes $\ell$ the group $\rho(I_{\ell})$ contains an element of order $3$. 

\prop[elliptic-curves-cartan-eltsord3]{Let $\ell$ be a prime, and let $J_{\ell}=\rho(I_{\ell})$. If $E$ has potentially multiplicative reduction at $\ell$, then $J_{\ell} \subset \{\pm \mrm{id}_2\}$ and $\ell \geq 2p+1$. If $E$ has potentially good reduction at $\ell$:
\begin{itemize}[noitemsep,label=\tiny$\bullet$]
\item If $\ell=p$, then $J_{\ell}$ is properly contained in $C$ if and only if either $p \equiv 2\mod{9}$ and $v_p(\Delta) \in \{4,10\}$, or $p \equiv 5 \mod{9}$ and $v_p(\Delta) \in \{2,8\}$.
\item Assume that $\ell \neq p$. If $3 \nmid v_{\ell}(\Delta)$, then $J_{\ell}$ contains an element of order $3$, and the converse holds when $\ell \geq 5$. 
\item If $J_{\ell}$ contains an element of order $3$ and $\ell \neq 3$, $K$ is unramified at $\ell$. 
\item When $\ell=3$, $J_{\ell}$ contains an element of order $3$ if and only if the conductor of $E$ is divisible by $27$, if and only if $\Q_3(E[2])/\Q_3$ is wildly ramified.  
\item When $\ell=2$, $J_{\ell}$ contains an element of order $3$ if and only if the ramification index of $\Q_2(E[3])/\Q_2$ is not a power of $2$. In this case, $J_{\ell}$ is cyclic of order $3$ or $6$. 
\end{itemize}
}

\demo{When $\ell=p$, $E$ has potentially good supersingular reduction by Proposition \ref{elliptic-curves-cartan-padic}, and this result also states that $3C \subset J_{p} \subset C$. In particular, $\rho_{|I_p}$ is irreducible. By \cite[(2.2)]{Conj-Serre}, the index of $J_p$ in $C$ is exactly the greatest common divisor of $k-1$ and $p^2-1$, where $k$ is the Serre weight of $\rho$. This weight is computed in \cite[Th\'eor\`eme 1]{Kraus-thesis}; since Proposition \ref{elliptic-curves-cartan-padic} implies that one has $4 \leq v_p(j(E)) = 3v_p(c_4)-v_p(\Delta)$, where $c_4$ is the classical invariant attached to a minimal model of $E/\Q$, the result follows from a simple case-by-case study. 

When $E$ has potentially multiplicative reduction at $\ell$, then $\ell \neq p$ by the above, so $\ell \equiv \pm 1 \mod{p}$ by Proposition \ref{elliptic-curves-cartan-nram}. Since $p-1,p+1 > 2$ are even and $2p-1 > 3$ is divisible by $3$, none of them equals $\ell$, so one has $\ell \geq 2p+1$. 

We assume from now on that $E$ has potentially good reduction at $\ell$. By \cite[\S 5.6]{Serre-image}, $|J_{\ell}|$ is divisible by $\frac{12}{(12,v_p(\Delta))}$, so if $3 \nmid v_{\ell}(\Delta)$, $J_{\ell}$ contains an element of order $3$. Assume now that $\ell \neq 3$ and $J_{\ell}$ contains an element of order $3$; by \cite[Cor. 2 to Th. 2]{GoodRed}, $J_{\ell} \subset N$ is isomorphic to a subgroup of $\SL{\F_3}$. Since $\SL{\F_3}$ does not embed in $N$, $J_{\ell}$ is isomorphic to a proper subgroup of $\SL{\F_3}$ with an element of order $3$, so one checks that $J_{\ell}$ is cyclic of order $3$ or $6$. Since $p \equiv -1\mod{3}$, one checks that cyclic subgroups of $N$ with order $3$ or $6$ are contained in $C$ -- so that $J_{\ell} \subset C$ and $K$ is unramified at $\ell$. 

When $\ell \geq 5$, $J_{\ell}$ is cyclic with order $|J_{\ell}|=\frac{12}{(12,v_p(\Delta))} < 12$ properly dividing $12$ by \cite[\S 5.6 $a_1$)]{Serre-image}, so $|J_{\ell}|$ is divisible by $3$ if and only if $3 \nmid v_{\ell}(\Delta)$. When $\ell=2$, one concludes in a similar way using the above and \cite[\S 5.6 $a_3$)]{Serre-image}. When $\ell=3$, the claim follows from \cite[Cor. 2 to Th. 2]{GoodRed}, which implies in particular that whether $\Q(E[n])$ is wildly ramified at $3$ does not depend on the choice of $n \geq 3$ prime to $3$. To conclude, it is enough to note that $\mrm{Gal}(\Q(E[4])/\Q(E[2]))$ is a $2$-group (hence there is no wild ramification in the extension above any odd prime), and that if $J_3$ contains an element of order $3$, the tame conductor exponent of $E[p]$ (hence of $E$) at $3$ is equal to $2$. }

\prop[elliptic-curve-cartan-plus-irreducible]{The representation $\rho'$ is irreducible if and only if $\rho: G_{\Q} \rar N$ is surjective. This is the case if any of the following conditions is satisfied. 
\begin{itemize}[noitemsep,label=\tiny$\bullet$]
\item there exists a prime $\ell$ such that $\rho(I_{\ell})$ contains an element of order $3$, 
\item there exists a prime $\ell \neq p$ at which $E$ has good reduction and such that $a_{\ell}(E)^2\equiv \ell \pmod{p}$,
\item $E$ is not CM and $p > 37$.
\end{itemize}
The second condition is necessary. If the class number of $K$ is prime to $3$ (for instance if $E$ has CM by $K$), the first condition is necessary. 
}

\demo{Every proper subgroup of $N/3C \simeq \mfk{S}_3$ is abelian, so the restriction of the irreducible representation $\mrm{Ind}_{C/3C}^{N/3C}{\psi}$ to any proper subgroup of $N/3C$ is reducible. Hence $\rho'$ is irreducible if and only if $\rho: G_{\Q} \rar N$ is surjective, since by Proposition \ref{elliptic-curves-cartan-padic} one has $\rho(G_{\Q}) \supset 3C$. When the last condition holds, $\rho: G_{\Q} \rar N$ is surjective by the main theorems of \cite{LFL,FL}.

Let $M$ be the image of $\rho$. Since $3C \subset M \not\subset C$ and $p \equiv 2,5\mod{9}$, $\rho$ is surjective if $M$ contains a scalar multiple of some element of order $3$, so the first condition is sufficient.  
For $g \in N$, $g$ (resp. some scalar multiple of $g$) has order $3$ if and only if $(\mrm{Tr}(g),\det{g})=(-1,1)$ (resp. $\mrm{Tr}(g)^2=\det{g}$). Hence Chebotarev's theorem implies that the second condition is necessary and sufficient.

If the class number of $G_K$ is prime to $3$, then the index of the closed subgroup of $G_K^{ab}$ generated by its inertia subgroups is prime to $3$: therefore, $\rho(G_{\Q})$ contains an element of order $3$ if and only if $\rho(I_{\lambda})$ contains an element of order $3$ for some prime ideal $\lambda$ of $\OO_K$, if and only if $\rho(I_{\ell})$ contains an element of order $3$ for some prime $\ell$.} 

We now determine the Artin conductor of $\rho'$.

\prop[artin-conductor-rhoprime]{Let $N_E$ denote the conductor of $E$. Let $M$ be the Artin conductor of $\rho'$ and $\ell$ be any prime.
\begin{itemize}[noitemsep,label=$-$]
\item If $E$ has good or potentially multiplicative reduction at $\ell \neq p$, then $v_{\ell}(M)=0$.
\item If $E$ has potentially good reduction at $\ell \neq 3$, then $v_{\ell}(M)$ equals $2$ if $\rho(I_{\ell})$ contains an element of order $3$, $1$ if $K$ is ramified at $\ell$, and $0$ otherwise.  
\item If $\ell=3$ is a prime of potentially good reduction for $E$, $v_{\ell}(M)$ equals $v_{\ell}(N_E)$ if $27 \mid N_E$, $1$ if $27 \nmid N_E$ and $K$ is ramified at $3$, and $0$ otherwise.   
\end{itemize}
Note that $v_{\ell}(M)$ does not change if $E$ is replaced with any quadratic twist. 
}

\demo{Let $\ell$ be any prime, and let $\psi' = \psi \circ \rho^{\ast}_{|G_K}: G_K \rar \C^{\ast}$, whose image has cardinality dividing $3$. For any prime $\lambda$ of $K$, let $n(\lambda)$ be the conductor exponent of $\psi'$ at $\lambda$. By \cite[Proposition VII.11.2 (ii)]{Neukirch-ANT}, the conductor exponent of $\rho'$ at $\ell$ is
\begin{itemize}[noitemsep,label=\tiny$\bullet$]
\item $2n(\lambda)$ if $\lambda$ is the unique prime above $\ell$, and $\ell$ is inert in $K$,
\item $n(\lambda)+v_{\ell}(D)$ if $\ell$ is ramified in $K$ and $\lambda$ is the unique prime above $\ell$, 
\item $n(\lambda)+n(\lambda')=2n(\lambda)=2n(\lambda')$ if $\lambda,\lambda'$ are the two distinct prime ideals of $K$ above $\ell$ (because $\psi'$ is conjugate to its inverse). 
\end{itemize}
 
It is clear that $K$ does not change if $E$ is replaced with a quadratic twist; because $\psi(-\mrm{id})=1$, neither does $n(\lambda)$. Hence $M$ does not change if $E$ is replaced with a quadratic twist.  

The result follows directly from Propositions \ref{elliptic-curves-cartan-padic}, \ref{elliptic-curves-cartan-nram}, \ref{elliptic-curves-cartan-eltsord3}, except in the case where $\ell=3$ is a prime of bad and potentially good reduction. Let $\Gamma \leq G_{\Q}$ be a pro-$3$-group; then $\rho(\Gamma)$ is a subgroup of $N$ which is a $3$-group, so it is contained in the $3$-torsion subgroup of $C$, which is cyclic of order $3$. Therefore, for $\gamma \in \Gamma$, $\rho(\gamma)$ has a non-trivial fixed vector if and only if $\rho(\gamma)$ is the identity, if and only if $\rho'(\gamma)$ is the identity if and only if $\rho'(\gamma)$ has a non-trivial fixed vector. 

Hence, if $27 \mid N_E$, $\Q_3(E[p])$ is wildly ramified by Proposition \ref{elliptic-curves-cartan-potgood}, and therefore $E[p](\overline{\Q}_3)$ and $\rho'$ have the same conductor exponents at $3$, so the claim holds. 

If $27 \nmid N_E$, then $\Q_3(E[p])$ is tamely ramified, so $\rho_{|I_3}$ factors through the tame quotient, hence $\rho(I_3)\cap C \leq C$ has order prime to $3$. Therefore $\psi'$ is unramified above $3$, whence the conclusion. }

\prop[ep-is-good]{Let $n$ be a non-zero integer. Let $E_{n}$ be the elliptic curve with equation $y^2=x^3+n$. 
Assume that: 
\begin{itemize}[noitemsep,label=\tiny$\bullet$]
\item $n$ is not divisible by $16$, $81$, or the sixth power of a prime number,
\item if $p \equiv 2 \mod{9}$, $v_p(n) \in \{2,5\}$,
\item if $p \equiv 5 \mod{9}$, $v_p(n) \in \{1,4\}$.
\end{itemize}

Then $E_{n}$ satisfies the assumptions at the beginning of the Section (whose notations we adopt) with $K=\Q(\sqrt{-3})$. Moreover,
\begin{itemize}[noitemsep,label=\tiny$\bullet$]
\item the above equation is minimal with discriminant $-2^43^3n^2$, 
\item $\rho$ is surjective (equivalently, $\rho'$ is irreducible) if and only if $n$ is not of the form $2p^{\alpha}a^3$,
\item $\rho(I_p) \subset 3C$,
\item the conductor of $\rho'$ is $2^{e_{2,n}} \cdot 3^{e_{3,n}} \cdot L^2$, where $L$ is the product of the primes $\ell \mid n$ with $\ell \nmid 6p$ such that $3 \nmid v_{\ell}(n)$, $e_{3,n}$ is the $3$-adic valuation of the conductor of $E_n$ if it is divisible by $27$ (equivalently, if $n \notin \pm 27^{\Z}(1+9\Z_3)$), and $e_{3,n}=1$ otherwise, and $e_{2,n}$ is $0$ if $v_2(n)=1$ and $2$ otherwise.   
\end{itemize}}

\demo{The assumptions of the beginning of the Section are satisfied because $E$ has (potential) complex multiplication by $K=\Q(\sqrt{-3})$ and $p$ is inert in $K$. In particular, $j(E)=0$, so $E$ has potentially good reduction everywhere. The value of the discriminant is a direct computation; the assumptions imply that it is not divisible by the $12$-th power of a prime, so the equation is minimal by \cite[Remark VII.1.1]{AEC1}. The claim about the image of an inertia subgroup at $p$ follows from Proposition \ref{elliptic-curves-cartan-eltsord3}. 

We claim that $\rho(I_2)$ is contains an element of order $3$ if and only if $v_2(n)=1$ (this also implies the claimed value for $e_{2,n}$ by Proposition \ref{artin-conductor-rhoprime}). Indeed, if $v_2(n)\neq 1$, then $v_2(n) \in \{0,2,3\}$, so $v_2(-16\cdot 27n^2)$ is not divisible by $3$, so Proposition \ref{elliptic-curves-cartan-eltsord3} implies that $\rho(I_2)$ contains an element of order $3$. If instead $v_2(n)=1$, we can write $n=2n'$ with $n' \in \Z_2^{\times}$. Because $\Z_2^{\times}$ is a pro-$2$-group, $n'$ is a cube in $\Z_2^{\times}$, so $(E_n)_{\Q_2}$ is a quadratic twist of $y^2=x^3+16$, which has good reduction at $2$. Therefore $\rho(I_2) \subset \{\pm \mrm{id}\}$, so $\rho(I_2)$ does not contain an element of order $3$. 

By Proposition \ref{elliptic-curves-cartan-eltsord3}, $\rho(I_3)$ contains an element of order $3$ if and only if the conductor of $E_n$ is divisible by $27$, if and only if $\Q_3(E_n[2])/\Q_3$ is wildly ramified, if and only if $\Q_3(\tau,\sqrt[3]{n})/\Q_3$ is wildly ramified (where $\tau$ is a primitive third root of unity), if and only if $\Q_3(\sqrt[3]{n})/\Q_3$ is totally ramified of degree $3$. This is true if and only if $n \in \Q_3^{\times}$ is not a cube, which is true if and only if $n \in \pm 27^{\Z}(1+9\Z_3)$. The conductor exponent of $\rho'$ at $3$ is correct by Proposition \ref{artin-conductor-rhoprime}.  

The conductor exponent of $\rho'$ at primes $\ell \geq 5$ is correct by Proposition \ref{artin-conductor-rhoprime}.  

Only the claim about the surjectivity of $\rho$ remains. By Proposition \ref{elliptic-curve-cartan-plus-irreducible}, and since $\Q(\sqrt{-3})$ has class number one, $\rho$ is surjective if and only if there is a prime $\ell$ (of bad reduction) for $E_n$ such that $\rho(I_{\ell})$ contains an element of order $3$. We have already established that this condition is equivalent to 
\begin{itemize}[noitemsep,label=\tiny$\bullet$]
\item that for $\ell = p$, it is never satisfied by construction,
\item for $\ell \notin \{2,3,p\}$, that $v_{\ell}(n)$ is not divisible by $3$,
\item for $\ell=3$, that $n$ is not a cube in $\Q_3$,
\item for $\ell=2$, that $v_{2}(a) \neq 1$. 
\end{itemize}
By construction, $2p^{v_p(n)} \equiv \pm 1\mod{9}$, so no $\rho(I_{\ell})$ contains an element of order $3$ if and only if $n$ is of the form $2p^{v_p(n)}a^3$ for some integer $a$. 
}

\cor[ep-compute-g]{Suppose that $E$ is a quadratic twist of the elliptic curve $E_p$ (if $p \equiv 5 \mod{9}$) or $E_{p^2}$ (if $p \equiv 2 \mod{9}$) as in Proposition \ref{ep-is-good}. Then the irreducible odd Artin representation $\rho'$ is attached to the unique weight one newform $g\in \mathcal{S}_1(\Gamma_1(108))$. It has rational coefficients, CM by $\Q(\sqrt{-3})$, and its character has conductor $3$.}

\demo{This is a direct consequence of Proposition \ref{ep-is-good}, and a verification using the LMFDB \cite{lmfdb}.  }

\subsection{Non-vanishing of the $L$-value}
\label{subsect:lvalue-quaternion-supersing}

In this Section, we are interested in exploiting an explicit Waldspurger formula of Cai, Shu and Tian \cite{CaiShuTian} to prove, in certain cases, a non-vanishing result for the special $L$-values appearing in Proposition \ref{fimm-main}. The construction of \cite{CaiShuTian} (and some later results in this Section) involve some arithmetic of quaternion algebras, and we refer the reader to \cite{Voight} for a thorough discussion of their properties that we will freely use.

\lem[anticyc-chars]{Let $K$ be an imaginary quadratic field with absolute discriminant $D$ and $\theta$ be a continuous complex character of $G_K$ such that the composition $G_{\Q}^{ab} \overset{\mrm{Ver}}{\rar} G_K^{ab} \overset{\theta}{\rar} \C^{\times}$ is trivial. 
Let $\mfk{a}_K: \mathbb{A}_K^{\times}/K^{\times} \rar G_K^{ab}$ the class field theory map, normalized so as to send a uniformizer to an arithmetic Frobenius. 
Then the conductor ideal $\mfk{c}$ of $\theta$ is generated by a positive integer $c$, the Artin conductor of $\mrm{Ind}_K^{\Q}{\theta}$ is $|D|c^2$, and $\theta \circ \mfk{a}_K$ factors as a character $\frac{\mathbb{A}_K^{\times}}{K^{\times}[(\Z+c\OO_K) \otimes \hat{\Z}]^{\times}} \rar \C^{\times}$. }

\demo{This is standard, so we only sketch the proof. By \cite[Proposition VII.11.7 (iii)]{Neukirch-ANT}, the Artin conductor of $\mrm{Ind}_K^{\Q}{\theta}$ is $|D|\mathbf{N}(\mfk{c})$. The assumption implies that $\theta \circ \mfk{a}_K$ vanishes on $\mathbb{A}_{\Q}^{\times}$. 

Because $\theta$ is $\mrm{Gal}(K/\Q)$-conjugate to $\theta^{-1}$, $\mfk{c}$ is stable under $\mrm{Gal}(K/\Q)$. What remains to be shown is that $\mfk{c}$ is generated by some positive integer $c > 0$; it is enough to show that if $K$ is ramified at $\ell$ and $\lambda$ is the unique prime of $K$ above $\ell$, the conductor exponent at $\lambda$ of $\theta$ is even. This reduces to the following elementary claim: if $K/\Q_{\ell}$ is a quadratic ramified extension, $\varpi$ is a uniformizer and $n \geq 0$, then $\Z_{\ell}^{\times}(1+\varpi^{2n+1}\OO_L) = \Z_{\ell}^{\times}(1+\varpi^{2n}\OO_L)$. }

Let $K$ be an imaginary quadratic field and $p \geq 11$ a prime number inert in $K$. Let $\psi: G_K \rar \C^{\times}$ be a character unramified at $p$ with order exactly $3$ and such that $G_{\Q}^{ab} \overset{\mrm{Ver}}{\rar} G_K^{ab} \overset{\psi}{\rar} \C^{\times}$ is zero. 
Then $\psi$ defines a representation $\rho = \mrm{Ind}_K^{\Q}{\psi}: G_{\Q} \rar \GL{\C}$ which is odd, irreducible with image $\mfk{S}_3$, and unramified at $p$, so it is attached to some weight one newform $g \in \mathcal{S}_1(\Gamma_1(N))$ of level prime to $p$.

Let $\OO \subset K$ be an order with class field $H$ such that $\psi$ factors through $\mrm{Gal}(H/K)$, and let $\mfk{a}: \mrm{Pic}(\OO) \rar \mrm{Gal}(H/K)$ be the class field theory isomorphism, normalized so that it maps a prime ideal $\lambda$ not dividing the conductor of $\OO$ to the arithmetic Frobenius $\Fr_{\lambda}$.

Note that $\OO=\Z+c\OO_K$, where $c$ is a generator for the conductor ideal of $\psi$, by Lemma \ref{anticyc-chars}, so that $-N$ is the discriminant of $\OO$.  

\prop[epn-caishutian]{ Let $B/\Q$ be the unique quaternion algebra exactly ramified at $p\infty$, and fix an embedding $K \rar B$. Let $R \subset B$ be a maximal order\footnote{Equivalently, by \cite[Theorems 14.1.3, 15.5.5]{Voight}, with reduced discriminant $p$.} such that $R \cap K=\OO$. Let $\varphi: \mrm{Pic}(\OO) \rar \mrm{Cls}(R)$ (i.e. the class set of invertible right fractional ideals of $R$ up to equivalence, cf. \cite[Definition 17.3.4]{Voight}) be given by $I \mapsto IR$. 

The following are equivalent:
\begin{itemize}[noitemsep,label=\tiny$\bullet$]
\item There exists $x \in \mrm{Cls}(R)$ such that the three subsets $\varphi^{-1}(x) \cap (\psi \circ \mfk{a})^{-1}(\omega)$, where $\omega$ runs through all third roots of unity, do not all have the same size. 
\item There exists $f \in \mathcal{S}_2(\Gamma_0(p))$ such that $L(f,g,1) \neq 0$.
\end{itemize}
}

\demo{This is a translation of \cite[Theorem 1.2]{CaiShuTian}. More precisely, \emph{loc.cit.} is the following result. 
The set $X :=B^{\times} \backslash (B \otimes \hat{\Z})^{\times}/(R \otimes \hat{\Z})^{\times}$ is finite (by \cite[Corollary 27.6.20]{Voight}). The free $\C$-vector space $\C[X]$ over $X$ has a linear form $\deg: \C[X] \rar \C$ mapping any $x \in X$ to $1$, a Hermitian pairing under which the $x \in X$ form an orthogonal basis, an action of operators $T_q$ and $S_q$ (for $q \neq p$ prime). The chosen embedding $K \rar B$ induces a map 

\[\varphi_0: \operatorname{Pic}(\OO) = K^{\times} \backslash (K \otimes \hat{\Z})^{\times} / (\OO \otimes \hat{\Z})^{\times} \rar B^{\times} \backslash (B \otimes \hat{\Z})^{\times} / (R \otimes \hat{\Z})^{\times} = X.\] 

Let $v_{\psi} = \sum_{t \in \mrm{Pic}(\OO)}{(\psi\circ \mfk{a})(t)[\varphi_0(t)]} \in \C[X]^0 := \ker{\deg}$. Then, for any $f \in \mathcal{S}_2(\Gamma_0(p))$, there is a unique line $L_f \subset \C[X]^0$ on which, for every $q \neq p$ prime, $S_q$ acts trivially and $T_q$ acts by multiplication by $a_q(f)$, and the central $L$-value $L(f,g,1)$ is a positive multiple of the norm of the orthogonal projection of $v_{\psi}$ on $L_f$. \footnote{If we regard the anticyclotomic character of \cite{CaiShuTian} as coming directly from $\mrm{Pic}(\OO)$, we do not need any normalization for the class field theory reciprocity map required to state Theorem 1.2 of \emph{loc.cit.}. Regardless, since $f,g$ have real coefficients, the outcome is independent from the choice of normalization.}

Now, $X$ identifies with $\mrm{Cls}(R)$ (cf. e.g. \cite[Lemma 27.6.8]{Voight}), and a careful reading of \cite[\S\S 1--5]{Gross-Heights}\footnote{Or, more bluntly, the Jacquet--Langlands correspondence (see e.g. \cite[Theorem 10.5]{Gelbart}).} shows that $\C[X]^0$ identifies as a $\C[(T_q)_{q \neq p},(S_q)_{q \neq p}]$-module with $\mathcal{S}_2(\Gamma_0(p))$ (where $T_q$ corresponds to the classical Hecke operator and $S_q$ acts trivially), whence the conclusion.}

\prop[injective-map]{Let $R$ be a maximal order in a quaternion algebra $B$ over $\Q$ ramified only at $p,\infty$. Let $\OO$ be a quadratic order of discriminant $-D$ optimally embedded in $R$. Assume that $p > D$. Then the map $I \in \mrm{Pic}(\OO) \mapsto IR \in \mrm{Cls}(R)$ is injective.}

\demo{The assumptions imply that $L := \mrm{Frac}(\OO)$ is inert at $p$. Let $I,J \subset L$ be two fractional ideals in two distinct ideal classes such that $\alpha IR=JR$ for some $\alpha \in B^{\times}$. Exactly as in the proof of the Minkowski bound (e.g. \cite[Theorem V.4]{Lang-ANT}), we may assume that $J=I'I$, where $I,I' \subset \OO$ are invertible ideals, and $I'$ is non-principal with norm at most $\frac{2}{\pi}\sqrt{D}$. Now, $\OO$ and $\alpha \OO \alpha^{-1}$ are two quadratic orders contained in the left order $R'$ of $I'IR=\alpha IR$ (which is a maximal order of $B$). By \cite[Theorem 2']{Kaneko}, if their respective saturations $\tilde{\OO},\tilde{\OO}'$ in $R'$ are distinct, one has \[p \leq \sqrt{\mrm{disc}(\tilde{\OO})\mrm{disc}(\tilde{\OO}')} \leq \sqrt{\mrm{disc}(\OO)\mrm{disc}(\alpha\OO\alpha^{-1})} = D,\] 
a contradiction. Hence $\tilde{\OO}=\tilde{\OO}'$, hence $\alpha \OO \alpha^{-1}=\OO$ (since they have the same discriminant). Therefore, either $\alpha \in L$, or one has $\alpha z \alpha^{-1}=\overline{z}$ for every $z \in L$. If $\alpha \in L$, then one finds $\alpha^{-1}I'R=R$, hence $\alpha^{-1}I'=\OO$, a contradiction. 

Thus, one has $\alpha z \alpha^{-1} =\overline{z}$ for every $z \in L$. Since $\alpha (IR) =I'IR \subset IR$, it is an integral element. Moreover, $\alpha^2$ commutes to $L$, so $\alpha^2 \in L \cap \Q(\alpha)=\Q$, hence $\alpha^2=-c$ for some positive integer $c$. Hence $c^2=\mrm{nrd}(\alpha)^2=\left|\frac{IR}{I'IR}\right|=(\mathbf{N}I')^2 \leq \frac{4}{\pi^2}D$, thus $c \leq \frac{2\sqrt{D}}{\pi}$.  

Now, we argue as in the proof of \cite[Theorem 2']{Kaneko} for the orders $\OO$ and $\Z[\alpha]$. Let $1,u=\frac{D+\sqrt{-D}}{2}$ be a $\Z$-basis of $\OO$, and let $R' \subset B$ be the order with $\Z$-basis $1,u,\alpha,\alpha u=\overline{u}\alpha$. Since the reduced trace of $\alpha$ and $\alpha u$ is zero, the discriminant of $R'$ is 
\[\det\begin{pmatrix} 2 & -D & 0 & 0\\-D & \frac{D^2-D}{2} & 0 & 0\\0 & 0 & -2c & -cD \\ 0 & 0 & -cD & -c\frac{D^2-D}{2}\end{pmatrix} = D^2c^2,\] which is not divisible by $p$, which is impossible since $B$ is ramified at $p$.}

\rem{Proposition \ref{injective-map} shows that the equivalent conditions of Proposition \ref{epn-caishutian} are satisfied when $p > N$. The non-vanishing was already known to hold when $-N$ is an odd fundamental discriminant (cf. \cite[Corollary 1]{Michel-Ramakrishnan}), but this is almost never the case in our applications (in Corollary \ref{ep-compute-g}, one has $N=108$). }

This result can also be interpreted in terms of supersingular elliptic curves through the Deuring correspondence.

\prop[supersing-optimal]{The group $\mrm{Gal}(H/K)$ acts simply transitively on the set of primes of $\OO_H$ dividing $p$. Moreover, for every intermediate extension, $K \subset H' \subset H$ and every prime $\mfk{p}$ of $H$ dividing $p$, $\OO/p\OO \rar \OO_K/p\OO_K \rar \OO_{H'}/\mfk{p}\cap \OO_{H'} \rar \OO_H/\mfk{p}\OO_H$ is a sequence of isomorphisms.

Let $F \supset H$ be a number field contained in $\C$, $\mfk{p} \subset \OO_F$ a prime ideal above $p$, and $E/\OO_F$ an elliptic curve with complex multiplication by $\OO$, i.e. $\mrm{End}_{\OO_F}(E) = \mrm{End}_{F}(E)=\mrm{End}_{\C}(E_{\C})=\OO$. \footnote{This is possible by CM theory (e.g. \cite[Theorem 11.1]{Cox}), classical approximation results \cite[Th\'eor\`eme 8.8.2]{EGA-IV3}, Lemma \ref{endo-abvar-spe}, and the existence of N\'eron models \cite[\S 1.4/3]{BLR}.} We have an embedding $\OO/p\OO \simeq \OO_H/\mfk{p}\OO_H \rar \OO_F/\mfk{p}\OO_F \rar \overline{\F_p}$; where the last map is fixed once and for all. Then $\mrm{End}(E_{\overline{\F_p}})$ is a maximal order in a quaternion algebra over $\Q$ ramified at $p$ and $\infty$. Moreover, the embedding $r: \OO = \mrm{End}_{\OO_F}(E)= \mrm{End}_{\OO_{F,\mfk{p}}}(E_{\OO_{F,\mfk{p}}}) \rar \mrm{End}(E_{\overline{\F_p}})$ is optimal.

Finally, the following diagram, where both vertical maps are isomorphisms, is commutative: 

\[
\begin{tikzcd}[ampersand replacement=\&]
\mrm{Pic}(\OO) \arrow{d}{\mfk{a}} \arrow{rrr}{\overline{\varphi}: I \mapsto \overline{I}R} \& \& \&\mrm{Cls}(R)\\
\mrm{Gal}(H/K) \arrow{rrr}{\sigma \mapsto \sigma(j(E)) \mod{\mfk{p}}} \& \& \&\mathcal{S} \arrow{u}{\gamma}
\end{tikzcd}
\]
where $\mathcal{S} \subset \overline{\F_p}$ is the set of supersingular $j$-invariants, and $\gamma$ sends $j_1 \in \mathcal{S}$ to $\mrm{Hom}_{\overline{\F_p}}(E_{\overline{\F_p}},E_1)$, where $E_1/\overline{\F_p}$ is the elliptic curve with $j$-invariant $j_1$.
}

\demo{First, $p$ is inert in $K$, so $p\OO_K \cap \OO=p\OO$ is a prime ideal of $\OO$, and therefore the prime ideal $p\OO_K$ splits totally in $H$. As a consequence, $\mrm{Gal}(H/K)$ acts simply transitively on the set of primes of $\OO_H$ lying above $p\OO_K$, and $\OO/p\OO \rar \OO_K/p\OO_K \rar \OO_H/\mfk{p}\OO_H$ is a composition of isomorphisms. Since $p$ is inert in $K$, $E_{\overline{\F_p}}$ is supersingular, so, by \cite[Theorem 42.1.9]{Voight}, $R := \mrm{End}(E_{\overline{\F_p}})$ is a maximal order in the quaternion algebra $B := R \otimes \Q$ over $\Q$ ramified exactly at $p$ and $\infty$. 

The map $r$ is injective by Lemma \ref{endo-abvar-spe}. Its cokernel contains no $p$-torsion, because $\OO \otimes \Z_p$ is a discrete valuation ring while $R \otimes \Z_p$ is a finite free $\Z_p$-module without $\OO \otimes \Z_p$-torsion (since $B$ is ramified at $p$). If $\ell \neq p$ is a prime, the specialization map $\Tate{\ell}{E_F} \rar \Tate{\ell}{E_{\OO_F/\mfk{p}}}$ is an isomorphism of $\Z_{\ell}$-modules, so the cokernel of $r$ has no $\ell$-torsion by \cite[Theorem 12.5]{MilAb}.

Now, we discuss the diagram appearing in the statement of the Proposition. The map $\gamma$ is well-defined and bijective by \cite[Theorem 42.3.2, Remark 42.3.3]{Voight}. By the usual reciprocity law (e.g. \cite[Theorem 11.36]{Cox}), $j(E) \in H$ and if $I \subset \OO$ denotes an invertible ideal, then the elliptic curve with $j$-invariant $\mfk{a}(I)(j(E))$ is precisely $E/E[I]$, so we are done if $\mrm{Hom}_{\overline{\F_p}}(E_{\overline{\F_p}},(E/E[I])_{\overline{\F_p}})$ and $\overline{I}R$ are isomorphic as right $R$-modules. Let $\pi: E \rar E/E[I]$ be the projection: by \cite[Lemma 42.2.7]{Voight}, since $E[I]_{\overline{\F_p}} = E_{\overline{\F_p}}[RI]$, the following map is an isomorphism of left $R$-modules: 
\[u \in \mrm{Hom}((E/E[I])_{\overline{\F_p}},E_{\overline{\F_p}}) \mapsto u \circ \pi \in RI.\]

Because taking dual isogenies in $R = \mrm{End}(E_{\overline{\F_p}})$ is the standard involution, the following map is an isomorphism of right $R$-modules: 
\[u\in \mrm{Hom}(E_{\overline{\F_p}},(E/E[I])_{\overline{\F_p}}) \mapsto \pi^{\vee} \circ u^{\vee} \in (RI)^{\vee}=\overline{I}R,\]
 so we are done. 
}

\cor[epn-supersing]{Let $j \in H$ be the $j$-invariant of any elliptic curve with complex multiplication by the order $\OO$, and fix a prime $\mfk{p} \subset \OO_H$ above $p$. Let $K \subset H' \subset H$ be the field such that $G_{H'}=\ker{\psi}$ and $\mfk{p}'=\mfk{p}\cap \OO_{H'}$. The following conditions are equivalent: 

\begin{itemize}[noitemsep,label=\tiny$\bullet$]
\item There exists $f \in \mathcal{S}_2(\Gamma_0(p))$ such that $L(f,g,1) \neq 0$.
\item $\sum_{\sigma \in \mrm{Gal}(H/K)}{\psi(\sigma)[\sigma(j) \mod{\mfk{p}}]} \neq 0$ (where the sum is computed in the $\C$-vector space with basis $\OO_K/p\OO_K$), 
\item There exists $j_p \in \OO_K/p\OO_K$ such that the three sets $\{\sigma \in \mrm{Gal}(H/K) \mid \psi(\sigma)=\omega, \sigma(j) \mod{\mfk{p}}=j_p\}$ (where $\omega$ runs through all third roots of unity) do not have the same cardinality.    
\end{itemize}

Let $f_K \in \OO_K[t]$ and $f_{H'} \in \OO_{H'}[t]$ the minimal polynomials of $j(E)$ over $K$ and $H'$ respectively. Then the above conditions are satisfied if $(f_K(t) \mod{p}) \in \frac{\OO_K}{p\OO_K}[t]$ and $(f_{H'}(t) \mod{\mfk{p}'})^3 \in \frac{\OO_{H'}}{\mfk{p}'}[t] \simeq \frac{\OO_K}{p\OO_K}[t]$ are distinct. In particular, they are satisfied if the reduction modulo $p$ of the Hilbert class polynomial attached to the order $\OO$ (see \cite[\S 2.2, 2.3]{Cremona-Sutherland}) is not a cube.  
}

\demo{The equivalence is a direct consequence of Proposition \ref{epn-caishutian} thanks to Proposition \ref{supersing-optimal} (because the only non-trivial additive relation between the third roots of unity is $1+\tau+\tau^2=1$, where $\tau$ is a primitive third root of unity). Furthermore, the third condition is equivalent to the following claim: the formal linear combinations $\sum_{\sigma \in C}{[\sigma(j) \mod{\mfk{p}}]}$ (over the three cosets $C$ of $\mrm{Gal}(H/H')$ in $\mrm{Gal}(H/K)$) are equal. It implies in particular that \[f_K(t) \mod{\mfk{p}} = \left(\prod_{\sigma \in \mrm{Gal}(H/H')}{(t-\sigma(j))}\right)^3 \mod{\mfk{p}} = f_{H'}(t)^3 \mod{\mfk{p}'},\] and therefore we are done. }

\cor[factor-hilbclasspoly]{Let $-D < 0$ be a discriminant whose associated fundamental discriminant corresponds to the imaginary quadratic field $K$. Let $p \geq 11$ be a prime inert in $K$. Assume that the Hilbert class polynomial $H_{-D}(t)$ is not separable modulo $p$. Then $p \leq D$. }

\demo{This is a direct consequence of Proposition \ref{supersing-optimal} and Proposition \ref{injective-map}. }

\rem{Let us consider couples $(p,M)$ of a following form: $p \equiv 2,5 \mod{9}$ is a prime not dividing $M$, $-M$ is the discriminant of some order in an imaginary quadratic field (of discriminant $-D$) in which $p$ is inert. We say that $(p,M)$ is \emph{bad} if the Hilbert class polynomial of discriminant $-M$ is a cube modulo $p$. We say that $(p,M)$ is very bad if there exists a Galois representation $G_{\Q} \rar \GL{\C}$ with Galois image $\mfk{S}_3$, conductor $M$, complex multiplication by $\Q(\sqrt{-M})$ and the same character as $\Q(\sqrt{-M})$ such that $L(f,\rho,1)=0$ for every $f \in \mathscr{S}$. 

The results of this Section imply that a very bad pair is bad, and that if a pair $(p,M)$ is bad, one has $p \leq M$. Much stronger results seem to hold empirically: for instance, 
\begin{itemize}[noitemsep,label=$-$]
\item when $M \leq 10000$, there are at most $585$ bad pairs (out of $560380$ possible pairs with $p \leq M$), all of which satisfy $p \leq 3\sqrt{M}$, but there are at most $29$ very bad pairs, all of which (except the pair $(47,2039\cdot 2^2)$) verify $p \leq 29$.
\item when $D=3$, $M= Dc^2$ with $c=3^ic' \leq 300$, where $i \leq 1$ and $c'$ is square-free (cf. Proposition \ref{result-for-en} for the significance of this situation), the set of bad pairs $(p,M)$ is $\{(11,3 \cdot 271^2)\} \cup \Big[\{11,23\} \times \{3 \cdot (6 \cdot 29^2), 3 \cdot (9 \cdot 29)^2, 3 \cdot (10 \cdot 29)^2\}\Big]$. 
\end{itemize}
This suggests that our results are not optimal, but we do not currently know of a way to improve them. }

\subsection{Proofs of the concrete results}
\label{subsect:concrete-results}

\prop[explicit-fimm]{In the situation of Corollary \ref{fimm-invariants}, if in addition $p \leq 10^4$, then the conductor of $E$ is $p^2N_{\rho}$ and one has $j(E) \in \Z$. 
}

\demo{We already saw that there were no exceptional primes for $p$, since $11 \leq p \leq 10^4$.}

\rem{In the notation of Section \ref{subsect:explicit-rhoprime}, the assumptions of Lemma \ref{anticyc-chars} and more generally of Section \ref{subsect:lvalue-quaternion-supersing} hold for the character $\psi': \psi \circ \rho^{\ast}_{|G_K}: G_K \rar \C^{\times}$. Indeed, $\psi^3=1$, so $\psi'$ is continuous with finite image of cardinality $1$ or $3$ (hence odd). Moreover, conjugation by $\mrm{Gal}(K/\Q)$ acts on $\rho: G_K \rar C$ by taking the $p$-th power so the $\mrm{Gal}(K/\Q)$-conjugate of $\psi'$ is, $(\psi')^p=(\psi')^{-1}(\psi')^{p+1}=(\psi')^{-1}$ since $p \equiv -1 \mod{3}$. }

\prop[result-for-ep]{Let $p \geq 11$ be a prime congruent to $2$ mod $9$ (resp. to $5$ mod $9$) and $E$ be a quadratic twist of the elliptic curve $y^2=x^3+p^2$ (resp. of $y^2=x^3+p$). Let $F/\Q$ be an elliptic curve such that $E[p](\Qbar)$ and $F[p](\Qbar)$ are isomorphic as $G_{\Q}$-modules. Then any prime $\ell$ of potentially multiplicative reduction for $F$ is exceptional for $p$ and congruent to $\pm 1 \mod{3p}$. In particular, if $p \leq 10000$, then $j(F) \in \Z$ and $E$ and $F$ have the same conductor. 
}

\demo{It is enough to check that Corollary \ref{fimm-invariants} applies. By Proposition \ref{ep-is-good} and Corollary \ref{ep-compute-g}, if $g$ is the unique newform in $\mathcal{S}_1(\Gamma_1(108))$, then it is enough to check that there is $f \in \mathcal{S}_2(\Gamma_0(p))$ such that $L(f,g,1) \neq 0$. By Propositions \ref{epn-supersing} and Corollary \ref{factor-hilbclasspoly}, it is enough to show that the Hilbert class polynomial $H_{-108}$ is not a cube modulo $p$ when $p \leq 108$. This is checked by a quick MAGMA calculation \cite{magma}.}

\cor[complete-result-for-ep]{In Proposition \ref{result-for-ep}, if one furthermore has $p \in \{23,29,41,47,59\}$, then $E$ and $F$ are isogenous.}

\demo{This follows from \cite[Theorem 1.3]{Cremona-Freitas}. Indeed, Tate's algorithm \cite[Chapter IV.9]{AEC2} shows that the conductor of the quadratic twist $y^2=x^3+\varepsilon p^{\alpha}$ of $E$ (where $\alpha =1$ if $p \in \{23,41,59\}$ and $\alpha=2$ otherwise, and $\varepsilon$ is a sign congruent to $p^{\alpha}$ modulo $4$) is exactly $108p^2 \leq 500000$ (we are stopping at Step 3 of \emph{loc.cit.} at $3$ -- i.e. Kodaira type II -- and at Step 5 at $2$ -- i.e. Kodaira type IV). }

\prop[result-for-en]{Let $p \geq 11$ be a prime and $\alpha \in \{1,2\}$ be an integer such that either $p \equiv 2 \mod{9}$ and $\alpha=2$, or $p \equiv 5 \mod{9}$ and $\alpha=1$. Let $n'$ be a cube-free integer prime to $p$, $n=n'p^{\alpha}$, and $E$ be a quadratic twist of the elliptic curve $E_n$ with equation $y^2=x^3+n$. Let $c=c'_6c_23^{e_3}$, where
\begin{itemize}[noitemsep,label=\tiny$\bullet$]
\item $c'_6$ is the product of the primes $\ell \nmid 6$ dividing $n'$,
\item $c_2 =1$ if $v_2(n')=1$ and $c_2=2$ otherwise,
\item $3^{1+2e_3}$ is the $3$-part of the conductor of $E_n$ if $n$ is not a cube in $\Q_3$, and $e_3=0$ otherwise.
\end{itemize}
Suppose that either $3c^2 \leq p$, or $c \leq 300$ and $(p,c)$ is distinct from $(11,271)$ and not contained in $\{11,23\} \times \{6 \cdot 29,9 \cdot 29,10 \cdot 29\}$. 

 Let $F/\Q$ be an elliptic curve such that $E[p](\Qbar)$ and $F[p](\Qbar)$ are isomorphic as $G_{\Q}$-modules. Then any prime $\ell$ of potentially multiplicative reduction for $F$ is exceptional for $p$ and congruent to $\pm 1 \mod{3p}$. In particular, if $p \leq 10000$, then $j(F) \in \Z$ and $E$ and $F$ have the same conductor.}

\demo{By Corollary \ref{fimm-invariants}, Proposition \ref{artin-conductor-rhoprime}, Lemma \ref{anticyc-chars}, and finally Proposition \ref{epn-caishutian}, it is enough to show that if $\OO \rar R$ is an optimal embedding of an order of discriminant $-3c^2$ into a maximal order of the quaternion algebra (over $\Q$) only ramified at $p\infty$, then the fibres of the map $\varphi$ of \emph{loc.cit.} do not equidistribute in any quotient of order $3$ of $\mrm{Pic}(\OO)$. By Corollary \ref{epn-supersing}, this holds if the Hilbert polynomial of discriminant $-3c^2$ is not a cube modulo $p$, which is the case in particular $p \geq 3c^2$ by Proposition \ref{injective-map}. We then perform an exhaustive search for the couples $(p,c)$ with $p \leq 3c^2 \leq 3 \cdot 300^2$ such the Hilbert class polynomial of discriminant $-3c^2$ is a cube modulo $p$, and check using MAGMA that the condition at the beginning of this proof is satisfied. }

\newpage
\nocite{*}
\bibliographystyle{amsplain}

{\scriptsize \bibliography{biblio-cartan-3}}

\appendix
\renewcommand{\thepropo}{\thesection.\arabic{propo}}

\section{Relative curves, relative Jacobians and functoriality}
\label{appendix:algeom}

The purpose of this appendix is to recall some results about the algebraic geometry of curves and their Jacobians. We need to extend some well-known constructions to the setting where our ``curves'' (contrary to what seems the standard convention) do not have geometrically connected fibres. 

We start with brief reminders on Picard groups.  

If $X$ is a scheme, $\mrm{Pic}(X)$ denotes the abelian group of isomorphism classes of line bundles over $X$. The rule $X \rar \mrm{Pic}(X)$ is a contravariant functor for morphisms of schemes: if $f: X \rar Y$ is a morphism of schemes, then $\mathcal{L} \mapsto f^{\ast}\mathcal{L}$ is a group homomorphism $\mrm{Pic}(Y) \rar \mrm{Pic}(X)$. In fact, $\mrm{Pic}$ is also a covariant functor for finite locally free maps: if $f: X \rar Y$ is a finite locally free morphism of schemes, the rule $f_{\ast}: \mathcal{L} \in \mrm{Pic}(X) \mapsto (\det{f_{\ast}\mathcal{L}})\otimes (\det{f_{\ast}\OO_X})^{\otimes (-1)} \in \mrm{Pic}(Y)$ is a group homomorphism, and $f \mapsto f_{\ast}$ is functorial (in the subcategory of $\Sch$ with the same objects, but where the arrows are exactly finite locally free morphisms). If $f: X \rar Y$ is finite locally free, the endomorphism $f_{\ast}f^{\ast}$ of $\mrm{Pic}(Y)$ is the $\deg{f}$-th power operator ($\deg{f}$ being a Zariski-locally constant function on $Y$). In particular, if $f$ is an isomorphism, $f_{\ast}$ and $f^{\ast}$ are inverses of each other. 

For any scheme $X$, $\mrm{Pic}(X)$ identifies to $H^1_{Zar}(X,\OO_X^{\times})$.\footnote{As pointed out in \cite[p. 284]{Conrad-Duality}, there is a choice of sign involved. We make one such choice once and for all.} The pull-back by a morphism $f: X \rar Y$ corresponds to $H^1(Y,\OO_Y^{\times}) \rar H^1(X,f^{-1}\OO_Y^{\times}) \rar H^1(X,\OO_X^{\times})$, where the first arrow is a general construction in sheaf cohomology, and the second arrow comes from the map $f^{-1}\OO_Y^{\times} \rar \OO_X^{\times}$ induced by $f$ itself. If $f: X \rar Y$ is finite locally free, the push-forward by $f$ in cohomology is constructed as follows. If $\mathcal{L}$ is a line bundle over $\OO_X$, then $f_{\ast}\mathcal{L}$ is a line bundle over the ringed space $(Y,f_{\ast}\OO_X)$, so this construction produces a map $f': H^1(X,\OO_X^{\times}) \rar H^1(Y,f_{\ast}\OO_X^{\times})$. Then the map in cohomology induced by $f_{\ast}$ is the composition of $f'$ with the norm map $f_{\ast}\OO_X^{\times} \rar \OO_Y^{\times}$. This has the following implication: 

\lem[mixed-p-functoriality]{Consider the following commutative diagram of schemes, where $g,g'$ are finite locally free. 
\[
\begin{tikzcd}[ampersand replacement=\&]
X' \arrow{r}{g'} \arrow{d}{f'}\& Y' \arrow{d}{f}\\
X \arrow{r}{g} \& Y
\end{tikzcd}
\]

Suppose that, for any open subscheme $U \subset Y$ and any $s \in \OO_X(g^{-1}(U))$, one has 
\[N_{\OO_{X'}((g \circ f')^{-1}(U))/\OO_{Y'}(f^{-1}(U))}((f')^{\sharp})(s) = f^{\sharp}(N_{\OO_X(f^{-1}(U))/\OO_Y(U)}s).\] 

Then the two maps $f^{\ast}g_{\ast}, (g')_{\ast}(f')^{\ast}: \mrm{Pic}(X) \rar \mrm{Pic}(Y')$ are equal. 
}

Let now $S$ be a scheme, and $X$ be a proper flat finitely presented $S$-scheme of relative dimension one. We call $P_{X/S}$ the functor $\Sch_S \rar \Ab$ sending a $S$-scheme $T$ to $\mrm{Pic}(X \times T)$, and $P^0_{X/S}$ its subfunctor sending a $S$-scheme $T$ on the subset of $\mrm{Pic}(X \times_S T)$ made only with those line bundles $\mathcal{L}$ such that, for every $t \in T$, the line bundle $\mathcal{L}_t$ on $X \times_S \Sp{\kappa(t)}$ has degree zero on every connected component, in the sense of \cite[Definition 0AYR]{Stacks}\footnote{Note that we apply the definition of \emph{loc.cit.} to every connected component of $X \times_S \Sp{\kappa(t)}$ separately.}. We denote by $\mrm{Pic}_{X/S}$ (resp. $\mrm{Pic}^0_{X/S}$) the sheaf over the big fppf site associated to $P_{X/S}$ (resp. to $P^0_{X/S}$). Because the quotient pre-sheaf $P_{X/S}/P^0_{X/S}$ is separated\footnote{Which amounts to the following claim: let $X$ be a proper scheme over some field $k$, pure of dimension one. Let $k'$ be a field extension and let $\mathcal{L}$ be a line bundle over $X$ such that the degree of $\mathcal{L}_{k'}$ over every connected component of $X_{k'}$ is zero. Then the degree of $\mathcal{L}$ over every connected component of $X$ is zero.}, the morphism of sheaves $\mrm{Pic}^0_{X/S} \rar \mrm{Pic}_{X/S}$ is injective. The formation of $\mrm{Pic}_{X/S}$ and $\mrm{Pic}^0_{X/S}$ commutes with base change of $S$. 

If $f: X \rar Y$ is a morphism of proper flat finitely presented $S$-schemes of relative dimension one, $f^{\ast}$ induces a morphism $P_{Y/S} \rar P_{X/S}$ (sending $\mathcal{L} \in \mrm{Pic}(Y \times_S T)=P_{Y/S}(T)$ to $(f \times_S \mrm{id})^{\ast}\mathcal{L} \in \mrm{Pic}(X \times T)=P_{X/S}(T)$) which maps $P^0_{Y/S}$ into $P^0_{X/S}$, so in particular $f^{\ast}$ induces a morphism $\mrm{Pic}^0_{Y/S} \rar \mrm{Pic}^0_{X/S}$ of fppf-sheaves. Its formation commutes with base change, and the rule $f \mapsto f^{\ast}$ is functorial contravariant.  

Suppose furthermore that $f$ is finite locally free, then the morphism $f_{\ast}: P_{X/S} \rar P_{Y/S}$ (sending $\mathcal{L} \in \mrm{Pic}(X \times_S T)=P_{X/S}(T)$ to $(f \times_S \mrm{id})_{\ast}\mathcal{L} \in \mrm{Pic}(Y \times T)=P_{Y/S}(T)$) sends the pre-sheaf $P^0_{X/S}$ into $P^0_{Y/S}$, so induces a map $\mrm{Pic}^0_{X/S} \rar \mrm{Pic}^0_{Y/S}$ of fppf-sheaves. The formation of this $f_{\ast}$ commutes with $f_{\ast}$ and the rule $f \mapsto f_{\ast}$ is functorial covariant.  

Lemma \ref{mixed-p-functoriality} then has the following consequence:

\lem[mixed-pic0-functoriality]{Let $S$ be a scheme and consider the following commutative diagram of proper flat finitely presented $S$-schemes of relative dimension one. 
\[
\begin{tikzcd}[ampersand replacement=\&]
X' \arrow{r}{g'} \arrow{d}{f'}\& Y' \arrow{d}{f}\\
X \arrow{r}{g} \& Y
\end{tikzcd}
\]

Assume that $f,f',g,g'$ are finite locally free, and that there exists an affine open immersion $\iota: U \rar Y$ such that
\begin{itemize}[noitemsep,label=\tiny$\bullet$]
\item the diagram is Cartesian above $U$,
\item for any $S$-scheme $T$, $\OO_{Y_T} \rar (\iota \times \mrm{id}_T)_{\ast}\OO_{U_T}$ is injective.
\end{itemize}
Then the two maps $f^{\ast}g_{\ast}, (g')_{\ast}(f')^{\ast}: P_{X/S} \rar P_{Y'/S}$ are equal, hence their sheafifications are equal as maps $\mrm{Pic}^i_{X/S} \rar \mrm{Pic}^i_{Y'/S}$, where $i$ is either empty or zero. 

 }

\demo{Let $U \rar Y$ be an affine open immersion of schemes such that $\OO_Y \rar \iota_{\ast}\OO_U$ is injective. It is elementary that for any flat $Y$-scheme $Y'$, the base change $\iota': U' \rar Y'$ of $\iota$ by $Y' \rar Y$ is an affine open immersion such that $\OO_{Y'} \rar \iota'_{\ast} \OO_{U'}$ is injective.

Therefore, for $Z \in \{X,X',Y',Y\}$, for any $S$-scheme $T$ and open subscheme $V \subset Z_T$, $\OO(V) \rar \OO(V \times_{Y_T} U_T)$ is injective. Since the diagram is Cartesian above $U$, this implies that, for every $S$-scheme $T$ and every affine open subscheme $V \subset Y_T$, the maps \[N_{\OO(X'_T \times_{Y_T} V)/\OO(Y' \times_{Y_T} V)} \circ (f'_T)^{\sharp},\, (f_T)^{\sharp} \circ N_{\OO(X_T \times_{Y_T} V)/\OO(V)}: \OO(X_T \times_{Y_T} V) \rar \OO(Y'_T \times_{Y_T} V)\] are equal. By Lemma \ref{mixed-p-functoriality} (applied to a suitable base change of the diagram), it follows that $f^{\ast}g_{\ast}, (g')_{\ast}(f')^{\ast}: P_{X/S} \rar P_{Y'/S}$ are equal, and we can then pass to the sheafification and restrict to the $\mrm{Pic}^0$ subsheaf.}

If $S$ is a scheme, a \emph{relative curve} over $S$ is a sequence $X \overset{g}{\rar} T \overset{\pi}{\rar} S$, where $T \rar S$ is finite \'etale, $X$ is a smooth projective $S$-scheme of relative dimension one, and $g_{\ast}\OO_X = \OO_T$ holds universally. 

\lem[loc-free-omega1-curve]{Let $X \rar T \rar S$ be a relative curve such that the fibres of $X \rar T$ have constant genus $g$ (e.g. if $T$ is connected) and let $f: X \rar S$ denote the structure map. Then $f_{\ast}\Omega^1_{X/S}$ is a finite locally free $f_{\ast}\OO_X$-module of rank $g$, hence a locally free $\OO_S$-module, and the formation of $f_{\ast}\Omega^1_{X/S}$ as well as $f_{\ast}\OO_X$ commutes with base change. 
}

\demo{We may assume that $S$ is affine. Since $T \rar S$ is finite \'etale, it is finite locally free and the obvious map $\Omega^1_{X/S} \rar \Omega^1_{X/T}$ is an isomorphism: therefore, it is enough to show the result when $T=S$, so that $\OO_S \rar f_{\ast}\OO_X$ is by definition a universal isomorphism. By assumption, the formation of $f_{\ast}\OO_X$ commutes with base change; furthermore, by \cite[Lemmas 0EX3, 0E7D]{Stacks} for every $i \geq 2$, $R^if_{\ast}\OO_X=0$ holds universally and the formation of $R^1f_{\ast}\OO_X$ commutes with base change. Then by cohomology and base change \cite[Lemma 5.1.1]{Conrad-Duality}, $R^1f_{\ast}\OO_X$ is locally free, and the conclusion follows from Grothendieck--Serre duality \cite[Theorem 5.1.2]{Conrad-Duality}.  }

\lem[loc-free-omega1-jac]{Let $X \rar S$ a smooth proper scheme of relative dimension one with geometrically connected fibres and $\alpha: J \rar S$ be its relative Jacobian, which represents $\mrm{Pic}^0_{X/S}$ (this fppf sheaf thus satisfies the sheaf property in the fpqc topology). There is a canonical difference map $\delta: X \times_S X \rar J$ sending any pair $(x,y)$ of $U$-points of $X$ (where $U$ is a $S$-scheme) to the class of the line bundle of $X_U$ attached to $[x]-[y]$. One has a canonical isomorphism $\iota: \alpha_{\ast}\Omega^1_{J/S} \rar f_{\ast}\Omega^1_{X/S}$ (whose formation commutes with base change) satisfying the following property: for any $S$-scheme $U$ and any $U$-point $x: U \rar X$, for any $\omega \in H^0(J_U,\Omega^1_{J_U/U})$, $\iota(\omega) = \delta(\cdot,x)^{\ast}\omega$.
}

\demo{Note that $J$ exists and its formation commutes to base change by \cite[9.4/4]{BLR}. First, assume that $X$ has a $S$-point $s_0$: then $j_{s_0}: x \mapsto [x]-[s_0]$ defines a $S$-morphism $X \rar J$. Then $\delta_{s_0}: (x,y) \mapsto j_{s_0}(x)-j_{s_0}(y)$ defines a map $X \times_S X \rar J$ which clearly does not depend on the choice of $s_0$. Pulling-back differentials along $j_{s_0}$ defines a map $\iota_{s_0}: \alpha_{\ast}\Omega^1_{J/S} \rar f_{\ast}\Omega^1_{X/S}$: since $\alpha_{\ast}\Omega^1_{J/S} \simeq 0^{\ast} \Omega^1_{J/S}$ by \cite[Proposition 1.1 (b)]{LLR}, $\iota_{s_0}$ is a morphism of locally free $\OO_S$-modules of rank $g$ whose formation commutes with base change. By \cite[Proposition 2.2]{MilJac}, $\iota_{s_0}$ is fibrewise surjective, hence it is an isomorphism. 

We finally want to show that $\iota_{s_0}$ does not depend on the choice of $s_0$: by unpacking the definitions, it is enough to show that the global differentials on $J/S$ are translation-invariant. This is a consequence of the following claim: let $m,p_1,p_2: J \times_S J \rar J$ be the product and the two projections respectively, then $m^{\ast}=p_1^{\ast}+p_2^{\ast}$ as morphisms $\alpha_{\ast}\Omega^1_{J/S} \rar (\alpha \times \alpha)_{\ast} \Omega^1_{J \times J/S}$. This statement holds by \cite[Proposition 1.1 (b)]{LLR}, since it is dual to the (clear) claim that the two maps $m, p_1+p_2: \mrm{Lie}(J \times J) \rar \mrm{Lie}(J) \times \mrm{Lie}(J)$ (of fpqc-sheaves over $\Sch_S$) agree. 

In general, $X \rar S$ has sections \'etale-locally over $S$ by \cite[Corollaire 17.16.3]{EGA-IV4}, and the constructions above are canonical, so they glue in the \'etale topology, so $\delta$ and $\iota$ exist by descent.}

\prop[generalized-jac-exists]{Let $X$ be a relative curve over a scheme $S$. Then $\mrm{Pic}^0_{X/S}$ is a sheaf in the fpqc topology and it is represented by an abelian scheme $J$ over $S$, the \emph{relative Jacobian} of $X \rar S$. Moreover, there is a natural difference map $X \times_T X \rar J$ sending any pair $(x,y)$ of $U$-points of $X$ (where $U$ is a $S$-scheme) with the same image in $T$ to the class of the line bundle of $X_U$ attached to $[x]-[y]$. 

If $\alpha: J \rar S$ and $f: X \rar S$ are the structure maps, the $\OO_S$-modules $f_{\ast}\Omega^1_{X/S}$, $\alpha_{\ast}\Omega^1_{J/S}$ are locally free, their formation commutes with base change, and one has a canonical isomorphism $\iota: \alpha_{\ast}\Omega^1_{J/S} \rar f_{\ast}\Omega^1_{X/S}$ (whose formation commutes with base change) satisfying the following property: for any $S$-scheme $U$ and any $S$-point $x: U \rar X$ (so that $g(x)$ endows $U$ into a $T$-scheme structure), for any $\omega \in H^0(J_U,\Omega^1_{J_U/U})$, $\iota(\omega)_{|X \times_{T,g(x)} U} = \delta(\cdot,x)^{\ast}\omega$.}

\demo{Apart from the representability of $\mrm{Pic}^0_{X/S}$, every property is \'etale-local (because $\iota$ is canonical). So by \cite[Lemma 04HN]{Stacks}, granting temporarily the representability of $\mrm{Pic}^0_{X/S}$ in general, we may assume that $T \rar S$ is the disjoint reunion of copies $T_i$ of $S$. Let $X_i=X \times_T T_i$, then $\mrm{Pic}^0_{X/S}$ is clearly isomorphic to the direct sum of the $\mrm{Pic}^0_{X_i/S}$, and $f_i: X_i \rar S$ is smooth proper of relative dimension one such that $\OO_S \rar (f_i)_{\ast}\OO_{X_i}$ is universally an isomorphism. We can then apply Lemmas \ref{loc-free-omega1-curve} and \ref{loc-free-omega1-jac}. 

Hence only the representability of $\mrm{Pic}^0_{X/S}$ remains to show. This is Zariski-local over $S$, so it is enough to consider the case where $S$ is affine. By classical approximation results \cite[Th\'eor\`eme 8.8.2, Th\'eor\`eme 8.10.5, Corollaire 11.2.6.1, Th\'eor\`eme 12.2.4]{EGA-IV3} and \cite[Proposition 4.5.6, Corollaire 6.5.2]{EGA-IV2}, we may assume that $S$ is furthermore Noetherian. 

Note that $P_{X/S}$ is exactly the Weil restriction along the finite \'etale map $T \rar S$ of $P_{X/T}$. Moreover, if $U$ is a $S$-scheme and $u \in U$, then $X \times_S \Sp{\kappa(u)}$ is the disjoint reunion of the connected smooth proper one-dimensional schemes $X \times_T \Sp{\kappa(u')}$, where $u' \in U \times_S T$ runs over points above $u$: this means that $P^0_{X/S}$ is the Weil restriction of $P^0_{X/T}$ along the finite \'etale map $T \rar S$. We are then done if we can show:
\begin{itemize}[noitemsep,label=\tiny$\bullet$]
\item The formation of the Weil restriction of pre-sheaves (on the big fppf site) along a finite \'etale map commutes with arbitrary base change,
\item The formation of the Weil restriction of a pre-sheaf along a finite \'etale map commutes with fppf-sheafification,
\item The Weil restriction of an abelian scheme along a finite \'etale map exists and is an abelian scheme. 
\end{itemize}  

The first bullet point is formal (cf. the adjunction formula \cite[7.6/1]{BLR}). The Weil restriction of an abelian scheme $A \rar T$ along a finite \'etale map $T \rar S$ is a scheme by \cite[7.6/4]{BLR}, it is smooth proper by \cite[7.6/5]{BLR}, is a group scheme by functoriality of the Weil restriction (see the discussion of group functors on p. 192 of \emph{op.cit.}). To check that it is an abelian scheme, we need to check that its geometric fibres are connected, so we are reduced to the situation where $S$ is an algebraically closed field and $T$ is a disjoint reunion of copies of $S$, and the claim is clear in this case, since the Weil restriction is a product of copies of the abelian varieties. 

The second bullet point is a formal consequence of the following claim: let $T \rar S$ be a finite \'etale map, $U$ be a $S$-scheme, and $Z$ be a flat $U \times_S T$-scheme of finite presentation. Then $Z$ has a fppf-cover by $U \times_S T$-schemes of the form $V \times_S T$, where $V$ is a flat finitely presented $U$-scheme (and the structure map to $T$ is the second projection). To do this, it is enough to consider the case where $U=S$ and $Z$ are both affine, and $\OO(T)$ is a free $\OO(S)$-module of finite rank. Then $Z$ has a Weil restriction $Z_1$ along the map $T \rar S$, which is affine (this is the first part of the proof of \cite[7.6/4]{BLR}), flat of finite presentation over $S$ by \cite[7.6/5]{BLR}, and we wish to show that the canonical map $Z_1 \times_S T \rar Z$ (see the first paragraph of p. 192 of \emph{op.cit.}) is faithfully flat and finitely presented. Since Weil restriction commutes with base change, we can work \'etale-locally over $S$ and hence assume that $T$ is the disjoint reunion of a finite number of copies of $S$, and the result is trivial in this case. 
}

We discuss how the push-forward and pull-back functorialities behave with respect to differentials. 

\prop[omega1-functorialities]{Let $\alpha: X \rar S, \beta: Y \rar S$ be two relative curves with relative Jacobians $\alpha_J: J_X \rar S, \beta_J: J_Y \rar S$, and $f: X \rar Y$ be a finite locally free morphism. If the factorizations are $X \rar T \rar S$ and $Y \rar T' \rar S$, then $f$ induces a $S$-morphism $f_S: T \rar T'$. Then the following diagrams commute:
\[
\begin{tikzcd}[ampersand replacement=\&]
X \times_T X \arrow{d}{f \times f} \arrow{r}{\delta_X} \& J_X \arrow{d}{f_{\ast}} \& (\beta_J)_{\ast}\Omega^1_{J_Y/S} \arrow{r}{\iota_Y} \arrow{d}{\Omega^1(f_{\ast})} \& (\beta)_{\ast}\Omega^1_{Y/S} \arrow{d}{\Omega^1(f)}\\
Y \times_{T'} Y \arrow{r}{\delta_Y} \& J_Y \& (\alpha_J)_{\ast}\Omega^1_{J_X/S} \arrow{r}{\iota_X} \& \alpha_{\ast}\Omega^1_{X/S}
\end{tikzcd}
\]

Let $\mrm{tr}: f_{\ast}\Omega^1_{X/S} \rar \Omega^1_{Y/S}$ be the canonical trace map, whose formation commutes to base change, and whose restriction to $f_{\ast}f^{\ast}\Omega^1_{Y/S} \simeq f_{\ast}\OO_X \otimes \Omega^1_{Y/S}$ is given by the classical trace $f_{\ast}\OO_X \rar \OO_Y$. Then the following diagram commutes:

\[
\begin{tikzcd}[ampersand replacement=\&]
(\alpha_J)_{\ast}\Omega^1_{J_X/S} \arrow{d}{\Omega^1(f^{\ast})} \arrow{r}{\iota_X} \& \alpha_{\ast}\Omega^1_{X/S} \arrow{d}{\beta_{\ast}\mrm{tr}}\\
(\beta_J)_{\ast}\Omega^1_{J_Y/S} \arrow{r}{\iota_Y} \& \beta_{\ast}\Omega^1_{Y/S} 
\end{tikzcd}
\]

}

\demo{If $X \rar T \rar S$ is a relative curve, then it is the Stein factorization of $X \rar S$, and in particular $T$ is the relative normalization of $S$ in $X$ by \cite[Theorem 03H2]{Stacks}. This is why $f$ induces a $S$-morphism $T \rar T'$. The first diagram commutes by construction of the difference map and by the explicit description of the push-forward for a line bundle attached to a section. That the second diagram commutes can be checked \'etale-locally over $S$, so that we may assume by \cite[Lemma 04HN]{Stacks} and \cite[Corollaire 17.16.3]{EGA-IV4} that $T,T'$ are both disjoint reunions of copies of $S$, and every $X \times_T S$, $Y \times_{T'} S$ admits a $S$-point, in which case the conclusion follows from the definitions of $\iota_X,\iota_Y$. 

We now discuss the commutativity of the last diagram. The existence and properties of $\mrm{tr}$ are given by Grothendieck duality by \cite[(2.7.36), (2.7.41)]{Conrad-Duality}. 

By working \'etale-locally over $S$ and breaking up $X,Y$ into pieces that are smooth $S$-schemes with geometrically connected fibres, we may assume as previously that $T \rar S, T' \rar S$ are isomorphisms and $X,Y$ both have $S$-points. The $\OO_S$-dual of the diagram, under \cite[Proposition 1.1 (b)]{LLR} on the left, and Grothendieck duality on the right, is exactly (by \cite[Theorem B.4.1]{Conrad-Duality}) the diagram of \cite[Proposition 1.3 (c)]{LLR}, and the conclusion follows. 
}

\end{document}